\documentclass[reqno]{amsart}
\numberwithin{equation}{section}
\usepackage{graphicx}

\newtheorem{corollary}{\sc Corollary}[section]
\newtheorem{theorem}{\sc Theorem}[section]
\newtheorem{lemma}{\sc Lemma}[section]

\newtheorem{remark}{\sc Remark}
\newtheorem{definition}{\sc Definition}[section]

\def\bdy #1{\partial #1}

\def\bbR{{\mathbb R}}

\def\div{\operatorname{div}}
\def\Def{\operatorname{Def}}
\def\id{{\text{Id}}}
\def\supp{{\text{supp}}}
\def\eps1{{\epsilon_1}}
\def\L{{\mathcal L}}
\def\bL{\bar{\mathcal L}}
\def\G{{\mathcal M}}
\def\bG{\bar{\mathcal M}}
\def\th{{\tilde{h}}}
\def\tth{{\bar{h}}}
\def\tA{{\tilde{A}}}
\def\ttA{{\bar{A}}}

\def\ta{{\tilde{a}}}
\def\tta{{\bar{a}}}
\def\tu{{\tilde{u}}}
\def\tv{{\tilde{v}}}
\def\tw{{\tilde{w}}}
\def\ttv{{\bar{v}}}
\def\tq{{\tilde{q}}}

\def\tF{{\tilde{F}}}

\def\teta{{\tilde{\eta}}}
\def\tteta{{\bar{\eta}}}

\def\tTheta{{\tilde{\Theta}}}
\def\ttTheta{{\bar{\Theta}}}
\def\tmem{{\mathfrak t}_{\text{mem}}}

\def\mv{{\mathfrak v}}
\def\mh{{\mathfrak h}}
\def\mq{{\mathfrak q}}
\def\nv{{v_{\kappa}}}
\def\nh{{h_{\kappa}}}
\def\nq{{q_{\kappa}}}
\def\ov{{v_{\epsilon}}}
\def\oh{{h_{\epsilon}}}
\def\oq{{q_{\epsilon}}}
\def\H1H2{{H^{1;2}(\Omega;\Gamma)}}

\begin{document}
\title[Navier-Stokes interacting with a fluid shell]
{Navier-Stokes equations interacting with a nonlinear elastic fluid
shell}

\author{C.H. Arthur Cheng}
\email{cchsiao@math.ucdavis.edu}
\author{Daniel Coutand}
\email{coutand@math.ucdavis.edu}
\author{Steve Shkoller}
\email{shkoller@math.ucdavis.edu}

\address{Department of Mathematics, University of California, Davis, CA 95616}

\maketitle
\begin{abstract}
We study a moving boundary value problem consisting of a viscous
incompressible fluid moving and interacting with a nonlinear elastic
fluid shell.  The fluid motion is governed by the Navier-Stokes
equations, while the fluid shell is modeled by a bending energy
which extremizes the Willmore functional and a membrane energy that
extremizes the surface area of the shell. The fluid flow and shell
deformation are coupled together by continuity of displacements and
tractions (stresses) along the moving material interface.  We prove
existence and uniqueness of solutions in Sobolev spaces.
\end{abstract}

%In our companian paper \ref{ChCoSh2006}, we list a number of references for
%analytical studies of fluid-structure interaction problems.
\section{Introduction}
\subsection{The problem statement and background}

We are concerned with establishing the existence and uniqueness of
solutions to the time-dependent incompressible Navier-Stokes
equations interacting with a nonlinear elastic fluid shell
(bio-membrane). Recently, there have been many experimental and
analytic studies on the configurations and deformations of elastic
bio-membranes (see, for example, \cite{RSCh2002}, \cite{YCFu1981},
\cite{WMGe1995}, \cite{MiSeWoDo1994}, \cite{OuHe1989},
\cite{USe1991}, and \cite{USe1993}), but the basic analysis of the
coupled fluid-structure interaction remains open.  The fundamental
difficulties arise from the degenerate elliptic operators that arise
as the shell tractions.  As we detail below, the bending energy of
the shell is the well-known Willmore function, the integral over the
moving surface of the square of the mean curvature.  The degenerate
elliptic operator arising from the first variation of this
functional is a fourth order nonlinear operator that smoothes only
in the direction which is normal to the moving domain.  Our analysis
will provide a well-posedness theorem and explain the interesting
interaction between the shape of the shell and the flow of the
fluid.

Fluid-structure interaction problems involving moving material
interfaces have been the focus of active research since the
nineties. The first problem solved in this area was for the case of
a rigid body moving in a viscous fluid (see \cite{DeEs}, \cite{GrMa}
and also the early works of \cite{Wein1972} and \cite{Serre} for a
rigid body moving in a Stokes flow in the full space). The case of
an elastic body moving in a viscous fluid was considerably more
challenging because of some apparent regularity incompatibilities
between the parabolic fluid phase and the hyperbolic solid phase.
The first existence results in this area were for regularized
elasticity laws, such as in \cite{DeEsGrTa} for a {\it finite}
number of elastic modes, or in \cite{DaVe}, \cite{ChDeEsGr}, and
\cite{Bo2005} for hyperviscous elasticity laws, or in \cite{LiWa} in
which a phase-field regularization ``fattens'' the sharp interface
via a diffuse-interface model.

The treatment of classical elasticity laws for the solid phase,
without any regularizing term, was only considered recently in
\cite{CoSh2005} for the three-dimensional linear St.
Venant-Kirchhoff constitutive law and in \cite{CoSh2006} for
quasilinear elastodynamics coupled to the Navier-Stokes equations.
Some of the basic new ideas introduced in those works concerned a
functional framework that scales in a hyperbolic fashion (and is
therefore driven by the solid phase), the introduction of
approximate problems either penalized with respect to the
divergence-free constraint in the moving fluid domain, or smoothed
by an appropriate parabolic artificial viscosity in the solid phase
(chosen in an asymptotically convergent and consistent fashion), and
the tracking of the motion of the interface by difference quotients
techniques.

In our companion paper \cite{ChCoSh2006}, we study the interaction
of the Navier-Stokes equations with an elastic solid shell. Herein,
we treat the case of a fluid shell or bio-membrane. This is a moving
boundary problem that models the motion of a viscous incompressible
Newtonian fluid inside of a deformable elastic fluid structure.

Let $\Omega\subset\bbR^3$ denote an open bounded domain with
boundary $\Gamma := \bdy\Omega$. For each $t\in (0,T]$, we wish to
find the domain $\Omega(t)$, a divergence-free velocity field
$u(t,\cdot)$, a pressure function $p(t,\cdot)$ on $\Omega(t)$, and a
volume-preserving transformation $\eta(t,\cdot):\Omega \to \bbR^3$
such that
\begin{subequations}\label{NSequation}
\begin{alignat}{2}
\Omega(t) &= \eta(t,\Omega) \,, && \\
\eta_t(t,x) &= u(t,\eta(t,x)) \,, && \\
u_t + \nabla_u u - \nu\Delta u &= -\nabla p + f &&\qquad \text{in }\
\Omega(t)\,,
\label{NSequation.a} \\
\div u &= 0 &&\qquad \text{in }\ \Omega(t)\,,
\label{NSequation.b}\\
(\nu\Def u -p\id)n &= {\mathfrak t}_{shell}  &&\qquad \text{on } \
\Gamma(t)\,,
\label{NSequation.c} \\
u(0,x) &= u_0(x) &&\qquad \forall x\in \Omega \,, \\
\eta(0,x) &= x &&\qquad \forall x\in \Omega \,,
\end{alignat}
\end{subequations}
where $\nu$ is the kinematic viscosity, $n(t,\cdot)$ is the outward pointing unit normal to $\Gamma(t)$,
$\Gamma(t) := \partial \Omega(t)$ denotes the boundary of $\Omega(t)$,
$\Def u$ is twice the rate of deformation tensor of $u$, given in coordinates by
$u^i_{,j} + u^j_{,i}$, and ${\mathfrak t}_{shell}$ is the traction imparted
onto the fluid by the elastic shell, which we describe next.

We shall consider a thin elastic shell modeled by the nonlinear Saint
Venant-Kirchhoff constitutive law. With $\epsilon$ denoting the
thickness of the shell, the hyperelastic stored energy function has the
asymptotic expansion
$$E_{shell} = \epsilon E_{mem} + \epsilon^3 E_{ben} + {\mathcal O}(\epsilon^4).$$
The membrane energy satisfies
\begin{align}
E_{mem} = \gamma \int_{\Gamma(t)} dS = \text{$\gamma$ times the
surface area of
$\Gamma(t)$}%\Big[\frac{\mu}{4} \sum_{\alpha,\beta=1}^2 (g_{\alpha\beta} -
%g_{0\alpha\beta})^2
%+ \frac{\mu\lambda}{4(2\mu + \lambda)}\Big(\sum_{\alpha=1}^2 (g_{\alpha\alpha} - g_{0\alpha\alpha})\Big)^2\Big] dS
\label{membrane}
\end{align}
where $\gamma>0$ is the surface tension, while the bending energy
$E_{ben}$ is given by
\begin{align}
E_{ben} = \int_{\Gamma(t)}\Big[(4\mu + 2\lambda) H^2 - 2\mu K\Big] dS,
\label{bending}
\end{align}
where $H$, $K$ denote the mean and Gauss curvatures on $\Gamma(t)$,
respectively, and $\lambda/2$ and $\mu/2$ are the
Lam$\acute{\text{e}}$ constants (see, for example,
\cite{GeKrMa1996}).

The traction vector
$${\mathfrak t}_{shell} = \epsilon {\mathfrak t}_{mem} + \epsilon^3 {\mathfrak t}_{ben} + {\mathcal O}(\epsilon^4)$$
is computed from the first variation of the energy function
$E_{shell}$; the traction vector associated to the membrane energy
is well-known to be
\begin{align}
{\mathfrak t}_{mem} = \gamma H n \,,
\end{align}
while the traction associated to the bending energy has a simple
intrinsic characterization given by
\begin{align}
{\mathfrak t}_{ben} = \sigma ({\Delta}_g H - 2HK + 2H^3)n \,,
\label{bendedstress}
\end{align}
where $\sigma$ is a function of the Lam$\acute{\text{e}}$ constants and
$\Delta_g$ denotes the Laplacian with respect to the induced metric $g$ on
$\Gamma(t)$:
\begin{align*}
\Delta_g f = \frac{1}{\sqrt{\det(g)}}\frac{\partial}{\partial x^\alpha}\Big(\sqrt{\det(g)} g^{\alpha\beta}
\frac{\partial f}{\partial x^\beta}\Big) \,.
\end{align*}

\subsection{Outline of the paper}
In Section \ref{Lagrangianformulation}, in addition to the use of
Lagrangian variables, we introduce a new coordinate system near the
boundary (shell) and three new maps, $\eta^\nu$, $\eta^\tau$, and
$h$, which facilitate the computation of the membrane and bending
tractions $\tmem$ and ${\mathfrak t}_{\text{ben}}$. A key
observation is the symmetry relation (\ref{gHsymmetry1}) which
reduces the derivative count on the tangential reparameterization
map $\eta^\tau$.

The space of solutions is introduced in Section \ref{notation}, and
the main theorem is stated in Section \ref{mainthm}. Section
\ref{convexset} defines our notation, and Section
\ref{preliminaryresult} provides some useful lemmas.

In Section \ref{linearizedprob}, we introduce the linearized and
regularized problems.  The regularization requires smoothing certain
variables as well as the introduction of a certain artificial viscosity
term on the boundary of the fluid domain.  Weak solutions of this
linear problem are established via a penalization scheme which approximates
the incompressibility of the fluid.

In Section \ref{energyestimatesection1},
we establish a regularity theory for our weak solution using energy estimates
for the mollified problem with constants that depend on the mollification
parameters.
In Section \ref{energyestimatesection2}, we improve these estimates so that
the constants are
independent of the artificial viscosity as well as other
regularization parameters.
This requires an elliptic estimate, arising from
the boundary condition (\ref{NSequation}e), which provides additional
regularity for the shape of the boundary.

In Section \ref{fixedpointuniqueness}, the Tychonoff fixed-point
theorem is used to prove the existence of solutions to the original
nonlinear problem (\ref{NSequation}) Uniqueness, following required
compatibility conditions, is established in Sections \ref{mainthm}
and Section \ref{fixedpointuniqueness}.

The inclusion of the inertial term $\epsilon \eta_{tt}$ into the
membrane traction $\tmem$ will be studied in a future publication.

\section{Lagrangian formulation}\label{Lagrangianformulation}
\subsection{A new coordinate system near the shell}
Consider the isometric immersion $\eta_0:(\Gamma,g_0) \to
(\bbR^3,\id)$. Let ${\mathcal B}=\Gamma\times(-\epsilon,\epsilon)$
where $\epsilon$ is chosen sufficiently small so that the map
$$B:{\mathcal B}\to\bbR^3: (y,z) \mapsto y+zN(y)$$
is itself an immersion, defining a tubular neighborhood of $\Gamma$
in $\bbR^3$. We can choose a coordinate system
$\frac{\partial}{\partial y^\alpha}$, $\alpha=1,2$ and
$\frac{\partial}{\partial z}$ on ${\mathcal B}$ where
$\frac{\partial}{\partial y^\alpha}$ denotes the tangential
derivative and $\frac{\partial}{\partial z}$ denotes the normal
derivative.

Let $G=B^*(\id)$ denote the induced metric on ${\mathcal B}$ from $\bbR^3$ so that
$$G(y,z)=G_z(y)+dz\otimes dz,$$
where $G_z$ is the metric on the surface $\Gamma\times\{z\}$; note
that $G_0=g_0$.

\begin{remark} By assumption,
${g_{0\alpha\beta} = \frac{\partial}{\partial y^\alpha} \cdot
\frac{\partial}{\partial y^\beta}}$, where $\cdot$ denotes the usual
Cartesian inner-product on $\bbR^n$. Let $C_{\alpha\beta}$ denote
the covariant components of the second fundamental form of the base
manifold $\Gamma$, so that $C_{\alpha\beta} = - N_{,\alpha}\cdot
\frac{\partial}{\partial y^\beta}$. Then, $G_z$ is given by
\begin{align*}
(G_z)_{\alpha\beta} = (g_0)_{\alpha\beta} - 2z C_{\alpha\beta} + z^2 g_0^{\gamma\delta} C_{\alpha\gamma}C_{\beta\delta}.
\end{align*}
\end{remark}
\vspace{.1 in}

Let $h:\Gamma\to(-\epsilon,\epsilon)$ be a smooth height function
and consider the graph of $h$ in ${\mathcal B}$, parameterized by
$\phi:\Gamma\to{\mathcal B}:y\mapsto (y,h(y))$. The tangent space to
graph($h$), considered as a submanifold of ${\mathcal B}$, is
spanned at a point $\phi(x)$ by the vectors
$$\phi_*(\frac{\partial}{\partial y^\alpha}) = \frac{\partial\phi}{\partial y^\alpha}=\frac{\partial}{\partial y^\alpha} +
\frac{\partial h}{\partial y^\alpha}\frac{\partial}{\partial z},$$
and the normal to graph($h$) is given by
\begin{align}
n(y)=J_h^{-1}(y)\Big(-G^{\alpha\beta}_{h(y)}\frac{\partial h}{\partial y^\alpha}\frac{\partial}{\partial y^\beta} +
\frac{\partial}{\partial z}\Big)\label{normal}
\end{align}
where $J_h=(1+h_{,\alpha}G^{\alpha\beta}_{h(y)} h_{,\beta})^{1/2}$. The mean curvature $H$ of graph($h$) is defined to be the
trace of $\nabla n$ where
\begin{align*}
(\nabla n)_{ij}=G(\nabla^{\mathcal B}_{\frac{\partial}{\partial w^i}} n,\frac{\partial}{\partial w^j}) \qquad\text{for
$i,j=1,2,3$}
\end{align*}
where $\frac{\partial}{\partial w^\alpha} = \frac{\partial}{\partial y^\alpha}$ for $\alpha=1$, $2$ and $\frac{\partial}{\partial
w^3}=\frac{\partial}{\partial z}$, and $\nabla^{\mathcal B}$ denotes the covariant derivative. Using (\ref{normal}),
\begin{align*}
(\nabla n)_{\alpha\beta} =&\ G\Big(\nabla^{\mathcal B}_{\frac{\partial}{\partial y^\alpha}}\Big[-J_h^{-1}
G_h^{\gamma\delta}h_{,\gamma}\frac{\partial}{\partial y^\delta}+J_h^{-1}\frac{\partial}{\partial z}\Big],\frac{\partial}
{\partial y^\beta}\Big)\\
=&\ -(G_h)_{\delta\beta}\Big[(J_h^{-1}G_h^{\gamma\delta}h_{,\gamma})_{,\alpha}
+ J_h^{-1}(-G_h^{\gamma\sigma}h_{,\gamma}\Gamma_{\alpha\sigma}^\delta
+ \Gamma_{\alpha 3}^\delta)\Big]\,; \\
(\nabla n)_{33} =&\ G\Big(\nabla^{\mathcal B}_{\frac{\partial}{\partial z}}\Big[-J_h^{-1}
G_h^{\gamma\delta}h_{,\gamma}\frac{\partial}{\partial y^\delta}+J_h^{-1}\frac{\partial}{\partial z}\Big],\frac{\partial}
{\partial z}\Big) \\
=&\ J_h^{-1}(-G_h^{\gamma\delta}h_{,\gamma}\Gamma_{3\delta}^3 + \Gamma_{33}^3)\,,
\end{align*}
where $\Gamma_{ij}^k$ denotes the Christoffel symbols with respect to the metric $G$.
It follows that the curvature of graph($h$) (in the divergence form) is
\begin{align}
H = -(J_h^{-1}G_h^{\gamma\delta}h_{,\gamma})_{,\delta} + J_h^{-1}
(-G_h^{\gamma\delta}h_{,\gamma}\Gamma_{j\delta}^j + \Gamma_{j 3}^j),
\label{divform}
\end{align}
or (in the quasilinear form)
\begin{align}
H = -J_h^{-1}G_h^{\alpha\beta}\Big[\delta_{\beta\gamma} - J_h^{-2} G^{\gamma\delta}_h h_{,\beta}h_{,\delta}\Big]h_{,\alpha\gamma}
+G_h^{\alpha\beta}F_{\alpha\beta}(y,h,\nabla h),\label{quasilinearform}
\end{align}
where $F_{\alpha\beta}$ denotes a smooth generic function of $y$, $h$ and $\nabla h$.

\begin{remark} Note that $G_h$ denotes the metric
$G_{z=h(y)}$, and not the metric on the submanifold graph$(h)$.
\end{remark}

\begin{remark}
If the initial height function is zero, i.e., $h(0)=0$, then
$H(0)=\Gamma^j_{j3}(0)$, which is the mean curvature of the base
manifold $\Gamma$ as required.
\end{remark}

%\begin{example} In the case that the initial surface $\Gamma$ is flat, so that
% $\Gamma \subset \bbR^2\times\{0\}$, then $G_{\alpha\beta} =
%\delta_{\alpha\beta}$, and
%$$H = -\frac{1}{\sqrt{1+|\nabla h|^2}}\Big[\delta_{\alpha\beta} - \frac{h_{,\alpha}h_{,\beta}}{1+|\nabla h|^2}\Big]h_{,\alpha\beta}.$$
%\end{example}
%
%\begin{example} In the case that the initial surface  $\Gamma=S^2$,
%with local coordinates $(y_1,y_2) = (\theta,\phi)$, then the metric
%$$G = (1+z)^2 \sin^2\phi d\theta^2 + (1+z)^2 d\phi^2 + dz^2\,,$$
%and the mean curvature is given by the formula
%\begin{align*}
%H=-\Big[\Big(J_h^{-1}\frac{h_\theta}{(1+h)^2\sin^2\phi} \Big)_{,\theta} +
%\Big(J_h^{-1}\frac{h_\phi}{(1+h)^2}\Big)_{,\phi}\Big] + J_h^{-1}\Big[G_h^{\gamma 2}h_{,\gamma}\cot\phi + \frac{2}{1+h}\Big].
%\end{align*}
%\end{example}
%\vspace{.1 in}

\subsection{Tangential reparameterization symmetry}
Let ${\mathcal N}$ denote the normal bundle to $\Gamma$, so that for
each $y \in \Gamma$, we have the Whitney sum ${\mathbb R}^3= T_y
\Gamma \oplus {\mathcal N}_y$.

Given a signed height function $h: \Gamma \times [0,T) \rightarrow
{\mathbb R}$, for each $t\in [0,T)$, define the {\it normal} map
$$
\eta^\nu : \Gamma \times [0,T) \rightarrow \Gamma(t), \ \ \ (y,t)
\mapsto y+ h(y,t) N(y), \ \ \ N(y) \in {\mathcal N}_y\,.
$$
Then, there exists a unique {\it tangential} map $\eta^\tau: \Gamma
\times [0,T)\rightarrow \Gamma $ (a diffeomorphism as long as $h$
remains a graph) such that the diffeomorphism $\eta(t)$ has the
decomposition
$$\eta(\cdot, t) = \eta^\nu(\cdot,t) \circ \eta^\tau(\cdot, t), \ \ \
\eta(y,t) = \eta^\tau(y,t) + h(\eta^\tau(y,t),t) N(\eta^\tau(y,t))\,.$$
\begin{figure}
\begin{center}
\includegraphics[scale = 0.45]{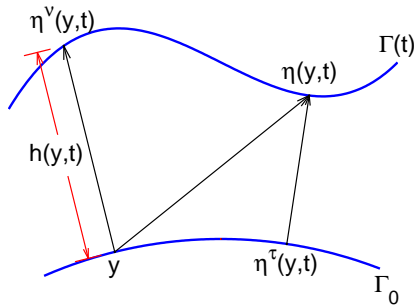}
\caption{The maps
$\eta^\tau$ and $\eta^\nu$}
\end{center}
\end{figure}
The tangent vector $\eta,_\alpha$ to $\Gamma(t)$ can be decomposed with
respect to the Whitney sum as
$
\eta_{,\alpha}(y,t) = \eta^\kappa_{,\alpha}(y,t)\frac{\partial}
{\partial y^\kappa} + h_{,\kappa}(\eta^\tau(y,t),t)
\eta^\kappa_{,\alpha}\frac{\partial}{\partial z}$
and hence the induced metric $g_{\alpha\beta} = \eta,_\alpha \cdot \eta,_\beta$
may be expressed as
\begin{align}
g_{\alpha\beta} = \Big[\Big((G_h)_{\kappa\sigma} + h_{,\kappa}h_{,\sigma}\Big)\circ\eta^\tau\Big]
\eta^\kappa_{,\alpha}\eta^\sigma_{,\beta}
:= \Big[{\mathcal G}_{\kappa\sigma}\circ\eta^\tau\Big]\eta^\kappa_{,\alpha}\eta^\sigma_{,\beta}. \label{gGrelation}
\end{align}
Note that %by the identity ${\mathcal G}_{\kappa\sigma} = (G_h)_{\kappa\sigma} + h_{,\kappa}h_{,\sigma}$,
${\mathcal G}_{\kappa\sigma}$ is the induced metric with respect to the
{\it normal} map $\eta^\nu$.
Furthermore, we have the following useful relationship between the determinant
of the two induced metrics:
\begin{align}
\det(g) = \det(\nabla_0\eta^\tau)^2\Big[\det(G_h) J_h^2\Big]\circ\eta^\tau
= \det(\nabla_0\eta^\tau)^2\Big[\det({\mathcal G})\Big]\circ\eta^\tau \label{detg}
\end{align}
where $\nabla_0$ denotes the surface gradient.

\begin{remark}
The identity (\ref{gGrelation}) can also be read as $(\eta^\tau)^*g = {\mathcal G}$.
\end{remark}

\vspace{.1 in}

Let $y$ and $\tilde{y} = \varphi(y)$ denote two different coordinate
systems on $\Gamma$ with associated metrics
\begin{align*}
g_{\alpha\beta} = \frac{\partial \eta^i}{\partial y^\alpha}\frac{\partial \eta^i}{\partial y^\beta}\,,\qquad
\tilde{g}_{\alpha\beta} = \frac{\partial \eta^i}{\partial \tilde{y}^\alpha}\frac{\partial \eta^i}{\partial \tilde{y}^\beta}.
\end{align*}
It follows that $\varphi^*\tilde{g} = g$. Let $H$, $\tilde{H}$, $K$,
$\tilde{K}$, $n$ and $\tilde{n}$ denote the mean curvature, Gauss
curvature, and the unit normal vector computed with respect to $y$
and $\tilde y$, respectively. Since $H$, $K$, and $n$ depend only on
the shape of $\Gamma(t)$, these geometric quantities are invariant
to tangential reparameterization; thus, the identity
\begin{align}
\tilde{H} = H\circ\varphi\,,\qquad \tilde{K} = K\circ\varphi\,,
\qquad \tilde{n} = n\circ\varphi. \label{HKnsymmetry}
\end{align}

Similarly, computing the first variation of
$\displaystyle{\int_{\Gamma(t)} H^2 dS}$ in our two coordinate systems yields
\begin{align*}
\Big[\Big(\Delta_g H + H(H^2 - K)\Big)n\Big](y) = \Big[\Big(\Delta_{\tilde{g}} \tilde{H}
+ \tilde{H}(\tilde{H}^2 - \tilde{K})\Big)\tilde{n}\Big](\tilde{y}) \qquad\forall\ \tilde{y}=\varphi(y).
\end{align*}
By (\ref{HKnsymmetry}), we have the following important identity
\begin{align}
\Big[\Delta_{\varphi^* \tilde{g}} H\Big](y) = \Big[\Delta_{\tilde{g}} (H\circ\varphi)\Big](\tilde{y})
\qquad\forall\ \tilde{y}=\varphi(y) \label{gHsymmetry1}
\end{align}
and hence
\begin{align}
[\Delta_{\mathcal G} (H\circ\eta^{-\tau})]\circ\eta^\tau = \Delta_g H \label{gHsymmetry2}
\end{align}
where by (\ref{quasilinearform}),
\begin{align}
H\circ\eta^{-\tau} =
-J_h^{-1}G_h^{\alpha\beta}\Big[\delta_{\beta\gamma} - J_h^{-2}
G^{\gamma\delta}_h h_{,\beta}h_{,\delta}\Big]h_{,\alpha\gamma}
+G_h^{\alpha\beta}F_{\alpha\beta}(y,h,\nabla h). \label{Hcompetatau}
\end{align}
\subsection{Bounds on $\eta^\tau$}\label{etatauremark}
Let $u^\tau$ denote the tangential velocity defined by
$\eta^\tau_t = u^\tau\circ\eta^\tau$.
Time-differentiating the relation $\eta=\eta^\nu\circ\eta^\tau$ and using
the definition of $\eta^\nu$, we find that
%\begin{align*}
%\eta^\tau_t = (\nabla_0 \eta^\nu)^{-1}\circ\eta^\tau (u\circ\eta - \eta^\nu_t\circ\eta^\tau)
%=\Big[(\nabla_0 \eta^\nu)^{-1} \Big(u\circ\eta^\nu - h_t \frac{\partial}{\partial z}\Big)\Big] \circ\eta^\tau
%\end{align*}
%for all $y\in \Gamma$, $t\in [0,T)$, where the inverse of $\nabla_0 \eta^\nu$ exists for a short time (since $\nabla_0
%\eta^\nu(0) =\id$). Therefore,
\begin{align}
u^\tau = (\nabla_0 \eta^\nu)^{-1} \Big[u\circ\eta^\nu - h_t \frac{\partial}{\partial z}\Big]\,. \label{utaudefn}
\end{align}
>From the trace theorem, it follows that
\begin{align}
\|u^\tau\|_{H^{2.5}(\Gamma)} \le C {\mathcal
P}(\|h\|_{H^{3.5}(\Gamma)},\|\eta\|_{H^3(\Omega)})
\Big[\|v\|_{H^3(\Omega)} + \|h_t\|_{H^{2.5}(\Gamma)}\Big]
\label{utauestimate}
\end{align}
for some polynomial ${\mathcal P}$. Since,
$ \eta^\tau(y,t) = y + \int_0^t (u^\tau\circ\eta^\tau)(y,s)ds $,
it follows that
\begin{align*}
\|\nabla_0 \eta^\tau(y,t)\|_{H^{1.5}(\Gamma)} \le C \Big[1+ \int_0^t
\|u^\tau\|_{H^{2.5}(\Gamma)} \Big(1 + \|\nabla_0
\eta^\tau\|_{H^{1.5}(\Gamma)}\Big)^4 ds\Big]
\end{align*}
and hence by Gronwall's inequality,
\begin{align}
\|\nabla_0 \eta^\tau(y,t)\|_{H^{1.5}(\Gamma)} \le C \Big[1+ \int_0^t
\|u^\tau\|_{H^{2.5}(\Gamma)} ds\Big] \label{etatauestimate}
\end{align}
for $t\in [0,T]$ sufficiently small. Furthermore, we also have
\begin{align}
\|\eta^\tau_t(y,t)\|_{H^{2.5}(\Gamma)} \le C
\|u^\tau\|_{H^{2.5}(\Gamma)} \Big[1+\|\nabla_0
\eta^\tau\|_{H^{1.5}(\Gamma)}\Big]^4. \label{etatautestimate}
\end{align}

\subsection{An expression for ${\mathfrak t}_{\bf ben}$ in terms of $h$ and $\eta^\tau$}
Now we can compute ${\mathfrak t}_{ben}$ in terms of $h$ and $\eta^\tau$: the highest order term of $\Delta_g H$ is
\begin{align*}
\Big\{\frac{1}{\sqrt{\det({\mathcal G})}}\frac{\partial}{\partial y^\gamma}\Big[\sqrt{\det({\mathcal G})}
{\mathcal G}^{\gamma\delta} \frac{\partial}{\partial
y^\delta} \Big(J_h^{-1}(G_h^{\alpha\beta} - J_h^{-2} G_h^{\alpha\kappa} G^{\beta\sigma}_h h_{,\kappa}h_{,\sigma})
h_{,\alpha\beta}\Big)\Big]\Big\}\circ\eta^\tau.
\end{align*}
Since ${\mathcal G}_{\alpha\beta} = (G_h)_{\alpha\beta} +
h_{,\alpha}h_{,\beta}$, the inverse of ${\mathcal G}_{\gamma\delta}$ is
\begin{align*}
\frac{1}{\det({\mathcal G})}\left[\begin{array}{cc}
(G_h)_{22} + h_{,2}^2 & -(G_h)_{12} - h_{,1} h_{,2} \\
-(G_h)_{12} - h_{,1} h_{,2} & (G_h)_{11} + h_{,1}^2
\end{array}
\right]
\end{align*}
which can also be written as
$${\mathcal G}^{\alpha\beta} = J_h^{-2} \Big[G_h^{\alpha\beta} -
(-1)^{\kappa+\sigma}\det(G_h)^{-1}(1-\delta_{\alpha\kappa})(1-\delta_{\beta\sigma})h_{,\kappa}h_{,\sigma}\Big].$$ Therefore,
the highest order term of $\Delta_g H$ can be written as
\begin{align*}
\frac{1}{\sqrt{\det(g_0)}}\Big[\sqrt{\det(g_0)}A^{\alpha\beta\gamma\delta} h_{,\alpha\beta}\Big]_{,\gamma\delta}\circ\eta^\tau
\end{align*}
where
\begin{align}
A^{\alpha\beta\gamma\delta} =&\ J_h^{-3} \Big[G_h^{\alpha\gamma} -
(-1)^{\kappa+\sigma}\det(G_h)^{-1}(1-\delta_{\alpha\kappa})(1-\delta_{\gamma\sigma})h_{,\kappa}h_{,\sigma}\Big]
\label{Aexpression}\\
&\times(G_h^{\beta\delta} - J_h^{-2} G_h^{\beta\kappa}G_h^{\delta\sigma}h_{,\kappa}h_{,\sigma}) \nonumber
\end{align}
is a fourth-rank tensor.

\subsection{Lagrangian formulation of the problem}
Let $\eta(t,x)=x+\int_0^t u(s,x)ds$  denote the Lagrangian particle
placement field, a volume-preserving embedding of $\Omega$ onto
$\Omega(t)\subset {\mathbb R}^3$,  and denote the cofactor matrix of
$\nabla \eta(x,t)$ by
\begin{align}
a(x,t) = [\nabla\eta(x,t)]^{-1} \label{adefn} \,.
\end{align}

Let $v=u\circ\eta$ denote the Lagrangian or material velocity field,
$q=p\circ\eta$ the Lagrangian pressure function, and $F=f\circ\eta$
the forcing function in the material frame. In the following
discussion, we also set $\epsilon = 1$. Then the system
(\ref{NSequation}) can be reformulated as
\begin{subequations}\label{NSequation1}
\begin{alignat}{2}
\eta_t &= v &&\quad\text{in } (0,T)\times\Omega,\label{NSequation1.a}\\
v_t^i - \nu(a^j_\ell D_\eta(v)_\ell^i)_{,j} &= -(a^k_i q)_{,k} + F^i &&\quad\text{in } (0,T)\times\Omega,\label{NSequation1.b}\\
a^k_i v^i_{,k} &= 0&&\quad\text{in } (0,T)\times\Omega,\label{NSequation1.c}\\
(\nu D_\eta(v)_\ell^i - q\delta_\ell^i )a^j_\ell N_j &= \sigma
\Theta
\Big[L(h)B_*(-G_h^{\alpha\beta}h_{,\alpha},1)\Big]\circ\eta^\tau
&&\quad\text{on } (0,T)\times\Gamma,
\label{NSequation1.d}\\
h_t &= B_*((-G_h^{\alpha\beta}h_{,\alpha},1))\cdot (v\circ\eta^{-\tau}) &&\quad\text{on}\ (0,T)\times\Gamma,\label{NSequation1.e}\\
v &= u_0 &&\quad\text{on } \{t=0\}\times\Omega,\label{NSequation1.f}\\
h &= 0 &&\quad\text{on }\{t=0\}\times\Gamma,\label{NSequation1.g}\\
\eta &= \id &&\quad\text{on }
\{t=0\}\times\Omega,\label{NSequation1.h}
\end{alignat}
\end{subequations}
where $D_\eta(v)_\ell^i:=(a^k_\ell v^i_{,k}+ a^k_i v^\ell_{,k})$,
$N$ denotes the outward-pointing unit normal to $\Gamma$, $\Theta$
is defined in Remark \ref{Thetaremark}, and $B_*$ is the
push-forward of $B$ defined as
$$B_* (\gamma'(0)) = (B\circ \gamma)'(0) \qquad\forall\ \gamma(t)\subset\Gamma.$$
$L(h)$ is the representation of
${\mathfrak t}_{shell}\cdot n$ using the height function $h$. It is defined as
follows
\begin{align*}
L(h) =&\
\frac{1}{\sqrt{\det(g_0)}}\Big[\sqrt{\det(g_0)}A^{\alpha\beta\gamma\delta}
h_{,\alpha\beta}\Big]_{,\gamma\delta}
+ L_1^{\alpha\beta\gamma}(y,h,Dh,D^2 h) h_{,\alpha\beta\gamma} \\
& + L_2(y,h,Dh,D^2 h)
\end{align*}
where $L_1$ and $L_2$ are polynomials of their variables with
$L_1(y,0)=0$, $g_0$ is the metric tensor on $\Gamma$. Note that
${\mathfrak t}_{mem}$ is included in $L_2$ since it is a second
order operator of $h$.

\begin{remark}\label{Thetaremark}
For a point $\eta(y,t)\in\Gamma(t)$, there are two ways of defining the
unit normal $n$ to $\Gamma(t)$:
\begin{itemize}
\item[1.] Let $n=\sqrt{g}^{-1} a^T N$ where $N$ is the unit normal to $\Gamma$.
\item[2.] Let $n=\displaystyle{
\Big[J_h^{-1}\Big(-G^{\alpha\beta}_h h_{,\alpha}\frac{\partial}{\partial y^\beta} +
\frac{\partial}{\partial z}\Big)\Big]\circ\eta^\tau}$ (denoted by $[J_h^{-1}(-\nabla_0 h,1)]\circ\eta^\tau$).
\end{itemize}
The function $\Theta$ is defined by
\begin{align*}
\Theta (-\nabla_0 h\circ\eta^\tau,1) = a^T N.
\end{align*}
Equating the modulus of both sides, by (\ref{detg}) we must have
\begin{align*}
\Theta = \sqrt{\det(g)}[(J_h^{-1})\circ\eta^\tau] = \det(\nabla_0\eta^\tau) \sqrt{\det(G_h)\circ\eta^\tau}.
\end{align*}
\end{remark}

\begin{remark}\label{simpleht}
An equivalent form of (\ref{NSequation1.e}) is given by
\begin{align*}
h_t =& -h_{,\alpha}(v\circ\eta^{-\tau})_{\alpha} +
(v\circ\eta^{-\tau})_z.
\end{align*}
This equation states that the shape of the boundary moves with the normal velocity of the fluid.
\end{remark}

\begin{remark}
For many of the nonlinear estimates that appear later, it is important that $L(h)$
is linear in the third derivative $h_{,\alpha\beta\gamma}$.
\end{remark}

\begin{remark}\label{importantremark}
Without using  the symmetry (\ref{gHsymmetry2}), we can still
compute $\Delta_g H$ in terms of $h$ and $\eta^\tau$ by using
(\ref{gGrelation}) and (\ref{detg}); however, $L_1$ would then
depend on $\nabla_0^2 \eta^\tau$ and thus lose one derivative of
regularity, preventing the closure of our energy estimate.
\end{remark}

\section{Notation and conventions}\label{notation}
\noindent
For $T>0$, we set
\begin{align*}
{\mathcal V}^1(T) =&\ \Big\{ v\in L^2(0,T;H^1(\Omega))\ \Big|\ v_t \in L^2(0,T;H^1(\Omega)')\Big\}; \\
{\mathcal V}^2(T) =&\ \Big\{ v\in L^2(0,T;H^2(\Omega))\ \Big|\ v_t \in L^2(0,T;L^2(\Omega))\Big\}; \\
{\mathcal V}^k(T) =&\ \Big\{ v\in L^2(0,T;H^k(\Omega))\ \Big|\ v_t \in L^2(0,T;H^{k-2}(\Omega))\Big\} \quad\text{for $k\ge 3$ }; \\
{\mathcal H}(T) =&\ \Big\{ h\in L^2(0,T;H^{5.5}(\Gamma))\ \Big|\ h_t
\in L^2(0,T;H^{2.5}(\Gamma)), h_{tt} \in
L^2(0,T;H^{0.5}(\Gamma))\Big\}
\end{align*}
with norms
\begin{align*}
\|v\|^2_{{\mathcal V}^1(T)} =&\ \|v\|^2_{L^2(0,T;H^1(\Omega))} + \|v_t\|^2_{L^2(0,T;H^1(\Omega)')} ;\\
\|v\|^2_{{\mathcal V}^2(T)} =&\ \|v\|^2_{L^2(0,T;H^2(\Omega))} + \|v_t\|^2_{L^2(0,T;L^2(\Omega))} ;\\
\|v\|^2_{{\mathcal V}^k(T)} =&\ \|v\|^2_{L^2(0,T;H^k(\Omega))} + \|v_t\|^2_{L^2(0,T;H^{k-2}(\Omega))} \quad\text{for $k\ge 3$ };\\
\|h\|^2_{{\mathcal H}(T)} =&\ \|h\|^2_{L^2(0,T;H^{5.5}(\Gamma))} +
\|h_t\|^2_{L^2(0,T;H^{2.5}(\Gamma))} +
\|h_{tt}\|^2_{L^2(0,T;H^{0.5}(\Gamma))}.
\end{align*}
We then introduce the space (of ``divergence free'' vector fields)
\begin{align*}
{\mathcal V}_v=\Big\{w\in H^1(\Omega)\ \Big|\ a_i^j(t)w^i_{,j} = 0\
\forall\ t\in[0,T]\Big\}
\end{align*}
and
\begin{align*}
{\mathcal V}_v(T)=\Big\{w\in L^2(0,T;H^1(\Omega))\ \Big|\
a_i^j(t)w^i_{,j} = 0\ \forall\ t\in[0,T]\Big\},
\end{align*}
where the cofactor matrix $a$ is defined by (\ref{adefn}). We use
$X_T$ to denote the space ${\mathcal V}^3(T)\times {\mathcal H}(T)$
with norm
\begin{align*}
\|(v,h)\|^2_{X_T} = \|v\|^2_{{\mathcal V}^3(T)} + \|h\|^2_{{\mathcal
H}(T)}
\end{align*}
and use $Y_T$, a subspace of $X_T$, to denote the space
\begin{align*}
Y_T = \Big\{(v,h)\in {\mathcal V}^3(T)\times {\mathcal H}(T)\ \Big|\
h_t\in L^\infty(0,T;H^2(\Gamma))\Big\}
\end{align*}
with norm
\begin{align*}
\|(v,h)\|^2_{Y_T} =&\ \|(v,h)\|^2_{X_T} + \|v\|^2_{L^\infty(0,T;H^2(\Omega))} + \|h\|^2_{L^\infty(0,T;H^4(\Gamma))} \\
& + \|h_t\|^2_{L^\infty(0,T;H^2(\Gamma))}.
\end{align*}

% Check with Steve if we need this remark.
%\begin{remark}
%By the Sobolev embedding theorem, if $(u,h)\in X(T)$, then $u\in L^\infty(0,T;H^2(\Omega))$ and $h\in
%L^\infty(0,T;H^4(\Gamma))$ with the following estimates:
%\begin{align*}
%\sup_{0\le t\le T} \|u\|^2_{H^2(\Omega)} \le \|u_0\|^2_{H^2(\Omega)} + \|u\|^2_{{\mathcal V}^3(T)}, \\
%\sup_{0\le t\le T} \|h\|^2_{H^4(\Omega)} \le \|h_0\|^2_{H^4(\Omega)} + \|h\|^2_{{\mathcal H}(T)},
%\end{align*}
%where $u_0$ and $h_0$ are the restriction of $u$ and $h$ to $\{t=0\}$, respectively.
%\end{remark}

We will solve (\ref{NSequation1}) by a fixed-point method in an appropriate subset of $Y_T$.

\section{The main theorem}\label{mainthm}
Before stating the main theorem, we define the following quantities. Let $q_0$ be defined by
\begin{subequations}\label{q0compatibility}
\begin{alignat}{2}
\Delta q_0 &= -\nabla u_0 : (\nabla u_0)^T + \nu [a^k_\ell D_\eta
(u_0)_\ell^i]_{,ki}(0) + \div F(0) &&\qquad \text{in}\quad\Omega,
\\
q_0 &= \nu (\Def u_0 \cdot N) \cdot N - \sigma L(0) &&\qquad
\text{on}\quad\Gamma
\end{alignat}
\end{subequations}
and
\begin{align}
u_1 = \nu\Delta u_0 - \nabla q_0 + F(0). \label{u1compatibility}
\end{align}
We also define the projection operator $P_{ij}(x): {\mathbb R}^3 \rightarrow
T_{\eta(x,t)} \Gamma(t)$ by
\begin{align*}
{\mathcal P}_{ij} (x) = [\delta_{ij} - (J_h^{-2}\circ\eta^\tau) a_i^k a_j^\ell N_k(x) N_\ell(x)]
=\Big[\delta_{ij} -
\frac{a_i^k N_k(x)}{|a_i^k N_k(x)|}\frac{a_j^\ell N_\ell(x)}{|a_j^\ell N_\ell(x)|}\Big].
\end{align*}

\begin{theorem}\label{maintheorem}
Let $\nu>0$, $\sigma>0$ be given, and
$$F\in L^2(0,T;H^2(\Omega)), F_t\in L^2(0,T;L^2(\Omega)), F(0)\in H^1(\Omega).$$
Suppose that the shell traction satisfies the compatibility
condition
\begin{align}
[\Def u_0\cdot N]_{\text{tan}} = 0.
\end{align}
There exists $T>0$ depending on $u_0$ and $F$ such that there exists
a solution $(v,h)\in Y_T$ of problem (\ref{NSequation1}). Moreover,
if $u_0\in H^{5.5}(\Omega)\cap H^{7.5}(\Gamma)$ and the associated
$u_1$, $q_0$ also satisfy the compatibility condition
\begin{align}
\text{CP} :=&\ \Big[g_0^{ki} u_{0,k}^j N_j N_\ell + g_0^{k\ell} u_{0,k}^j N_j N_i\Big]
\Big[\nu(\Def u_0)_i^j - q_0 \delta_i^j\Big]N_j \nonumber\\
&+ \nu(\delta_{i\ell} - N_i N_\ell)\Big[(\Def u_1)_i^j - \Big((\nabla u_0 \nabla u_0)
+ (\nabla u_0 \nabla u_0)^T \Big)_i^j \Big]N_j \label{compatibility1}\\
&- (\delta_{i\ell} - N_i N_\ell)\Big[\nu(\Def u_0)_i^j - q_0 \delta_i^j\Big] u_{0,j}^k N_k = 0 \nonumber
\end{align}
then the solution $(v,h)\in Y_T$ is unique.
\end{theorem}

\section{A bounded convex closed set of $Y_T$}\label{convexset}
\begin{definition}\label{CTM}
Given $M>0$. Let $C_T(M)$ denote the subset of $Y_T$ consisting of elements of $(v,h)$ in $Y_T$ such that
\begin{align}
\|(v,h)\|^2_{Y_T} \le M
\end{align}
and such that $v(0)=u_0$, $h(0)=0$ and $h_t(0)=(B_0)_*((0,1))\cdot u_0.$
\end{definition}

\begin{remark}\label{vtaudefn}
For $(v,h)\in C_T(M)$, define $u^\tau$ by (\ref{utaudefn}) and let $\eta^\tau$ be the associated flow map.
Also define $v^\tau$ as $u^\tau\circ\eta^\tau$. By (\ref{etatauestimate}) and (\ref{etatautestimate}), we have
\begin{align}
\sup_{t\in [0,T]}\|\nabla_0 \eta^\tau(t)\|_{H^{1.5}(\Gamma)} +
\|v^\tau\|^2_{L^2(0,T;H^{2.5}(\Gamma))} \le C(M)
\label{mainetatauestimate}
\end{align}
for some constant $C(M)$.
\end{remark}

We will make use of the following lemmas (proved in \cite{CoSh2005}):
\begin{lemma}
There exists $T_0\in (0,T)$ such that for all $T\in (0,T_0)$ and for all $v\in C_T(M)$, the matrix $a$ is
well-defined (by (\ref{adefn})) with the estimate (independent of $v\in C_T(M))$
\begin{align}
&\ \|a\|_{L^\infty(0,T;H^2(\Omega))} + \|a_t\|_{L^\infty(0,T;H^1(\Omega))} + \|a_t\|_{L^2(0,T;H^2(\Omega))} \nonumber\\
+&\ \|a_{tt}\|_{L^\infty(0,T;L^2(\Omega))} +
\|a_{tt}\|_{L^2(0,T;H^1(\Omega))} \le C(M). \label{aestimate1}
\end{align}
\end{lemma}

\begin{lemma} There exists $T_1\in (0,T)$ and a constant $C$ (independent of $M$) such that for all $T\in (0,T_1)$
and $v\in C_T(M)$, for all $\phi\in H^1(\Omega)$ and $t\in [0,T]$
\begin{align}
C\|\phi\|^2_{H^1(\Omega)} \le \int_{\Omega} \Big[|v|^2 +
|D_\eta(v)|^2\Big]dx  \label{equinorm}
\end{align}
where
$$|D_\eta(v)|^2:=D_\eta(v)^i_j D_\eta(v)^i_j=(a^k_j v^i_{,k}+ a^k_j v^i_{,k})(a^\ell_j v^i_{,\ell}+ a^\ell_i
v^j_{,\ell}).$$
\end{lemma}
\vspace{.1 in}

\noindent
In the remainder of the paper, we will assume that
$$0<T<\min\{T_0,T_1,\bar{T}\}$$
for some fixed $\bar{T}$ where the forcing $F$ is defined on the time interval $[0,\bar{T}]$.

\section{Preliminary results}\label{preliminaryresult}
\subsection{Pressure as a Lagrange multiplier}
In the following discussion, we use $\H1H2$ to denote the space
$H^1(\Omega)\cap H^2(\Gamma)$ with norm
$$\|u\|^2_{\H1H2} = \|u\|^2_{H^1(\Omega)} + \|u\|^2_{H^2(\Gamma)}$$
and $\bar{\mathcal V}_\ttv$ ($\bar{\mathcal V}_\ttv(T)$) to denote the space
\begin{align*}
\Big\{v\in {\mathcal V}_\ttv \ \Big|\ v\in H^2(\Gamma)\Big\} \Big(
\Big\{v\in {\mathcal V}_\ttv(T) \ \Big|\ v\in
L^2(0,T;H^2(\Gamma))\Big\} \Big).
\end{align*}

\begin{lemma}\label{divprob}
For all $p\in L^2(\Omega)$, $t\in [0,T]$, there exists a constant
$C>0$ and $\phi\in \H1H2$ such that $a_i^j(t)\phi_{,j}^i = p$ and
\begin{align}
\|\phi\|_{\H1H2} \le C\|p\|_{L^2(\Omega)}.\label{divestimate}
\end{align}
\end{lemma}
\begin{proof} We solve the following problem on the time-dependent domain $\Omega(t)$:
\begin{align*}
\div (\phi\circ\eta(t)^{-1}) = p\circ\eta(t)^{-1}
\qquad\text{in}\quad\eta(t,\Omega) := \Omega(t).
\end{align*}
The solution to this problem can be written as the sum of the solutions to the following two problems
\begin{alignat}{2}
\div (\phi\circ\eta(t)^{-1}) =&\ p\circ\eta(t)^{-1} -\bar{p}(t) &&\qquad\text{in}\quad\eta(t,\Omega), \label{divprob1}\\
\div (\phi\circ\eta(t)^{-1}) =&\ \bar{p}(t)
&&\qquad\text{in}\quad\eta(t,\Omega), \label{divprob2}
\end{alignat}
where $\displaystyle{\bar{p}(t)=\frac{1}{|\Omega|}\int_{\Omega}
p(t,x) dx}$. The existence of the solution to problem
(\ref{divprob1}) with zero boundary condition is standard (see, for
example, \cite{GPG1} Chapter 3), and the solution to problem
(\ref{divprob2}) can be chosen as a linear function (linear in $x$)
, for example, $\bar{p}(t) x_1$. The estimate (\ref{divestimate})
follows from the estimates of the solutions to (\ref{divprob1}).
\end{proof}

Define the linear functional on $\H1H2$ by
$(p,a_i^j(t)\varphi_{,j}^i)_{L^2(\Omega)}$ where $\varphi \in
\H1H2$. By the Riesz representation theorem, there is a bounded
linear operator $Q(t): L^2(\Omega)\to \H1H2$ such that for all
$\varphi\in\H1H2$,
\begin{align*}
(p,a_i^j(t)\varphi_{,j}^i)_{L^2(\Omega)} = (Q(t)p,\varphi)_{\H1H2}
:= (Q(t)p,\varphi)_{H^1(\Omega)} + (Q(t)p,\varphi)_{H^2(\Gamma)}.
\end{align*}
Letting $\varphi=Q(t)p$ shows that
$$\|Q(t)p\|_{\H1H2} \le C\|p\|_{L^2(\Omega)}$$
for some constant $C>0$. By Lemma \ref{divprob},
$$\|p\|^2_{L^2(\Omega)}\le \|Q(t)p\|_{\H1H2}\|\varphi\|_{\H1H2} \le C\|Q(t)p\|_{\H1H2}\|p\|_{L^2(\Omega)}$$
which shows that $R(Q(t))$ is closed in $\H1H2$. Since $\bar{\mathcal V}_v(t)\subset R(Q(t))^\perp$ and $R(Q(t))^\perp
\subset \bar{\mathcal V}_v(t)$, it follows that
\begin{align}
\H1H2(t)=R(Q(t))\oplus_{\H1H2} \bar{\mathcal V}_v(t).\label{H1decomp}
\end{align}
We can now introduce our Lagrange multiplier
\begin{lemma}\label{lagrangemultiplier}
Let ${\mathcal L}(t)\in \H1H2'$ be such that ${\mathcal L}(t)\varphi
= 0$ for any $\varphi\in \bar{\mathcal V}_v(t)$. Then there exist a
unique $q(t)\in L^2(\Omega)$, which is termed the pressure function,
satisfying
$$\forall\ \varphi\in \H1H2,\quad {\mathcal L}(t)(\varphi)=(q(t),a_i^j(t)\varphi_{,j}^i)_{L^2(\Omega)}.$$
Moreover, there is a $C>0$ (which does not depend on $t\in [0,T]$ and $\epsilon$ and on the choice of $v\in C_T(M)$) such that
$$\|q(t)\|_{L^2(\Omega)} \le C\|{\mathcal L}(t)\|_{\H1H2'}.$$
\end{lemma}
\begin{proof}
By the decomposition (\ref{H1decomp}), for given $\ta$, let $\varphi=v_1+v_2$, where $v_1\in {\mathcal V}_v(t)$ and $v_2\in
R(Q(t)$. It follows that
$${\mathcal L}(t)(\varphi)={\mathcal L}(t)(v_2) = (\psi(t),v_2)_{\H1H2} = (\psi(t),\varphi)_{\H1H2}$$
for a unique $\psi(t)\in R(Q(t))$.\\
>From the definition of $Q(t)$ we then get the existence of a unique
$q(t)\in L^2(\Omega)$ such that
$$\forall\ \varphi\in \H1H2,\quad{\mathcal L}(t)(\varphi)=(q(t),a_i^j(t) \varphi_{,j}^i)_{L^2(\Omega)}.$$
The estimate stated in the lemma is then a simple consequence of (\ref{divestimate}).
\end{proof}

% Check with Steve if we need to uncomment this section
%\subsection{Standard inequalities}
%\begin{lemma} Suppose $f \in H^{1.5}(\Gamma)$ and $g \in H^{0.5}(\Gamma)$, then $fg \in H^{0.5}(\Gamma)$ and satisfies
%\begin{align}
%\|fg\|_{H^{0.5}(\Gamma)} \le C\|f\|_{H^{1.5}(\Gamma)} \|g\|_{H^{0.5}(\Gamma)} \label{Hp5ineq}
%\end{align}
%for some constant $C$ depending on the geometry of $\Gamma$.
%\end{lemma}
%
%\begin{lemma}
%For space dimension $n=3$,
%\begin{align}
%\|f\|_{L^4(\Omega)} \le C\|f\|^{3/4}_{H^1(\Omega)} \|f\|^{1/4}_{L^2(\Omega)}\ , \label{interpolation1a}
%\end{align}
%\begin{align}
%\|f\|_{L^\infty(\Omega)} \le C\|f\|^{3/4}_{H^2(\Omega)} \|f\|^{1/4}_{L^2(\Omega)}\ , \label{interpolation2a}
%\end{align}
%\begin{align}
%\|f\|_{L^2(\Gamma)} \le C \|f\|^{1/2}_{H^1(\Omega)} \|f\|^{1/2}_{L^2(\Omega)}\ , \label{interpolation5}
%\end{align}
%while for space dimension $n=2$,
%\begin{align}
%\|f\|_{L^4(\Omega)} \le C\|f\|^{1/2}_{H^1(\Omega)} \|f\|^{1/2}_{L^2(\Omega)}\ , \label{interpolation1b}
%\end{align}
%\begin{align}
%\|f\|_{L^\infty(\Omega)} \le C\|f\|_{H^2(\Omega)}^{1/2}\|f\|_{L^2(\Omega)}^{1/2}\ , \label{interpolation2b}
%\end{align}
%\begin{align}
%\|f\|_{H^{0.5}(\Omega)} \le&\ C\|f\|^{1/2}_{H^1(\Omega)} \|f\|^{1/2}_{L^2(\Omega)}\ , \label{interpolation4}
%\end{align}
%and for either $n=2$ or $n=3$,
%\begin{align}
%\|v\|_{H^1(\Omega)} \le C \|v\|^{1/2}_{H^2(\Omega)} \|v\|^{1/2}_{L^2(\Omega)} \label{interpolation3}
%\end{align}
%where $C$ depends on $\Omega$.
%\end{lemma}

\subsection{Estimates for $a$ and $h$} We make use of near-identity transformations.
The following lemmas can be found in \cite{CoSh2002} and
\cite{CoSh2005}.
\begin{lemma}\label{a}
There exists $K>0$, $T_0>0$ such that if $0< t\le T_0$, then, for any $(\tv,\th)\in C_{T_0}(M)$,
\begin{subequations}\label{aestimate}
\begin{align}
\|\ta^T - \id\|_{L^\infty(0,T;{\mathcal C}^0(\overline{\Omega}_0))} &\le K\sqrt{t}\ ; \\
\|\ta - \id\|_{L^\infty(0,T;H^2(\Omega))} &\le K\sqrt{t}\ ; \\
\|\ta_t - \ta_t(0)\|_{L^\infty(0,T;H^1(\Omega))} &\le C(M) t\ ; \\
\|\ta_t\|_{L^\infty(0,T;H^1(\Omega))} &\le K.
\end{align}
\end{subequations}
\end{lemma}

We also need the following
\begin{lemma}\label{abc} For any $(\tv,\th)\in C_{T_0}(M)$,
\begin{align}
\|\th\|_{H^{3.5}(\Gamma)} \le C M t^{1/4} \label{hconvergence}
\end{align}
for all $0<t\le T_0$.
\end{lemma}
\begin{proof} For $(\tv,\th)\in C_T(M)$, $\|\th\|^2_{H^4(\Gamma)} + \|\th_t\|^2_{H^2(\Gamma)} \le M$.
By $\th(0)=0$,
\begin{align*}
\|\th(t)\|_{H^2(\Gamma)} \le \int_0^t \|\th_t\|_{H^2(\Gamma)} ds \le
\sqrt{M}t.
\end{align*}
Finally, the interpolation inequality
\begin{align}
\|\nabla_0^2 f(t)\|_{H^{1.5}(\Gamma)} \le C\|\nabla_0^4
f\|^{3/4}_{L^2(\Gamma)} \|\nabla_0^2 f\|^{1/4}_{L^2(\Gamma)},
\label{interpolation6}
\end{align}
implies
\begin{align*}
\|\th\|_{H^{3.5}(\Gamma)} \le C\|\th\|^{3/4}_{H^4(\Gamma)}
\|\th\|^{1/4}_{H^2(\Gamma)} \le C M t^{1/4}.
\end{align*}
\end{proof}

\begin{corollary}\label{Lconvergence}
$\|L_1(t)\|_{H^{1.5}(\Gamma)}$ and $\|L_2(t)\|_{H^{1.5}(\Gamma)}$
converge to zero as $t\to 0$, uniformly in $(v,h)\in C_{T_0}(M)$.
Furthermore, for $t \le 1$,
\begin{align*}
\|L_1(t)\|_{H^{1.5}(\Gamma)} + \|L_2(t)\|_{H^{1.5}(\Gamma)} \le
C(M)t^{1/4}.
\end{align*}
\end{corollary}

By the fact that $\|\th_t\|^2_{H^2(\Gamma)}\le M$ and
$\|\th_{tt}\|^2_{L^2(0,T;H^{0.5}(\Gamma))}\le M$ if $(\tv,\th)\in
C_T(M)$, similar computations lead to the following lemma.
\begin{lemma}
For all $(\tv,\th)\in C_{T}(M)$,
\begin{align}
\|\th_t(t)\|_{H^{1.5}(\Gamma)} \le C M t^{1/8} \label{htconvergence}
\end{align}
for all $0<t\le T$.
\end{lemma}

%\subsection{Inequality for $\Theta$}
%We also need the following inequality
%\begin{align}
%\|\nabla_0\ttTheta\|_{H^{1.5}(\Gamma)} + \|\nabla_0 b\|_{H^{1.5}(\Gamma)} + \|B_{,\gamma}\|_{H^{1.5}(\Gamma)}
%\le C(\epsilon')t^{1/2} \label{ThetaB}
%\end{align}
%where $b$ and $B$ are defined in Appendix \ref{L2H3inequality}.

\section{The linearized problem}\label{linearizedprob}
Suppose that $(\tv,\th)\in C_T(M)$ is given. Let $\teta(t)= \id +
\int_0^t \tv(s)ds$ and $\ta =(\nabla\teta)^{-1}$. We are concerned
with the following time-dependent linear problem, whose fixed-point
$v=\tv$ provides a solution to (\ref{NSequation1}):
\begin{subequations}\label{NSequation0}
\begin{alignat}{2}
v^i_t - \nu [\ta_\ell^k D_\teta(v)_\ell^i]_{,k} &= -(\ta_i^k q)_{,k} + F^i&&\ \ \text{in}\ \
(0,T)\times\Omega\,,\label{NSequation0.a}\\
\ta_i^j v^i_{,j} &= 0&&\ \ \text{in}\ \ (0,T)\times\Omega\,,\label{NSequation0.b}\\
[\nu D_\teta(v)^j_i -q\delta^j_i]\ta_j^\ell N_\ell &= \sigma \tTheta\Big[\L_\th(h)(-\nabla_0\th,1)\Big]\circ\teta^\tau &&
\ \ \text{on}\ \ (0,T)\times\Gamma\,,\label{NSequation0.c}\\
&\quad + \sigma \tTheta \Big[[\G(\th)(-\nabla_0\th,1)]\circ\teta^\tau\Big] \nonumber\\
h_t\circ\teta^\tau &= [\th_{,\alpha}\circ\teta^\tau] v_{\alpha} - v_z &&\ \ \text{on}\ \ (0,T)\times\Gamma\,,\label{NSequation0.d}\\
v &= u_0&&\ \ \text{on}\ \ \{t=0\}\times\Omega \,,\label{NSequation0.e}\\
h &= 0 &&\ \ \text{on}\ \ \{t=0\}\times\Gamma
\,.\label{NSequation0.f}
\end{alignat}
\end{subequations}
where
$D_\teta(v)_i^j=\ta_i^k v^j_{,k}+\ta_j^kv^i_{,k}$, $\tTheta =\det(\nabla_0\teta^\tau)$, and
\begin{align*}
\L_\th(h) =&\ \frac{1}{\sqrt{\det(g_0)}}\Big[\sqrt{\det(g_0)}\tA^{\alpha\beta\gamma\delta} h_{,\alpha\beta}\Big]_{,\gamma\delta}
\end{align*}
with
\begin{align*}
\tA^{\alpha\beta\gamma\delta} =&\ J_\th^{-3} \sqrt{\det(G_\th)}\Big[G_\th^{\alpha\gamma} -
(-1)^{\kappa+\sigma}\det(G_\th)^{-1}(1-\delta_{\alpha\kappa})(1-\delta_{\gamma\sigma})\th_{,\kappa}\th_{,\sigma}\Big] \\
&\times(G_\th^{\beta\delta} - J_\th^{-2}G_\th^{\beta\mu}G_\th^{\delta\nu}\th_{,\mu}\th_{,\nu})
\end{align*}
and
\begin{align*}
\G(\th) =&\ \sqrt{\det(G_\th)\circ\teta^\tau}\Big[L_1^{\alpha\beta\gamma}(y,\th,D\th,D^2 \th) \th_{,\alpha\beta\gamma}
+ L_2(y,\th,D\th,D^2 \th)\Big].
\end{align*}
Here the thickness $\epsilon$ is assumed to be 1.

We will also use $L_\th(h)$ to denote $\L_\th(h) + \G(\th)$.

\begin{remark}
$\L_\th$ is a coercive fourth order operator for small $\th \le
\delta$. Actually, it is easy to see that $\L_\th$ is coercive at
time $t=0$, and the coercivity of $\L_\th$ for $t>0$ (but
sufficiently small) follows from the continuity of $\th$ in time
into the space $H^2(\Gamma)$. Moreover, by Lemma \ref{abc}, we have
the following corollary.
\end{remark}

\begin{corollary}\label{ellipticconstant}
There exists a $\nu_1 > 0$ and $0<T\le T_0$ such that for all $0<t\le T$,
\begin{align*}
\nu_1 \|\nabla_0^2 f(t)\|^2_{L^2(\Gamma)} \le \int_{\Gamma}
\tA^{\alpha\beta\gamma\delta} f_{,\alpha\beta}(t)
f_{,\gamma\delta}(t) dS.
\end{align*}
for all $0<t\le T$. Later on we will denote the right-hand side quantity of this inequality by $E_\tth(f)$,
where the subscript $\tth$ indicates that $\ttA$ is a function of $\tth$.
\end{corollary}

\begin{remark}
Given $(\tv,\th)\in {\mathcal V}^3(T)\times {\mathcal H}(T)$, for
the corresponding $\teta^\tau$, we have
\begin{align*}
\|\teta^\tau\|^2_{L^\infty(0,T;H^{2.5}(\Omega))} +
\|\teta^\tau_t\|^2_{L^2(0,T;H^{2.5}(\Gamma))} \le C(M)
\end{align*}
where (\ref{etatautestimate}) and (\ref{etatauestimate}) are used to obtain this estimate.
\end{remark}
\vspace{.1 in}
\noindent
The solution of (\ref{NSequation0}) is found as a weak limit of a sequence
of regularized problems.
\begin{definition}{\bf (Mollifiers on $\Gamma$)}
For $\epsilon>0$, let
$$K_\epsilon^p := (1-\epsilon \Delta_0)^{-\frac{p}{2}}: H^s(\Gamma) \to H^{s+p}(\Gamma)$$
denote the usual self-adjoint Frederich mollifier on the compact
manifold $\Gamma$, where $\Delta_0$ is the surface
Laplacian defined on $\Gamma$.%, given by
%$$\Delta_0 f = \frac{1}{\sqrt{\det(g_0)}}\frac{\partial}{\partial y^\alpha}
%\Big(\sqrt{\det(g_0)}g^{\alpha\beta}\frac{\partial f}{\partial y^\beta}\Big).$$
\end{definition}

By the Sobolev extension theorem, there exist bounded extension operators
\begin{align*}
E_s&: H^s(\Omega) \to H^s(\bbR^n), \ \ s \ge 1\,.
\end{align*}
For fixed (but small) $\epsilon$ and $\epsilon_1>0$, let
$\rho_{\epsilon}$ be a (positive) smooth mollifier on $\bbR^n$. Set
$\ttv= \rho_{\epsilon}\ast E_1(\tv)$, $\tF = \rho_{\epsilon}\ast
E_2(F)$, $\tu_0= \rho_{\epsilon}\ast E_3(u_0)$, where $\ast$ denotes
the convolution in space, and $\bar{h} = K_{\epsilon}^m(\tilde{h})$
for large enough $m$. Define $\tteta$ and $\tta$ in the same fashion
as $\teta$ and $\ta$. Note that $\ttv \to\tv\in V(T)$, $\tF\to F$ in
${\mathcal V}^2(T)$, $\tilde{u_0} \to u_0$ in $H^{2.5}(\Omega)$ and
$\bar{h} \to \tilde{h}$ in ${\mathcal H}(T)$ as $\epsilon\to 0$.

The regularized problem takes the form
\begin{subequations}\label{NSequation2}
\begin{alignat}{2}
v^i_t - \nu [\tta_\ell^k D_\tteta(v)_\ell^i]_{,k} &= -(\tta_i^k q)_{,k} + \tF^i&&\ \ \text{in}\ \ (0,T)\times\Omega\,,\label{NSequation2.a}\\
\tta_i^j v^i_{,j} &= 0&&\ \ \text{in}\ \ (0,T)\times\Omega\,,\label{NSequation2.b}\\
[\nu D_\tteta(v)^j_i -q\delta^j_i]\tta_j^\ell N_\ell &= \sigma \L^\eps1_\tth(h^\eps1)(-\nabla_0\tth\circ\tteta^\tau,1)
\nonumber\\
&\quad + \sigma \G^\eps1_\tth(-\nabla_0\tth\circ\tteta^\tau,1) + \kappa\Delta_0^2 v &&
\ \ \text{on}\ \ (0,T)\times\Gamma\,,\label{NSequation2.c}\\
h_t\circ\tteta^\tau &= [(\tth_{,\alpha})\circ\tteta^\tau] v_{\alpha} - v_z &&\ \ \text{on}\ \ (0,T)\times\Gamma\,,\label{NSequation2.d}\\
v &= \tu_0&&\ \ \text{on}\ \ \{t=0\}\times\Omega \,,\label{NSequation2.e}\\
h &= 0 &&\ \ \text{on}\ \ \{t=0\}\times\Gamma
\,,\label{NSequation2.f}
\end{alignat}
\end{subequations}
where
\begin{align*}
&\bL^\eps1_\tth(f) = \frac{\ttTheta}{\sqrt{\det(g_0)}}\Big[\Big(\sqrt{\det(g_0)}\ttA^{\alpha\beta\gamma\delta}
f_{,\alpha\beta}\Big)_{,\gamma\delta}\Big]^\eps1\circ\tteta^\tau\,, \\
&\bG_\tth^\eps1 =
\ttTheta\Big[\Big(L_1^{\alpha\beta\gamma}(\cdot,\tth,D\tth,D^2\tth)\tth_{,\alpha\beta\gamma} + L_2(\cdot,\tth,D\tth)
\Big)^\eps1\Big]^\eps1\circ\tteta^\tau(y,t).
\end{align*}
Note that $\displaystyle{\bL_\tth^\eps1(f) + \bG_\tth^\eps1 = \ttTheta\Big[L_\tth(f)\Big]^\eps1\circ\tteta^\tau}$.
%\begin{align*}
%\bL_\tth^\eps1(f) + \bG_\tth^\eps1 = \ttTheta\Big[L_\tth(f)\Big]^\eps1\circ\tteta^\tau.
%\end{align*}

\subsection{Weak solutions}
\begin{definition} A vector $v\in \bar{\mathcal V}_{\ttv}(T)$ with $v_t\in \bar{\mathcal V}_\ttv(T)'$ for almost all $t\in (0,T)$ is a weak
solution of (\ref{NSequation2}) provided that
\begin{subequations}\label{weakform}
\begin{align}
\text{\rm (i)}&\ \langle v_t,\varphi\rangle +
\frac{\nu}{2}\int_\Omega D_\tteta v : D_\tteta\varphi dx +
\sigma\int_{\Gamma} \ttA^{\alpha\beta\gamma\delta}
h_{,\alpha\beta}^\eps1
\Big[-\tth_{,\sigma}(\varphi^\sigma\circ\tteta^{-\tau}) \\
& + (\varphi^z\circ\tteta^{-\tau})\Big]^\eps1_{,\gamma\delta}dS +
\kappa\int_{\Gamma} \Delta_0 v\cdot \Delta_0 \varphi dS = \langle
\tF, \varphi\rangle - \sigma \langle
\G_\tth^\eps1,\varphi\rangle_{\Gamma}
\nonumber\\
\text{\rm (ii)} &\ v(0,\cdot)=\tu_0
\end{align}
\end{subequations}
for almost all $t\in [0,T]$, where $\langle \cdot,\cdot\rangle$ denotes the duality product between $\bar{\mathcal V}_v(t)$ and
its dual $\bar{\mathcal V}_v(t)'$, and $h$ is given by the evolution equation (\ref{NSequation2.d}) and the initial condition
(\ref{NSequation2.f}):
\begin{align}
h(y,t) = \int_0^t \Big[-\tth_{,\alpha}(y,s)v^\alpha(\tteta^{-\tau}(y,s),0,s) + v^z(\tteta^{-\tau}(y,s),0,s)\Big] ds
\label{hevol}
\end{align}
\end{definition}

\subsection{Penalized problems}
Letting $\theta>0$ denote the penalized parameter, we define $w_\theta$
(with also $\epsilon$ and $\epsilon_1$ dependence in mind) to be
the ``unique'' solution of the problem (whose existence can be obtained via a modified Galerkin method which will be
presented in the following sections):
\begin{subequations}\label{penalizedprob1}
\begin{align}
&\text{(i)} \ \langle w_{\theta t},\varphi\rangle +
\frac{\nu}{2}\int_\Omega D_\tteta w_\theta : D_\tteta\varphi dx +
\sigma\int_{\Gamma} \ttA^{\alpha\beta\gamma\delta}
h_{,\alpha\beta}^\eps1
\Big[-\tth_{,\sigma}(\varphi^\sigma\circ\tteta^{-\tau}) \nonumber\\
&\ + (\varphi^z\circ\tteta^{-\tau})\Big]^\eps1_{,\gamma\delta}dS +
\kappa\int_{\Gamma} \Delta_0 v\cdot \Delta_0 \varphi dS
+ (\frac{1}{\theta}\tta_i^jv_{,j}^i, \tta_k^\ell\varphi_{,\ell}^k)_{L^2(\Omega)} \\
&\ = \langle \tF, \varphi\rangle - \sigma \langle \bG_\tth^\eps1 (-\nabla_0\tth\circ\tteta^\tau,1)
,\varphi\rangle_{\Gamma} \nonumber\\
&\text{(ii)} \ v(0,\cdot)=\tu_0
\end{align}
\end{subequations}
where $\langle \cdot,\cdot\rangle$ denotes the pairing between
$H^1(\Omega)$ and its dual, and $h$ in (\ref{penalizedprob1}a)
satisfies (\ref{hevol}) with $v$ replaced by $w_\theta$.

\subsection{Weak solutions for the penalized problem}
The goal of this section is to establish the existence of $v$ to the problem (\ref{NSequation2}) (or the weak formulation
(\ref{weakform})), as well as the energy inequality satisfied by $v$ and $v_t$. Before proceeding,
we introduce variable $\tq_0$ and $\tw_1$ as follows:
let $\tq_0$ be the solution of the following Laplace equation
\begin{subequations}\label{tildecompatibility1}
\begin{alignat}{2}
\Delta \tq_0 &= \nabla \tu_0:(\nabla \tu_0)^t - \div \tF(0) &&\qquad \text{in $\Omega$}\,, \\
\tq_0 &= \nu (\Def \tu_0)_i^j N_i N_j - \sigma \G_0^\eps1(0) +
\kappa\Delta_0^2 \tu_0\cdot N &&\qquad\text{on $\Gamma$}\,,
\end{alignat}
\end{subequations}
and $\tw_1$ be defined by
\begin{align}
\tw_1 = \nu\Delta\tu_0 - \nabla \tq_0 + \tF(0). \label{tildecompatibility2}
\end{align}
By elliptic regularity,
\begin{align*}
\|\tq_0\|^2_{H^1(\Omega)} &\le C\Big[ \|\tu_0\|^2_{H^2(\Omega)} +
\|\tF(0)\|^2_{L^2(\Omega)} +
\|\G_0^\eps1(0)\|^2_{H^{0.5}(\Gamma)} + \|\Delta_0^2 \tu_0\|^2_{H^{0.5}(\Gamma)} \Big]\\
&\le C(M)\Big[\|\tu_0\|^2_{H^2(\Omega)} +
\|\tu_0\|^2_{H^{4.5}(\Gamma)} + \|\tF(0)\|^2_{L^2(\Omega)} +
1\Big]\,,
\end{align*}
and hence
\begin{align*}
\|\tw_1\|^2_{L^2(\Omega)} \le C(M)\Big[\|\tu_0\|^2_{H^2(\Omega)} +
\|\tu_0\|^2_{H^{4.5}(\Gamma)} + \|\tF(0)\|^2_{L^2(\Omega)} + 1\Big].
\end{align*}

\begin{remark}
By (\ref{hconvergence}), the constant $C(M)$ in the estimates above
can also be refined as a constant independent of $M$ if $T$ is
chosen small enough.
\end{remark}
\vspace{.1 in}

By introducing a (smooth) basis $(e_\ell)_{\ell=1}^\infty$ of $\H1H2$, and taking the approximation at rank $m\ge 2$ under the
form $w_\ell(t,x) = \sum\limits_{k=1}^\ell d_k(t)e_k(x)$ with
\begin{align}
h_\ell (y,t) &= \int_0^t \Big[-\tth_{,\alpha}(y,s)w_\ell^\alpha(\tteta^{-\tau}(y,s),0,s) +
w_\ell^z(\tteta^{-\tau}(y,s),0,s)\Big] ds, \label{hevol2}
\end{align}
and satisfying on $[0,T]$,
\begin{subequations}\label{penalizedprob}
\begin{align}
&\text{(i) } \ (w_{\ell tt},\varphi)_{L^2(\Omega)} + \nu (\tta_i^j
w_{\ell t,j},\tta_i^k \varphi_{,k})_{L^2(\Omega)}
+ \nu ((\tta_i^j\tta_i^k)_t w_\ell, \varphi_{,k})_{L^2(\Omega)} \nonumber\\
&\ + \nu\int_{\Omega} \Big[\tta_r^j\tta^k_i w_{\ell t,j}^i +
(\tta_r^j\tta^k_i)_t w_{\ell,j}^i\Big] \varphi_{,k}^r dx + \kappa
\int_{\Gamma} \Delta_0 w_{\ell t}\cdot \Delta_0 \varphi dS -
((\tta^j_i q_\ell)_t,\varphi^i_{,j})_{L^2(\Omega)}
\nonumber\\
&\ + \sigma \int_{\Gamma} \ttA^{\alpha\beta\gamma\delta}
[-\tth_{,\sigma} (w_\ell^\sigma\circ\tteta^{-\tau}) +
w_\ell^z\circ\tteta^{-\tau}]^\eps1_{,\alpha\beta}[-\tth_{,\sigma}
(\varphi^\sigma\circ\tteta^{-\tau}) +
\varphi^z\circ\tteta^{-\tau}]^\eps1_{,\gamma\delta}dS \nonumber\\
&\ + \sigma\int_{\Gamma} (\ttA^{\alpha\beta\gamma\delta})_t
h^\eps1_{\ell,\alpha\beta} [-\tth_{,\sigma}
(\varphi^\sigma\circ\tteta^{-\tau}) +
\varphi^z\circ\tteta^{-\tau}]^\eps1_{,\gamma\delta}dS \\
&\ + \sigma\int_{\Gamma} \ttA^{\alpha\beta\gamma\delta}
h^\eps1_{\ell,\alpha\beta} [-\tth_{t,\sigma}
(\varphi^\sigma\circ\tteta^{-\tau}) + \tth_{,\sigma} \ttv^\kappa
(\varphi^\sigma_{,\kappa}
\circ\tteta^{-\tau}) + \ttv^\kappa (\varphi^z_{,\kappa}\circ\tteta^{-\tau})]^\eps1_{,\gamma\delta}dS \nonumber\\
&\ = \langle \tF_t,\varphi\rangle - \sigma \int_{\Gamma}
\Big[L_1^{\alpha\beta\gamma}\tth_{,\alpha\beta\gamma} +
L_2\Big]^\eps1_t \Big[\tth_{,\sigma}
(\varphi^\sigma\circ\tteta^{-\tau}) - \varphi^z\circ\tteta^{-\tau}\Big]^\eps1 dS \nonumber\\
&\ - \sigma\int_{\Gamma}
\Big[L_1^{\alpha\beta\gamma}\tth_{,\alpha\beta\gamma} +
L_2\Big]^\eps1 \Big[\tth_{t,\sigma}
(\varphi^\sigma\circ\tteta^{-\tau}) - \tth_{,\sigma} \ttv^\kappa
(\varphi^\sigma_{,\kappa}
\circ\tteta^{-\tau}) - \ttv^\kappa (\varphi^z_{,\kappa}\circ\tteta^{-\tau}) \Big]^\eps1 dS \nonumber\\
&\ \ \forall\ \varphi\in\text{span}(e_1,\cdots,e_\ell)\,,\nonumber\\
&\text{(ii) } \ w_{\ell t}(0)= (w_1)_\ell, w_\ell(0) = (u_0)_\ell
\text{ in $\Omega$}\,,
\end{align}
\end{subequations}
where $\displaystyle{q_\ell = \tq_0 -\frac{1}{\theta}\tta_i^j w_{\ell,j}^i}$,
and $(\tu_0)_\ell$ denote the respective $\H1H2$ projections of
$u_0$ on span($e_1,e_2,\cdots,e_\ell$).

\begin{remark}
The existence of $w_k$ follows from the solution of
\begin{align*}
d_k''(t) + d_\ell'(t) A_{k\ell}(t) + d_\ell(t) B_{k\ell}(t) + \int_0^t d_\ell(s) C_{k\ell}(s,t)ds = F(t)
\end{align*}
for functions $A$, $B$, $C$ and $F$; however, the existence of the solution $d_k$ does not immediately follow from the
fundamental theorem of ODE due to the presence of the time-integral. A straightforward fix-point argument can be implemented,
whose details we leave to interested reader.
\end{remark}
\vspace{.1 in}

\noindent

The use of the test function $\varphi=w_{\ell t}$ in this system of ODE gives us in turn the energy law
\begin{align}
&\ \ \frac{1}{2}\frac{d}{dt} \|w_{\ell t}\|^2_{L^2(\Omega)} +
\frac{\nu}{2} \|D_\tteta (w_{\ell t})\|^2_{L^2(\Omega)}
+ \frac{\sigma}{2}\frac{d}{dt} E_\tth(h_{\ell t,\alpha\beta}^\eps1) + \theta \|q_{\ell t}\|^2_{L^2(\Omega)} \nonumber\\
&\ + \nu ((\tta_i^j\tta_i^k)_t w_{\ell,j}, w_{\ell
t,k})_{L^2(\Omega)} + \nu\int_{\Omega} (\tta_r^j\tta^k_i)_t
w_{\ell,j}^i w_{\ell t,k}^r dx
+ \kappa \|\Delta_0 w_{\ell t}\|^2_{L^2(\Gamma)} \nonumber\\
&\ + (q_{\ell t}, \tta_{it}^j w_{\ell,j}^i)_{L^2(\Omega)} - (q_\ell,
\tta_{it}^j w_{\ell t,j}^i)_{L^2(\Omega)} -
\frac{\sigma}{2}\int_{\Gamma} (\ttA^{\alpha\beta\gamma\delta})_t
h^\eps1_{\ell t,\alpha\beta}
h^\eps1_{\ell t,\gamma\delta} dS \nonumber\\
&\ - \sigma \int_{\Gamma} (\ttA^{\alpha\beta\gamma\delta})_t
h^\eps1_{\ell,\alpha\beta} \Big[h_{\ell tt} + \tth_{t,\sigma}
(w_{\ell t}^\sigma\circ\tteta^{-\tau})\Big]^\eps1_{,\gamma\delta} dS
+ \sigma\int_{\Gamma} \ttA^{\alpha\beta\gamma\delta} h^\eps1_{\ell,\alpha\beta}\times \label{welltweak}\\
&\ \times\Big[-\tth_{t,\sigma} (w_{\ell
t}^\sigma\circ\tteta^{-\tau}) + \tth_{,\sigma} \ttv^\kappa (w_{\ell t,\kappa}^\sigma\circ\tteta^{-\tau})
+ \ttv^\kappa(w_{\ell t,\kappa}^z\circ\tteta^{-\tau})\Big]^\eps1_{,\gamma\delta} dS \nonumber\\
& = \langle \tF_t,w_{\ell t}\rangle - \sigma\int_{\Gamma}
\Big[(L_1^{\alpha\beta\gamma} \tth_{,\alpha\beta\gamma} +
L_2)(-\nabla_0 \tth,1)\Big]_t\cdot
(w_{\ell t}\circ\tteta^{-\tau}) dS \nonumber\\
&\ - \sigma\int_{\Gamma} (L_1^{\alpha\beta\gamma}
\tth_{,\alpha\beta\gamma} + L_2) \ttv^\kappa\Big[-\tth_{,\sigma}
(w^\sigma_{\ell t,\kappa} \circ\tteta^{-\tau}) + (w^z_{\ell
t,\kappa} \circ\tteta^{-\tau})\Big] dS. \nonumber
\end{align}
For the tenth term (the integral with
$\displaystyle{\frac{\sigma}{2}}$ as its coefficient), we have
\begin{align*}
\Big|\int_{\Gamma} (\ttA^{\alpha\beta\gamma\delta})_t h^\eps1_{\ell
t,\alpha\beta} h^\eps1_{\ell t,\gamma\delta} dS\Big| \le
C(M)\|\tth_t\|_{H^{2.5}(\Gamma)} \|\nabla_0^2 h_{\ell
t}\|^2_{L^2(\Gamma)}.
\end{align*}
By $\eps1$-regularization and the identity
\begin{align*}
\int_{\Gamma} (\ttA^{\alpha\beta\gamma\delta})_t
h^\eps1_{\ell,\alpha\beta} h^\eps1_{\ell tt,\gamma\delta} dS
=&\int_{\Gamma} \frac{1}{\sqrt{\det(g_0)}}
\Big[\sqrt{\det(g_0)}(\ttA^{\alpha\beta\gamma\delta})_t\Big]_{,\gamma\delta}
h^\eps1_{\ell,\alpha\beta} h^\eps1_{\ell tt} dS \\
& + \int_{\Gamma} \frac{2}{\sqrt{\det(g_0)}}
\Big[\sqrt{\det(g_0)}(\ttA^{\alpha\beta\gamma\delta})_t\Big]_{,\gamma}
h^\eps1_{\ell,\alpha\beta\delta} h^\eps1_{\ell tt} dS \\
& + \int_{\Gamma} (\ttA^{\alpha\beta\gamma\delta})_t
h^\eps1_{\ell,\alpha\beta\gamma\delta} h^\eps1_{\ell tt} dS,
\end{align*}
we find that
\begin{align*}
&\ \Big|\int_{\Gamma} (\ttA^{\alpha\beta\gamma\delta})_t h^\eps1_{\ell,\alpha\beta} h^\eps1_{\ell tt,\gamma\delta} dS\Big|\\
\le&\ C(\eps1)\Big[1+\|\tth_t\|_{H^{2.5}(\Gamma)}\Big]\|\nabla_0^2
h_\ell\|_{L^2(\Gamma)} \Big[\|w_\ell\|_{H^1(\Omega)} + \|w_{\ell
t}\|_{H^1(\Omega)}\Big].
\end{align*}
Similarly, the second part of the eleventh term and the last term of the left-hand side can be bounded by
\begin{align*}
C(\eps1)\|\tth_t\|_{H^{2.5}(\Gamma)} \|\nabla_0^2
h_\ell\|_{L^2(\Gamma)} \|w_{\ell t}\|_{H^1(\Omega)}
\end{align*}
where we also use the $\eps1$-regularization to control $\nabla_0^3 w_{\ell t}$.
It also follows that the last two terms on the right-hand side can be bounded by
\begin{align*}
C(M)\Big[1+\|\tth_t\|_{H^{2.5}(\Gamma)}\Big]\|w_{\ell
t}\|_{H^1(\Omega)}.
\end{align*}

With positive $\theta$, the fourth term of the left-hand side involving the square of $q_{\ell t}$ acts as a viscous energy
term. Integrating (\ref{welltweak}) in time from $0$ to $t$, we then get
\begin{align}
&\ \|w_{\ell t}\|^2_{L^2(\Omega)} + \|\nabla_0^2 h_{\ell
t}\|^2_{L^2(\Gamma)} + \int_0^t \Big[\|\nabla w_{\ell
t}\|^2_{L^2(\Omega)} + \kappa \|w_{\ell t}\|^2_{H^2(\Gamma)}
+ \theta \|q_{\ell t}\|^2_{L^2(\Omega)}\Big] ds \nonumber\\
\le&\ C(M)\Big[\|w_{\ell t}(0)\|^2_{L^2(\Omega)} +
\|w_\ell(0)\|^2_{H^1(\Omega)}
+ \|q_\ell(0)\|^2_{H^{0.5}(\Omega)}\Big] \label{ellL2H1vt}\\
& + C(\eps1) \int_0^t \Big[1+\|\tth_t(s)\|^2_{H^{2.5}(\Gamma)}\Big]\|\nabla_0^2 h_{\ell t}(s)\|^2_{L^2(\Gamma)} ds \nonumber\\
& + C(\theta)\int_0^t \|\ttv(t')\|^2_{H^3(\Omega)} \int_0^{t'}
\Big[\|\nabla w_{\ell t}(s)\|^2_{L^2(\Omega)} + \|q_{\ell
t}(s)\|^2_{L^2(\Omega)}\Big] ds dt', \nonumber
\end{align}
where $C(\eps1), C(\theta) \to \infty$ as $\eps1,\theta\to 0$, and we use
\begin{align*}
\|f(t)\|_X \le \|f(0)\|_X + \int_0^t \|f_t(s)\|_X ds \le \|f(0)\|_X + \sqrt{t} \int_0^t \|f_t(s)\|^2_X ds
\end{align*}
for $f=w_\ell$, $f=h_\ell$ and $f=g_\ell$ to obtain (\ref{ellL2H1vt}).

\begin{remark}\label{thetadependence}
The $\theta$-dependence follows from estimating the terms $(q_{\ell
t}, \tta_{it}^j w_{\ell,j}^i)_{L^2(\Omega)}$:
\begin{align*}
& \Big|(q_{\ell t}, \tta_{it}^j w_{\ell,j}^i)_{L^2(\Omega)}\Big| \le
\frac{\theta}{2} \|q_{\ell t}\|^2_{L^2(\Omega)} +
\frac{1}{2\theta} \|\tta_{it}^j\|^2_{L^\infty(\Omega)} \|w_{\ell,j}^i\|^2_{L^2(\Omega)} \\
\le&\ \frac{\theta}{2} \|q_{\ell t}\|^2_{L^2(\Omega)} +
\frac{C(M)}{\theta}\Big[\|\nabla w_\ell(0)\|^2_{L^2(\Omega)} + t
\int_0^t \|\nabla w_{\ell t}\|^2_{L^2(\Omega)}(s) ds\Big].
\end{align*}
\end{remark}
\vspace{.1 in}

By the Gronwall inequality, for $0\le t\le T$,
\begin{align}
&\ \|w_{\ell t}(t)\|^2_{L^2(\Omega)} + \|\nabla_0^2 h_{\ell
t}(t)\|^2_{L^2(\Gamma)} + \int_0^t \Big[\|\nabla w_{\ell
t}\|^2_{L^2(\Omega)}
+ \kappa \|w_{\ell t}\|^2_{H^2(\Gamma)} \nonumber \\
&\qquad\quad + \theta \|q_{\ell t}\|^2_{L^2(\Omega)}\Big] ds \le
C(\eps1,\theta) N_0(u_0,F) \label{thetaL2H1}
\end{align}
where
$$N_0(u_0,F):= \|u_0\|^2_{H^{2.5}(\Omega)} +\|u_0\|^2_{H^{4.5}(\Gamma)}
+ \|F_t\|^2_{L^2(0,T;H^1(\Omega)')} + \|F(0)\|^2_{H^{0.5}(\Omega)} +
1.$$ We can then infer that $w_\ell$ defined on $[0,T]$, and that
there is a subsequence, still denoted with the subscript $\ell$,
satisfying
\begin{subequations}\label{thetaweaklimits}
\begin{alignat}{2}
w_\ell &\rightharpoonup w_\theta &&\quad\text{in }L^2(0,T;\H1H2) \\
w_{\ell t} &\rightharpoonup w_{\theta t}&&\quad\text{in }L^2(0,T;\H1H2) \\
\nabla_0^2 h_\ell &\rightharpoonup \nabla_0^2 h_\theta &&\quad\text{in }L^2(0,T;L^2(\Gamma)) \\
\nabla_0^2 h_{\ell t} &\rightharpoonup \nabla_0^2 h_{\theta t}&&\quad\text{in }L^2(0,T;L^2(\Gamma)) \\
q_{\ell t} &\rightharpoonup q_{\theta t}&&\quad\text{in
}L^2(0,T;L^2(\Omega))
\end{alignat}
\end{subequations}
where $$q_\theta = \tq_0 - \frac{1}{\theta}\tta_i^j
w_{\theta,j}^i.$$ From the standard procedure for weak solutions, we
can now infer from these weak convergences and the definition of
$w_\ell$ that $w_{\ell tt}\in L^2(0,T;H^1(\Omega)')$. In turn,
$w_{\ell t}\in {\mathcal C}^0([0,T];H^1(\Omega)')$, $w_\ell\in
{\mathcal C}^0([0,T];L^2(\Omega))$ with $w_\theta(0)=u_0$,
$w_{\theta t}(0) = w_1$.

\noindent
Moreover, (\ref{thetaweaklimits}) implies that $w_\theta$ satisfies
\begin{subequations}\label{thetaweak2}
\begin{align}
&\text{(i) } \int_0^T \Big[(w_{\theta tt},\varphi)_{L^2(\Omega)} +
\nu (\tta_i^j w_{\theta t,j},\tta_i^k \varphi_{,k})_{L^2(\Omega)}
+ \nu ((\tta_i^j\tta_i^k)_t w_\theta, \varphi_{,k})_{L^2(\Omega)}\Big]dt \nonumber\\
&\ + \nu \int_0^T \Big[\int_{\Omega} \tta_r^j\tta^k_i w_{\theta
t,j}^i \varphi_{,k}^r dx + \nu\int_{\Omega} (\tta_r^j\tta^k_i)_t
w_{\theta,j}^i \varphi_{,k}^r dx\Big]dt
+ \sigma \int_0^T \int_{\Gamma} \ttA^{\alpha\beta\gamma\delta}\times \nonumber\\
&\qquad \times [-\tth_{,\sigma} (w_\theta^\sigma\circ\tteta^{-\tau}) +
w_\theta^z\circ\tteta^{-\tau}]^\eps1_{,\alpha\beta}[-\tth_{,\sigma} (\varphi^\sigma\circ\tteta^{-\tau}) +
\varphi^z\circ\tteta^{-\tau}]^\eps1_{,\gamma\delta}dS dt \nonumber\\
&\ + \sigma\int_0^T \int_{\Gamma} (\ttA^{\alpha\beta\gamma\delta})_t
h^\eps1_{\theta,\alpha\beta} [-\tth_{,\sigma}
(\varphi^\sigma\circ\tteta^{-\tau}) +
\varphi^z\circ\tteta^{-\tau}]^\eps1_{,\gamma\delta}dS dt\\
&\ + \sigma\int_0^T \int_{\Gamma} \ttA^{\alpha\beta\gamma\delta}
h^\eps1_{\theta,\alpha\beta} [-\tth_{t,\sigma}
(\varphi^\sigma\circ\tteta^{-\tau}) + \tth_{,\sigma} \ttv^\kappa
(\varphi^\sigma_{,\kappa}
\circ\tteta^{-\tau}) + \ttv^\kappa (\varphi^z_{,\kappa}\circ\tteta^{-\tau})]^\eps1_{,\gamma\delta}dS dt \nonumber\\
&\ + \kappa \int_0^T \int_{\Gamma} \Delta_0 w_{\theta t}\cdot
\Delta_0 \varphi dS dt
- \int_0^T ((\tta^j_i q_\theta)_t,\varphi^i_{,j})_{L^2(\Omega)} dt \nonumber\\
&\ = \int_0^T\Big\{ \langle \tF_t,\varphi\rangle - \sigma
\int_{\Gamma} \Big[L_1^{\alpha\beta\gamma}\tth_{,\alpha\beta\gamma}
+ L_2\Big]^\eps1_t \Big[\tth_{,\sigma}
(\varphi^\sigma\circ\tteta^{-\tau}) - \varphi^z\circ\tteta^{-\tau}\Big]^\eps1 dS \nonumber\\
&\qquad\qquad - \sigma\int_{\Gamma}
\Big[L_1^{\alpha\beta\gamma}\tth_{,\alpha\beta\gamma} +
L_2\Big]^\eps1 \Big[\tth_{t,\sigma}
(\varphi^\sigma\circ\tteta^{-\tau}) - \tth_{,\sigma} \ttv^\kappa
(\varphi^\sigma_{,\kappa}
\circ\tteta^{-\tau}) \nonumber\\
&\qquad\qquad\qquad\qquad - \ttv^\kappa (\varphi^z_{,\kappa}\circ\tteta^{-\tau}) \Big]^\eps1 dS\Big\}dt \nonumber\\
&\text{(ii) } \ w_{\theta t}(0)= \tilde{w}_1, w_\theta(0) = \tu_0
\text{ in $\Omega$}\,,
\end{align}
\end{subequations}
for all $\varphi\in L^2(0,T;\H1H2)$. Choosing $\varphi$ to be independent of time, we find that for all $t\in [0,T]$,
\begin{align*}
&\ (w_{\theta t},\varphi)_{L^2(\Omega)} + \frac{\nu}{2}
\int_{\Omega} D_\tteta(w_\theta):D_\tteta(\varphi) dx
+ \kappa \int_{\Gamma} \Delta_0 w_\theta\cdot \Delta_0\varphi dS \nonumber\\
& + \sigma \int_{\Gamma} \ttA^{\alpha\beta\gamma\delta}
h^\eps1_{\theta,\alpha\beta}[-\tth_{,\sigma}
(\varphi^\sigma\circ\tteta^{-\tau}) +
\varphi^z\circ\tteta^{-\tau}]^\eps1_{,\gamma\delta}dS
- (\tta^j_i q_\theta,\varphi^i_{,j})_{L^2(\Omega)} \\
=&\ \langle \tF,\varphi\rangle + \sigma \int_{\Gamma}
\Big[L_1^{\alpha\beta\gamma\delta} \tth_{,\alpha\beta\gamma} +
L_2\Big]^\eps1 \Big[-\tth_{,\sigma}\varphi^\sigma\circ\tteta^{-\tau}
+ \varphi^z\circ\tteta^{-\tau}\Big]^\eps1 dS + c(\varphi)
\end{align*}
for all $\varphi\in \H1H2$, where $c(\varphi)\in\bbR$ is given by
\begin{align*}
c(\varphi) =&\ (\tilde{w}_1,\varphi)_{L^2(\Omega)} +
\frac{\nu}{2}\int_{\Omega} \Def (\tu_0) :\Def\varphi dx
- (\tq_0 - \frac{1}{\theta}\div \tu_0,\div\varphi)_{L^2(\Omega)} \\
& - (\tF(0),\varphi)_{L^2(\Omega)} -
\sigma(\bG^\eps1_0(0)(0,1),\varphi)_{L^2(\Gamma)} + \kappa(\Delta_0
\tu_0 ,\Delta_0 \varphi)_{L^2(\Gamma)}\,.
\end{align*}
By compatibility conditions (\ref{tildecompatibility1}) and (\ref{tildecompatibility2}), $c(\varphi) = 0.$
Therefore, the weak limit ($w_\theta$, $h_\theta$) satisfies, for all $t\in [0,T]$,
\begin{align}
& (w_{\theta t},\varphi)_{L^2(\Omega)} + \frac{\nu}{2} \int_{\Omega}
D_\tteta(w_\theta):D_\tteta(\varphi) dx
+ \kappa \int_{\Gamma} \Delta_0 w_\theta\cdot \Delta_0\varphi dS \nonumber\\
& - (\tta^j_i q_\theta,\varphi^i_{,j})_{L^2(\Omega)} + \sigma
\int_{\Gamma} \ttA^{\alpha\beta\gamma\delta}
h^\eps1_{\theta,\alpha\beta}[-\tth_{,\sigma}
(\varphi^\sigma\circ\tteta^{-\tau}) + \varphi^z\circ\tteta^{-\tau}]^\eps1_{,\gamma\delta}dS \label{thetaweak1}\\
=&\ \langle \tF,\varphi\rangle - \sigma \int_{\Gamma}
\Big[L_1^{\alpha\beta\gamma\delta} \tth_{,\alpha\beta\gamma} +
L_2\Big]^\eps1 \Big[-\tth_{,\sigma}\varphi^\sigma\circ\tteta^{-\tau}
+ \varphi^z\circ\tteta^{-\tau}\Big]^\eps1 dS \,, \nonumber
\end{align}
for all $\varphi\in \H1H2$.

Since $w_\theta\in L^2(0,T;\H1H2)$, we can use it as a test function in (\ref{thetaweak1}) and obtain (after time integration)
\begin{align}
&\ \ \frac{1}{2}\|w_\theta\|^2_{L^2(\Omega)} + \frac{\sigma}{2}
E_\tth(h^\eps1_\theta) + \int_0^t \Big[\frac{\nu}{2} \|D_\tteta
w_\theta\|^2_{L^2(\Omega)} + \kappa\|\Delta_0
w_\theta\|^2_{L^2(\Gamma)}
\nonumber\\
&\ \ + \theta \|q_\theta\|^2_{L^2(\Omega)} \Big]ds - \theta \int_0^t
(q_\theta, \tq_0) dt - \frac{\sigma}{2} \int_0^t \int_{\Gamma}
(\ttA^{\alpha\beta\gamma\delta})_t
h^\eps1_{\theta,\alpha\beta} h^\eps1_{\theta,\gamma\delta} dS ds \label{thetaeq1}\\
&= \frac{1}{2}\|\tu_0\|^2_{L^2(\Omega)} + \int_0^t \langle
\tF,\varphi\rangle + \sigma\langle
\bG_\tth^\eps1(-\nabla_0\tth\circ\tteta^\tau,1) ,\varphi
\rangle_{\Gamma} dt. \nonumber
\end{align}
Consequently,
\begin{align*}
&\ \Big[\|w_\theta(t)\|^2_{L^2(\Omega)} + \|\nabla_0^2
h^\eps1_\theta(t)\|^2_{L^2(\Gamma)}\Big] + \int_0^t \|\nabla
w_\theta\|^2_{L^2(\Omega)} ds
+ \kappa\int_0^t \|w_\theta\|^2_{H^2(\Gamma)} ds \\
& + \theta\int_0^t \|q_\theta\|^2_{L^2(\Omega)} ds \\
\le&\ C(M)\Big[\|\tu_0\|^2_{L^2(\Omega)} + \theta
\|\tq_0\|^2_{L^2(\Omega)} + \|\tF\|^2_{H^1(\Omega)'} +
\|\bG^\eps1_\tth(-\nabla_0\tth\circ\tteta^\tau,1)\|^2_{L^2(\Gamma)} \Big] \\
& + C(M) \int_0^t \|\tth_t\|_{H^{2.5}(\Gamma)} \|\nabla_0^2 h^\eps1_\theta\|^2_{L^2(\Gamma)} ds \\
\le&\ C(M)\Big[N_1(u_0,F) + \int_0^t \|\tth_t\|_{H^{2.5}(\Gamma)}
\|\nabla_0^2 h^\eps1_\theta\|^2_{L^2(\Gamma)} ds\Big]
\end{align*}
where
\begin{align*}
N_1(u_0,F)=&\ \|u_0\|^2_{H^2(\Omega)} + \|u_0\|^2_{H^{4.5}(\Gamma)}
+ \|F\|^2_{L^2(0,T;H^1(\Omega)')} +
\|F_t\|^2_{L^2(0,T;H^1(\Omega)')} \\
& + \|F(0)\|^2_{H^1(\Omega)} + 1.
\end{align*}
By the Gronwall inequality,
\begin{align}
& \sup_{0\le t\le T} \Big[\|w_\theta(t)\|^2_{L^2(\Omega)} +
\|\nabla_0^2 h^\eps1_\theta(t)\|^2_{L^2(\Gamma)}\Big]
+ \int_0^T \Big[\|\nabla w_\theta\|^2_{L^2(\Omega)} + \theta \|q_\theta\|^2_{L^2(\Omega)}\Big]ds \nonumber\\
\le&\ C(M)N_1(u_0,F). \label{thetaL2H1ineq}
\end{align}

\subsection{Improved pressure estimates}
By $\eps1$-regularization, we can rewrite (\ref{thetaweak1}) as, for a.a. $t\in[0,T]$,
\begin{align*}
&\ (w_{\theta t},\varphi)_{L^2(\Omega)} + \frac{\nu}{2}
\int_{\Omega} D_\tteta(w_\theta):D_\tteta(\varphi) dx + \kappa
(\Delta_0 w_\theta,\Delta_0\varphi)_{L^2(\Gamma)}
- (\tta^j_i q_\theta,\varphi^i_{,j})_{L^2(\Omega)} \\
&\ + \sigma \int_{\Gamma} \bL^\eps1_\tth(h^\eps1_\theta)
\Big[-\tth_{,\sigma}\circ\tteta^\tau \varphi^\sigma + \varphi^z\Big]
dS = \langle \tF,\varphi\rangle + \sigma\langle
\bG_\tth^\eps1(-\nabla_0\tth\circ\tteta^\tau,1) ,\varphi
\rangle_{\Gamma}.
\end{align*}
Therefore, by the Lagrange multiplier lemma, we conclude that
\begin{align*}
\|q_\theta\|^2_{L^2(\Omega)} \le&\ C(M)\Big[\|w_{\theta
t}\|^2_{H^1(\Omega)'} + \|\nabla w_\theta\|^2_{L^2(\Omega)}
+ \|\tF\|^2_{H^1(\Omega)'} + \kappa\|\Delta_0^2 w_\theta\|^2_{H^{-2}(\Gamma)}\\
&\qquad\quad + \|[\bL^\eps1_\tth(h^\eps1_\theta)+\bG^\eps1_\tth]
(-\nabla_0\tth\circ\tteta^\tau,1)\|^2_{H^{-2}(\Gamma)}\Big]
\end{align*}
and hence
\begin{align}
\|q_\theta\|^2_{L^2(\Omega)} \le&\ C(M)\Big[\|w_{\theta
t}\|^2_{L^2(\Omega)} + \|\nabla w_\theta\|^2_{L^2(\Omega)} + \kappa
\|w_\theta\|^2_{H^2(\Gamma)} + \|\nabla_0^2
h_\theta\|^2_{L^2(\Gamma)}
\nonumber\\
&\qquad\quad + \|F\|^2_{H^1(\Omega)'} + 1\Big]. \label{eps1lagrange}
\end{align}

\subsection{Weak limits as $\theta\rightarrow 0$}
Since $w_{\theta t} \in L^2(0,T;\H1H2)$, we can use it as a test function in (\ref{thetaweak2}). Similar to the way we
obtain (\ref{ellL2H1vt}), we find that
\begin{align*}
&\ \frac{1}{2}\|w_{\theta t}\|^2_{L^2(\Omega)} +
\frac{\nu}{2}\int_0^t |D_\tteta w_{\theta t}\|^2_{L^2(\Omega)} ds
+ \frac{\sigma}{2} E_\tth(h^\eps1_{\theta t}) + \kappa\int_0^t \|\Delta_0^2 w_{\theta t}\|^2_{L^2(\Gamma)} ds \\
&+ \theta\int_0^t \|q_{\theta t}\|^2_{L^2(\Omega)} ds + \int_0^t
(q_{\theta t}, \tta_{it}^j w_{\theta,j}^i)_{L^2(\Omega)} ds
-\int_0^t (q_\theta, \tta_i^j w_{\theta t,j}^i) ds \\
&\le C(M)N_0(u_0,F) + C(M)\int_0^t \|\ttv(t')\|^2_{H^3(\Omega)}
\int_0^{t'}
\|\nabla w_{\theta t}(s)\|^2_{L^2(\Omega)} ds dt' \\
&+ C(\eps1) \int_0^t \Big[1+\|\tth_t\|_{H^{2.5}(\Gamma)}\Big]
\|\nabla_0^2 h^\eps1_{\theta t}\|^2_{L^2(\Gamma)} ds.
\end{align*}
By (\ref{eps1lagrange}),
\begin{align}
&\ \Big|\int_0^t (q_\theta, \tta_i^j w_{\theta t,j}^i) ds\Big| \le
C(M,\delta)\int_0^t \|q_\theta\|^2_{L^2(\Omega)} ds
+ \delta\int_0^t \|\nabla w_{\theta t}\|^2_{L^2(\Omega)} ds \nonumber\\
%\le&\ C(M)\int_0^t \Big[\|w_{\theta t}\|^2_{L^2(\Omega)} + \|\nabla w_\theta\|^2_{L^2(\Omega)}
%+ \kappa \|w_\theta\|^2_{H^2(\Gamma)} + \|\nabla_0^2 h_\theta\|^2_{L^2(\Gamma)} \Big] ds \nonumber\\
%& + \delta\int_0^t \|\nabla w_{\theta t}\|^2_{L^2(\Omega)} ds + C(M)N_0(u_0,F)\nonumber\\
\le&\ C(M)\Big[N_1(u_0,F) + \int_0^t \Big(\|w_{\theta
t}\|^2_{L^2(\Omega)} + \kappa \|w_\theta\|^2_{H^2(\Gamma)}
+ \|\nabla_0^2 h_\theta\|^2_{L^2(\Gamma)}\Big) ds\Big] \nonumber\\
& + \delta\int_0^t \|\nabla w_{\theta t}\|^2_{L^2(\Omega)} ds
\label{qtheta1}
\end{align}
where (\ref{thetaL2H1ineq}) is used to bound $\|\nabla
w_\theta\|^2_{L^2(0,T;L^2(\Omega))}$.

\noindent
Integrating by parts,
\begin{align*}
& \int_0^t (q_{\theta t}, \tta_{it}^j w_{\theta,j}^i)_{L^2(\Omega)}
ds = (q_\theta, \tta_{it}^j
w_{\theta,j}^i)_{L^2(\Omega)}(t) + (\tq_0,\tu_{0,i}^j \tu_{0,j}^i)_{L^2(\Omega)} \\
&\qquad\quad - \int_0^t (q_\theta, \tta_{itt}^j
w_{\theta,j}^i)_{L^2(\Omega)} ds - \int_0^t (q_\theta, \tta_{it}^j
w_{\theta t,j}^i)_{L^2(\Omega)} ds.
\end{align*}
By $\epsilon$-regularization, the last two term can be bounded by
\begin{align*}
C(M)\int_0^t \|q_\theta\|_{L^2(\Omega)}\Big[C(\epsilon)\|\nabla
w_\theta\|_{L^2(\Omega)} + \|\nabla w_{\theta
t}\|_{L^2(\Omega)}\Big] ds
\end{align*}
and hence
\begin{align}
& \Big|\int_0^t (q_\theta, \tta_{itt}^j
w_{\theta,j}^i)_{L^2(\Omega)} ds\Big|
+ \Big|\int_0^t (q_\theta, \tta_{it}^j w_{\theta t,j}^i)_{L^2(\Omega)} ds\Big| \nonumber\\
\le&\ C(M,\delta) \int_0^t \|q_\theta\|^2_{L^2(\Omega)} ds +
C(\epsilon)\int_0^t \|\nabla w_\theta\|^2_{L^2(\Omega)} ds
+ \delta\int_0^t \|\nabla w_{\theta t}\|^2_{L^2(\Omega)} ds \nonumber\\
\le&\ C(\epsilon,\delta)N_1(u_0,F) + C(M,\delta)\int_0^t \|w_{\theta
t}\|^2_{L^2(\Omega)} ds + C(\eps1) \int_0^t \|\nabla_0^2
h_\theta\|^2_{L^2(\Gamma)} ds \nonumber\\
& + \delta\int_0^t \|\nabla w_{\theta t}\|^2_{L^2(\Omega)} ds.
\label{qtheta2}
\end{align}

\noindent For $(q_\theta, \tta_{it}^j
w_{\theta,j}^i)_{L^2(\Omega)}(t)$, it is easy to see that
\begin{align*}
&\ \Big|(q_\theta, \tta_{it}^j w_{\theta,j}^i)_{L^2(\Omega)}(t)\Big|
\le \delta_1 \|w_{\theta t}\|^2_{L^2(\Omega)} +
C(\epsilon,\delta_1)\|\nabla w_\theta\|^2_{L^2(\Omega)} \\
\le&\  C(\epsilon,\delta_1)\|\nabla w_\theta\|^2_{L^2(\Omega)} +
\delta_1 C(\eps1) \|\nabla_0^2 h_\theta\|^2_{L^2(\Gamma)} +
\delta_1\Big[\|w_{\theta t}\|^2_{L^2(\Omega)} + \|F\|_{L^2(\Omega)}
+ 1\Big]
\end{align*}
while for $(\tq_0,\tu_{0,i}^j \tu_{0,j}^i)_{L^2(\Omega)}$, it is
bounded by $C(M) N_1(u_0,F)$. Combining (\ref{qtheta1}),
(\ref{qtheta2}) and the estimates above, by choosing $\delta>0$ and
$\delta_1>0$ small enough,
\begin{align*}
&\ \|w_{\theta t}\|^2_{L^2(\Omega)} + \|\nabla_0^2 h_{\theta
t}\|^2_{L^2(\Gamma)} +
\int_0^t \Big[\|\nabla w_{\theta t}\|^2_{L^2(\Omega)} + \kappa \|w_{\theta t}\|^2_{H^2(\Gamma)} + \theta\|q_{\theta t}\|^2_{L^2(\Omega)}\Big] ds \\
\le&\ C(\eps1,\epsilon)\Big[N_2(u_0,F) +
\int_0^t \Big(\|w_{\theta t}\|^2_{L^2(\Omega)} + (1+\|\tth_t\|_{H^{2.5}(\Gamma)}) \|\nabla_0^2 h_{\theta t}\|^2_{L^2(\Gamma)} \\
& + \|\ttv\|^2_{H^3(\Omega)} \int_0^{s} \|\nabla w_{\theta
t}\|^2_{L^2(\Omega)}dt' \Big) ds\Big] + C_1(\eps1,\epsilon) \|\nabla
w_\theta\|^2_{L^2(\Omega)}
\end{align*}
where $N_2(u_0,F) = N_1(u_0,F) +
\|F\|^2_{L^\infty(0,T;L^2(\Omega))}$. By the Gronwall inequality,
\begin{align}
&\ \|w_{\theta t}\|^2_{L^2(\Omega)} + \|\nabla_0^2 h_{\theta
t}\|^2_{L^2(\Gamma)} +
\int_0^t \Big[\|\nabla w_{\theta t}\|^2_{L^2(\Omega)} + \kappa \|w_{\theta t}\|^2_{H^2(\Gamma)}\Big] ds \nonumber\\
\le&\ C(\eps1,\epsilon)N_2(u_0,F) + C_1(\eps1,\epsilon) \|\nabla
w_\theta\|^2_{L^2(\Omega)}. \label{thetaL2H1vt2}
\end{align}
By using $\displaystyle{w_\theta(t) = \tu_0 + \int_0^t w_{\theta t} ds}$, we see that
\begin{align*}
&\ \|w_{\theta t}\|^2_{L^2(\Omega)} + \|\nabla_0^2 h_{\theta
t}\|^2_{L^2(\Gamma)} +
\int_0^t \Big[\|\nabla w_{\theta t}\|^2_{L^2(\Omega)} + \kappa \|w_{\theta t}\|^2_{H^2(\Gamma)}\Big] ds \\
\le&\ C(\eps1,\epsilon)N_2(u_0,F) + C_1(\eps1,\epsilon) t\int_0^t
\|\nabla w_{\theta t}\|^2_{L^2(\Omega)} ds.
\end{align*}
Therefore, for any $\displaystyle{0\le t\le t_1 = \min\Big\{T,\frac{1}{2C_1}\Big\}}$, we have
\begin{align*}
&\ \|w_{\theta t}\|^2_{L^2(\Omega)} + \|\nabla_0^2 h_{\theta
t}\|^2_{L^2(\Gamma)} + \frac{1}{2}
\int_0^t \Big[\|\nabla w_{\theta t}\|^2_{L^2(\Omega)} + \kappa \|w_{\theta t}\|^2_{H^2(\Gamma)}\Big] ds \\
\le&\ C(\eps1,\epsilon)N_2(u_0,F).
\end{align*}
By $\displaystyle{w_\theta(t_1) = \tu_0 + \int_0^{t_1} w_{\theta t} ds}$, we also have
\begin{align}
\|\nabla w_\theta(t_1)\|^2_{L^2(\Omega)} \le
C(\eps1,\epsilon)N_2(u_0,F). \label{thetat1}
\end{align}
For $t\ge t_1$, since $\displaystyle{w_\theta(t) = w_\theta(t_1) + \int_{t_1}^t w_{\theta t} ds}$, we have from
(\ref{thetaL2H1vt2}) and (\ref{thetat1}) that
\begin{align*}
&\ \|w_{\theta t}\|^2_{L^2(\Omega)} + \|\nabla_0^2 h_{\theta
t}\|^2_{L^2(\Gamma)} +
\int_0^t \Big[\|\nabla w_{\theta t}\|^2_{L^2(\Omega)} + \kappa \|w_{\theta t}\|^2_{H^2(\Gamma)}\Big] ds \\
\le&\ C(\eps1,\epsilon)N_2(u_0,F) +
C_1(\eps1,\epsilon)\Big[\|w_\theta(t_1)\|^2_{L^2(\Omega)} +
(t-t_1)\int_{t_1}^t \|\nabla_0 w_{\theta t}\|^2_{L^2(\Omega)} ds\Big] \\
\le&\ C(\eps1,\epsilon)N_2(u_0,F) + C_1(\eps1,\epsilon)
(t-t_1)\int_{t_1}^t \|\nabla_0 w_{\theta t}\|^2_{L^2(\Omega)} ds\Big]. \\
\end{align*}
Therefore, for any $t_1\le t\le 2t_1$, we also have
\begin{align*}
&\ \|w_{\theta t}\|^2_{L^2(\Omega)} + \|\nabla_0^2 h_{\theta
t}\|^2_{L^2(\Gamma)} + \frac{1}{2}
\int_0^t \Big[\|\nabla w_{\theta t}\|^2_{L^2(\Omega)} + \kappa \|w_{\theta t}\|^2_{H^2(\Gamma)}\Big] ds\\
\le&\ C(\eps1,\epsilon)N_2(u_0,F)
\end{align*}
which with $\displaystyle{w_\theta(2t_1) = \tu_0 + \int_0^{2t_1} w_{\theta t} ds}$ gives
\begin{align*}
\|\nabla w_\theta(2t_1)\|^2_{L^2(\Omega)} \le
C(\eps1,\epsilon)N_2(u_0,F).
\end{align*}
By induction, for any $t\in [0,T]$,
\begin{align}
&\ \|w_{\theta t}\|^2_{L^2(\Omega)} + \|\nabla_0^2 h_{\theta
t}\|^2_{L^2(\Gamma)} + \frac{1}{2}
\int_0^t \Big[\|\nabla w_{\theta t}\|^2_{L^2(\Omega)} + \kappa \|w_{\theta t}\|^2_{H^2(\Gamma)}\Big] ds \nonumber\\
\le& C(\eps1,\epsilon)N_2(u_0,F). \label{thetaL2H1vt}
\end{align}
We also get a $\theta$-independent bound for
$\|q_\theta\|^2_{L^2(0,T;L^2(\Omega))}$ by (\ref{eps1lagrange}):
\begin{align}
\|q_\theta\|^2_{L^2(0,T;L^2(\Omega))} \le
C(\eps1,\epsilon)N_2(u_0,F). \label{qthetaestimate}
\end{align}

Let $\theta = \frac{1}{m}$. Energy inequalities ({\ref{thetaL2H1ineq}), (\ref{thetaL2H1vt}) and (\ref{qthetaestimate})
show that there exists a
subsequence $w_{\frac{1}{m_\ell}}$ such that
\begin{subequations}\label{thetaconvergence}
\begin{alignat}{2}
w_{\frac{1}{m_\ell}} &\rightharpoonup \mv &&\qquad\text{in}\quad L^2(0,T;\H1H2) \\
w_{\frac{1}{m_\ell}t} &\rightharpoonup \mv_t &&\qquad\text{in}\quad L^2(0,T;\H1H2) \\
\nabla_0^2 h_{\frac{1}{m_\ell}} &\rightharpoonup \nabla_0^2 \mh &&\qquad\text{in}\quad L^2(0,T;L^2(\Omega)) \\
\nabla_0^2 h_{\frac{1}{m_\ell}t} &\rightharpoonup \nabla_0^2 \mh_t &&\qquad\text{in}\quad L^2(0,T;L^2(\Omega)) \\
q_{\frac{1}{m_\ell}} &\rightharpoonup \mq &&\qquad\text{in}\quad
L^2(0,T;L^2(\Omega))\,.
\end{alignat}
\end{subequations}
Moreover, (\ref{thetaL2H1ineq}) also shows that $\|\tta_i^j
w_{\frac{1}{m},j}^i\|_{L^2(0,T;L^2(\Omega))} \to 0$ as $m\to
\infty$. Therefore the weak limit ${\mathfrak v}$ satisfies the
``divergence-free'' condition (\ref{NSequation2}.b), i.e.,
\begin{align}
{\mathfrak v}\in {\mathcal V}_\ttv(T). \label{divfree}
\end{align}
Since (\ref{thetaL2H1ineq}) is independent of $\theta$ and $\eps1$, by the property of lower-semicontinuity of norms,
\begin{align}
& \sup_{0\le t\le T}\Big[\|\mv(t)\|^2_{L^2(\Omega)} + \|\nabla_0^2
\mh(t)\|^2_{L^2(\Gamma)}\Big]
+ \|\nabla \mv\|_{L^2(0,T;L^2(\Omega))}^2 + \kappa \|\mv \|^2_{H^2(\Gamma)} \nonumber\\
\le&\ C(M)N_1(u_0,F). \label{L2H1estimate}
\end{align}

By (\ref{thetaconvergence}) and $\eps1$-regularization, the weak limit $(\mv,\mh,\mq)$ satisfies,
for all $\varphi\in L^2(0,T;\H1H2)$,
\begin{align*}
&\ \int_0^T (\mv_t,\varphi)_{L^2(\Omega)}dt + \frac{\nu}{2}
\int_0^T\int_{\Omega} D_\tteta(\mv):D_\tteta(\varphi) dx dt
+ \kappa \int_0^T\int_{\Gamma} \Delta_0 \mv \cdot \Delta_0\varphi dS dt \\
&\ - \int_0^T (\tta^j_i \mq, \varphi^i_{,j})_{L^2(\Omega)} dt +
\sigma \int_0^T \int_{\Gamma} \ttA^{\alpha\beta\gamma\delta}
\mh^\eps1_{,\alpha\beta}[-\tth_{,\sigma}
(\varphi^\sigma\circ\tteta^{-\tau}) + \varphi^z\circ\tteta^{-\tau}]^\eps1_{,\gamma\delta}dS dt \\
&= \int_0^T \Big\{\langle \tF,\varphi\rangle - \sigma \int_{\Gamma}
\Big[L_1^{\alpha\beta\gamma\delta} \tth_{,\alpha\beta\gamma} +
L_2\Big]^\eps1 \Big[-\tth_{,\sigma}\varphi^\sigma\circ\tteta^{-\tau}
+ \varphi^z\circ\tteta^{-\tau}\Big]^\eps1 dS \Big\}dt\,.
\end{align*}
By the density argument, we find that for a.a. $t\in[0,T]$, $\varphi\in \H1H2$,
\begin{align}
&\ (\mv_t,\varphi)_{L^2(\Omega)} + \frac{\nu}{2} \int_{\Omega}
D_\tteta(\mv):D_\tteta(\varphi) dx + \kappa \int_{\Gamma} \Delta_0
\mv\cdot \Delta_0 \varphi dS
- (\tta^j_i \mq, \varphi^i_{,j})_{L^2(\Omega)} \nonumber\\
&\ + \sigma \int_{\Gamma} \ttA^{\alpha\beta\gamma\delta}
\mh^\eps1_{,\alpha\beta}[-\tth_{,\sigma}
(\varphi^\sigma\circ\tteta^{-\tau}) +
\varphi^z\circ\tteta^{-\tau}]^\eps1_{,\gamma\delta}dS
\label{eps1weakform}\\
&= \langle \tF,\varphi\rangle - \sigma \int_{\Gamma}
\Big[L_1^{\alpha\beta\gamma\delta} \tth_{,\alpha\beta\gamma} +
L_2\Big]^\eps1 \Big[-\tth_{,\sigma}\varphi^\sigma\circ\tteta^{-\tau}
+ \varphi^z\circ\tteta^{-\tau}\Big]^\eps1 dS \nonumber\,,
\end{align}
or after a change of variable $y'=\tteta^\tau(y,t)$,
\begin{align}
&\ (\mv_t,\varphi)_{L^2(\Omega)} + \frac{\nu}{2}
(D_\tteta\mv,D_\tteta \varphi)_{L^2(\Omega)} + \kappa \int_{\Gamma}
\Delta_0 \mv\cdot \Delta_0 \varphi dS
- (\tta^j_i \mq, \varphi^i_{,j})_{L^2(\Omega)} \label{eps1weakform1}\\
&+ \sigma \int_{\Gamma} \L_\tth^\eps1(\mh)
(-\nabla_0\tth\circ\tteta^\tau,1)\cdot\varphi dS = \langle
\tF,\varphi\rangle - \sigma\int_{\Gamma} \bG_\tth^\eps1
(-\nabla_0\tth\circ\tteta^\tau,1)\cdot\varphi dS\,. \nonumber
\end{align}
Furthermore, if $\varphi\in {\mathcal V}_\ttv$, then
\begin{align*}
&\ (\mv_t,\varphi)_{L^2(\Omega)} + \frac{\nu}{2}
(D_\tteta\mv,D_\tteta \varphi)_{L^2(\Omega)}
+ \kappa \int_{\Gamma} \Delta_0 \mv\cdot \Delta_0 \varphi dS \nonumber\\
& + \sigma \int_{\Gamma} \L_\tth^\eps1(\mh)
(-\nabla_0\tth\circ\tteta^\tau,1)\cdot\varphi dS = \langle
\tF,\varphi\rangle - \sigma\int_{\Gamma} \bG_\tth^\eps1
(-\nabla_0\tth\circ\tteta^\tau,1)\cdot\varphi^\eps1 dS
\end{align*}
for a.a. $t\in[0,T]$. In other words, $(\mv,\mh,\mq)$ is a weak solution of (\ref{NSequation2}).

\section{Estimates independent of $\eps1$}\label{energyestimatesection1}
\subsection{Partition of unity} Since $\Omega$ is compact, by partition of unity, we can choose two non-negative smooth
functions $\zeta_0$ and $\zeta_1$ so that
\begin{align*}
\zeta_0 + \zeta_1 &= 1 \quad\text{in}\quad \Omega \ ;\\
\supp(\zeta_0) &\subset\subset \Omega \ ;\\
\supp(\zeta_1) &\subset\subset
\Gamma\times(-\epsilon,\epsilon):=\Omega_1.
\end{align*}
We will assume that $\zeta_1=1$ inside the region $\Omega_1'\subset
\Omega_1$ and $\zeta_0=1$ inside the region $\Omega'\subset \Omega$.
Note that then $\zeta_1 = 1$ while $\zeta_0 = 0$ on $\Gamma$.

\subsection{Higher regularity}
\subsubsection{$\eps1$-independent bounds for $\mq$}
Similar to (\ref{eps1lagrange}), we have
\begin{align}
\|\mq\|^2_{L^2(\Omega)} \le&\ C(M)\Big[\|\mv_t\|^2_{L^2(\Omega)} +
\|\nabla \mv\|^2_{L^2(\Omega)} + \kappa \|\mv\|^2_{H^2(\Gamma)}
+ \|\nabla_0^2 \mh^\eps1\|^2_{L^2(\Gamma)} \nonumber\\
&\qquad\quad + \|F\|^2_{L^2(\Omega)} + 1\Big]. \label{lagrange1}
\end{align}
\subsubsection{Interior regularity} Converting the fluid equation (\ref{NSequation2}) into Eulerian variables by composing
with $\tteta^{-1}$, we obtain a Stokes problem in the domain
$\tteta(\Omega)$:
\begin{subequations}\label{stokes}
\begin{align}
-\nu\Delta {\mathfrak u} + \nabla {\mathfrak p} &= \tF\circ\tteta^{-1} - \mv_t\circ\tteta^{-1} +
\nu\tta_{\ell,j}^j\circ\tteta^{-1} {\mathfrak u}_{,\ell} - {\mathfrak p}\tta^j_{i,j}\circ\tteta^{-1}\,, \\
\div {\mathfrak u} &= 0\,,
\end{align}
\end{subequations}
where ${\mathfrak u} = \mv\circ \tteta^{-1}$ and ${\mathfrak p} = \mq\circ \tteta^{-1}$. By the regularity results for the Stokes problem,
\begin{align*}
&\ \|{\mathfrak u}\|^2_{H^2(\tteta(\Omega))} + \|{\mathfrak p}\|^2_{H^1(\tteta(\Omega))} \\
\le&\ C\Big[\|\tF\circ\tteta^{-1}\|^2_{L^2(\tteta(\Omega))} +
\|\mv_t\circ\tteta^{-1}\|^2_{L^2(\tteta(\Omega))}
+ \|\nabla {\mathfrak u}\|^2_{L^2(\tteta(\Omega))} + \|{\mathfrak p}\|^2_{L^2(\tteta(\Omega))} \\
&\quad + \|{\mathfrak u}\|^2_{H^{1.5}(\Gamma)} \Big]
\end{align*}
or
\begin{align*}
\|\mv\|^2_{H^2(\Omega)} + \|\mq\|^2_{H^1(\Omega)}
\le&\ C\Big[\|F\|^2_{L^2(\Omega)} + \|\mv_t\|^2_{L^2(\Omega)} + \|\mv\|^2_{H^{1.5}(\Gamma)}\Big] \\
& + C(M)\Big[\|\nabla \mv\|^2_{L^2(\Omega)} +
\|\mq\|^2_{L^2(\Omega)} \Big]
\end{align*}
for some constant $C$ independent of $M$, $\epsilon$. By (\ref{lagrange1}),
\begin{align}
\|\mv\|^2_{H^2(\Omega)} + \|\mq\|^2_{H^1(\Omega)} \le&\
C(M)\Big[\|\mv_t\|^2_{L^2(\Omega)} + \|\nabla \mv\|^2_{L^2(\Omega)}
+ \|\mv\|^2_{H^2(\Gamma)} \nonumber\\
&\qquad\quad + \|\nabla_0^2 \mh^\eps1\|^2_{L^2(\Gamma)} +
\|F\|^2_{L^2(\Omega)} + 1\Big] \label{regularity1temp}
\end{align}
Similarly,
\begin{align*}
\|\mv\|^2_{H^3(\Omega)} + \|\mq\|^2_{H^2(\Omega)} \le&\
C\Big[ \|F\|^2_{H^1(\Omega)} + \|\mv_t\|^2_{H^1(\Omega)} + \|\mv\|^2_{H^{2.5}(\Gamma)}\Big] \\
& + C(M)\Big[\|\nabla \mv\|^2_{H^1(\Omega)} +
\|\mq\|^2_{H^1(\Omega)} \Big]
\end{align*}
and therefore by (\ref{lagrange1}) and (\ref{regularity1temp}),
\begin{align}
\|\mv\|^2_{H^3(\Omega)} + \|\mq\|^2_{H^2(\Omega)} \le&\
C(M)\Big[\|\mv_t\|^2_{H^1(\Omega)}
+ \|\nabla \mv\|^2_{L^2(\Omega)} + \|\nabla_0^2 \mv\|^2_{H^1(\Omega_1)} \nonumber\\
&\qquad\quad + \|\nabla_0^2 \mh^\eps1\|^2_{L^2(\Gamma)} +
\|F\|^2_{H^1(\Omega)} + 1\Big]. \label{regularity2temp}
\end{align}

For the regularized problem, because
the $\epsilon$-regularization ensures that the forcing and the initial data are smooth, while
the $\eps1$-regularization ensures that the right-hand side of (\ref{NSequation2.c}) is smooth, by standard
difference quotient technique, it is also easy to see that
\begin{align}
\nabla_0^k \mv\in L^2(0,T;H^1(\Omega_1)\cap H^2(\Gamma))
\quad\text{for $k=1,2,3,4$} \label{bdyhigherregularity}
\end{align}
Since (\ref{thetaconvergence}b) implies that $\mv_t\in
L^2(0,T;H^1(\Omega))$, by $\eps1$-regularization and
(\ref{regularity2temp}) we conclude that
\begin{align}
\mv\in L^2(0,T;H^3(\Omega)),\quad \mq\in L^2(0,T;H^2(\Omega)).
\label{eps1higherregularity}
\end{align}

\subsection{Estimates for $\mv_t(0)$ and $\mq(0)$} By (\ref{eps1higherregularity}) and $\eps1$-regularization,
$(\mv,\mh,\mq)$ satisfies the strong from (\ref{NSequation2}).\
Taking the ``divergence'' of (\ref{NSequation2.a}) and then making use of condition (\ref{NSequation2}.b), we find that
\begin{align}
-\tta_{it}^k \mv_{,k}^i - \nu \tta_i^k[\tta_\ell^j D_\tteta(\mv)_\ell^i]_{,jk}
= - \tta_i^k(\tta_i^j \mq)_{,jk} + \tta_i^k \tF^i_{,k}. \label{qequation}
\end{align}
Let $t=0$, by the identity $\tta_{kt}^\ell = -\tta_k^i \ttv_{,i}^j \tta_j^\ell$,
\begin{align*}
\Delta \mq(0) = \nabla \tu_0 : (\nabla \tu_0)^T - \div (\tF(0))
\qquad\text{in}\quad\Omega
\end{align*}
with
\begin{align*}
\mq(0) = \nu (\Def \tu_0)_i^j N_i N_j - \sigma \G_0^\eps1(0) +
\kappa \Delta_0^2 \tu_0 \qquad\text{on}\quad\Gamma
\end{align*}
while (\ref{NSequation2.a}) gives us
\begin{align*}
\mv_t(0) = \nu \Delta \tu_0 - \nabla \mq(0) + \tF(0)
\qquad\text{in}\quad\Omega.
\end{align*}
By standard elliptic regularity result,
\begin{align}
\|\mv_t(0)\|^2_{L^2(\Omega)} + \|\mq(0)\|^2_{H^1(\Omega)} \le
CN_0(u_0,F) \label{zerovalue}
\end{align}
for some constant independent of $M$, $\epsilon$ and $\eps1$.

\subsection{$L^2_tL^2_x$-estimates for $\mv_t$} Since $\mv_t\in L^2(0,T;H^1(\Omega))$, we can use it as a test function in
(\ref{eps1weakform1}). By (\ref{divfree}), we find that
\begin{align*}
&\ \ \|\mv_t\|^2_{L^2(\Omega)} +
\frac{\nu}{4}\frac{d}{dt}\int_{\Omega} |D_\tteta \mv|^2 dx -
\frac{\nu}{2}\int_{\Omega} (D_\tteta \mv)_i^j\tta_{jt}^k \mv_{,k}^i
dx
+ \kappa \int_{\Gamma} \Delta_0 \mv\cdot\Delta_0 \varphi dS \\
&\ + \int_{\Omega} \mq \tta_{kt}^\ell \mv_{,\ell}^k dx
+ \sigma \int_{\Gamma} \L_\tth^\eps1(\mh) (-\nabla_0\tth\circ\tteta^\tau,1)\cdot\mv_t dS \\
&= \langle \tF, \mv_t\rangle - \sigma\int_{\Gamma} \G_\tth^\eps1
(-\nabla_0\tth\circ\tteta^\tau,1)\cdot\mv_t dS\,.
\end{align*}
By (\ref{aestimate1}),
\begin{align*}
\int_{\Omega} (D_\tteta \mv)_i^j \tta_{jt}^k \mv^i_{,k} dx \le
C(M)C(\delta) \|\nabla \mv\|^2_{L^2(\Omega)} + \delta
\|\mv\|^2_{H^2(\Omega)}
\end{align*}
and by (\ref{lagrange1}) and the interpolation inequality,
\begin{align*}
\Big|\int_{\Omega} \mq\tta_{kt}^\ell \mv_{,\ell}^k dx\Big| \le&\
C(M)C(\delta)\Big[\|\nabla \mv\|^2_{L^2(\Omega)} + \|\nabla_0^4
\mh^\eps1\|^2_{L^2(\Gamma)} + \|F\|^2_{L^2(\Omega)} + 1\Big] \\
& + \delta \|\mv\|^2_{H^2(\Omega)} +
\frac{1}{2}\|\mv_t\|^2_{L^2(\Omega)}
\end{align*}
for some $C(\delta)$. Also, the last term on the left hand side is bounded by
\begin{align*}
&\ C(M)\Big[\|\nabla_0^4 \mh^\eps1\|_{L^2(\Gamma)} + 1\Big] \|\mv_t\|_{H^1(\Omega)} \\
\le&\ C(M)C(\delta_1)\Big[\|\nabla_0^4 \mh^\eps1\|^2_{L^2(\Gamma)} +
1\Big] + \delta_1 \|\mv_t\|^2_{H^1(\Omega)}.
\end{align*}
Combining all the estimates above,
\begin{align*}
&\ \frac{1}{2}\|\mv_t\|^2_{L^2(\Omega)} +
\frac{\nu}{4}\frac{d}{dt}\int_{\Omega} |D_\tteta \mv|^2 dx
+ \frac{\kappa}{2} \frac{d}{dt} \int_{\Gamma} |\Delta_0 \mv|^2 dS \\
\le&\ C\Big[\|\nabla \mv\|^2_{L^2(\Omega)} + \|\nabla_0^4
\mh^\eps1\|^2_{L^2(\Gamma)} + \|F\|^2_{L^2(\Omega)} + 1\Big] +
\delta \|\mv\|^2_{H^2(\Omega)} + \delta_1 \|\mv_t\|^2_{H^1(\Omega)}
\end{align*}
for some constant $C$ depending on $M$, $\delta$ and $\delta_1$. Therefore by (\ref{L2H1estimate}),
\begin{align}
&\ \int_0^t\|\mv_t\|^2_{L^2(\Omega)}ds + \|\nabla
\mv(t)\|^2_{L^2(\Omega)} + \kappa \|\mv\|^2_{H^2(\Gamma)}
\label{L2L2vtineq}\\
\le&\ C\Big[N_2(u_0,F) + \int_0^t \|\nabla_0^4
\mh^\eps1\|^2_{L^2(\Gamma)} ds\Big] + \delta \int_0^t
\|\mv\|^2_{H^2(\Omega)} ds + \delta_1\int_0^t
\|\mv_t\|^2_{H^1(\Omega)} ds. \nonumber
\end{align}

\subsection{Energy estimates for $\nabla_0^2 v$ near the boundary}
Because of (\ref{bdyhigherregularity}),
$\nabla_0^2(\zeta_1^2 \nabla_0^2 \mv)$ in (\ref{eps1weakform}) can be used as a test function in (\ref{eps1weakform1}).
It follows that
\begin{align*}
& \Big| \int_{\Gamma} \Big[\bL^\eps1_\tth(\mh^\eps1) +
\bG_\tth^\eps1\Big]
(-\nabla_0 \tth\circ\tteta^\tau,1)\cdot \nabla_0^4 \mv dS \Big| \\
\le&\ C(M) \Big[\|\nabla_0^2 \mh^\eps1\|_{H^2(\Gamma)} + 1\Big]\|\mv\|_{H^4(\Gamma)} \\
\le&\ C(M,\delta_3) \Big[1+ \|\mh\|^2_{H^4(\Gamma)}\Big] + \delta_3
\|\mv\|^2_{H^4(\Gamma)}.
\end{align*}
By (\ref{hevol}), we find that
\begin{align*}
\|\mh\|^2_{H^4(\Gamma)} \le C(\epsilon) \Big[\int_0^t
\|\tth\|_{H^5(\Gamma)} \|\mv\|_{H^4(\Gamma)} ds\Big]^2 \le
C(\epsilon) \int_0^t \|\mv\|^2_{H^4(\Gamma)} ds
\end{align*}
and hence
\begin{align*}
& \Big| \int_{\Gamma} \Big[\bL^\eps1_\tth(\mh^\eps1) +
\bG_\tth^\eps1\Big]
(-\nabla_0 \tth\circ\tteta^\tau,1)\cdot \nabla_0^4 \mv dS \Big| \\
\le&\ \bar{C} \Big[1+ \int_0^t \|\mv\|^2_{H^4(\Gamma)}\Big] +
\delta_3 \|\mv\|^2_{H^4(\Gamma)}.
\end{align*}
for some constant $\bar{C}$ depending on $M$, $\epsilon$ and $\delta_3$. Since
\begin{align*}
\Delta_0 f = \frac{1}{\sqrt{\det(g_0)}}\frac{\partial}{\partial y^\alpha}\Big[ \sqrt{\det(g_0)}
g_0^{\alpha\beta} \frac{\partial}{\partial y^\beta} f\Big],
\end{align*}
by the regularity on $\Gamma$ (and hence on $g_0$),
\begin{align*}
\int_{\Gamma} |\Delta_0 \nabla_0^2 \mv|^2 dS \le& \int_{\Gamma}
\Delta_0^2 \mv\cdot (\nabla_0^4 v) dS +
C\|\mv\|_{H^3(\Gamma)} \|\mv\|_{H^4(\Gamma)} \\
\le& \int_{\Gamma} \Delta_0^2 v\cdot (\nabla_0^4 v) dS +
C(\delta)\|\mv\|^2_{H^1(\Omega)} + \delta \|\mv\|^2_{H^4(\Gamma)}
\end{align*}
which implies, by choosing $\delta>0$ small enough, that
\begin{align*}
\nu_2 \|v\|^2_{H^4(\Gamma)} \le \int_{\Gamma} \Delta_0^2 v\cdot
(\nabla_0^4 v) dS + C\|v\|^2_{H^1(\Omega)}.
\end{align*}
By the identity
\begin{align}
&\ (\mq,\tta_k^\ell (\nabla_0^2 (\zeta_1^2 \nabla_0^2 \mv^k)_{,\ell}) \nonumber\\
=&\ (\mq, \nabla_0^2 \tta_k^\ell (\zeta_1^2 \nabla_0^2
\mv^k)_{,\ell}) + 4(\zeta_1\nabla_0 \mq,\nabla_0\tta_k^\ell
\zeta_{1,\ell}\nabla_0^2 \mv^k)
+ 2(\nabla_0 \mq,\zeta_1^2 \nabla_0\tta_k^\ell \nabla_0^2 \mv^k_{,\ell}) \nonumber\\
& - 2(\zeta_1 \nabla_0 \mq,\nabla_0 (\tta_k^\ell \zeta_{1,\ell} \nabla_0^2 \mv^k)
+ 2(\mq,\nabla_0(\tta_k^\ell \zeta_{1,\ell} \nabla_0\zeta_1 \nabla_0^2 \mv^k)) \label{qidentity}\\
& + (\nabla_0 \mq, \nabla_0(\zeta_1^2 \nabla_0\tta_k^\ell \nabla_0 \mv^k_{,\ell})), \nonumber
\end{align}
(\ref{aestimate1}) and (\ref{regularity1temp}) imply that
\begin{align*}
& (\mq,\tta_k^\ell (\nabla_0'^2 (\zeta_1^2 \nabla_0^2 \mv^k)_{,\ell})
\le C(M)\|\mq\|_{H^1(\Omega)} \|\mv\|_{H^3(\Omega)} \\
\le&\ C(M)C(\delta)\Big[\|\mv_t\|^2_{L^2(\Omega)} + \|\nabla
\mv\|^2_{L^2(\Omega)}
+ \|\nabla\nabla_0 \mv\|^2_{L^2(\Omega_1)} + \kappa \|\mv\|^2_{H^2(\Gamma)} \\
&\qquad\qquad\quad + \|\nabla_0^2 \mh^\eps1\|^2_{L^2(\Gamma)} +
\|F\|^2_{L^2(\Omega)} + 1\Big] + \delta \|\mv\|^2_{H^3(\Omega)}.
\end{align*}
For the viscosity term,
\begin{align*}
& \int_{\Omega} D_\tteta \mv: D_\tteta (\nabla_0^2 (\zeta_1^2 \nabla_0^2 \mv)) dx \\
=&\ \|\zeta_1 D_\tteta \nabla_0^2 \mv\|^2_{L^2(\Omega)} +
\frac{1}{2}\int_{\Omega}
\Big[\nabla_0^2(\tta_i^k\tta_i^\ell)\mv_{,\ell}^j +
\nabla_0^2(\tta_i^k\tta_j^\ell)\mv_{,\ell}^i\Big]
(\zeta_1^2 \nabla_0^2 \mv^j)_{,k} dx \\
& +
\int_{\Omega}\Big[\nabla_0(\tta_i^k\tta_i^\ell)\nabla_0\mv_{,\ell}^j
+ \nabla_0(\tta_i^k\tta_j^\ell)\nabla_0
\mv_{,\ell}^i\Big](\zeta_1^2 \nabla_0^2 \mv^j)_{,k} dx \\
& + \int_{\Omega} D_\tteta (\nabla_0^2 \mv)_i^j \tta_i^k
\zeta_1\zeta_{1,k}\nabla_0^2 \mv^j dx
\end{align*}
and hence by interpolation %(\ref{interpolation1a}) (or (\ref{interpolation1b}) if $n=2$),
\begin{align*}
&\frac{1}{2}\|\zeta_1 D_\tteta \nabla_0'^2 \mv\|^2_{L^2(\Omega)} \le
\int_{\Omega} D_\tteta \mv:
D_\tteta (\nabla_0^2 (\zeta_1^2 \nabla_0^2 \mv)) dx \\
&\qquad\qquad\quad + C(M)C(\delta)\Big[\|\nabla
\mv\|^2_{L^2(\Omega)} + \|\nabla \nabla_0
\mv\|^2_{L^2(\Omega_1')}\Big] + \delta\|\mv\|^2_{H^3(\Omega)}.
\end{align*}
Summing all the estimates, by letting $\delta_3 = \frac{\nu_2 \kappa}{2}$, we conclude that
\begin{align*}
&\ \frac{1}{2}\frac{d}{dt}\|\zeta_1 \nabla_0^2 \mv\|^2_{L^2(\Omega)}
+ \frac{\nu}{4} \|\zeta_1 D_\tteta \nabla_0^2 \mv\|^2_{L^2(\Omega)} + \frac{\nu_2\kappa}{2}\|\mv\|^2_{H^4(\Gamma)} \\
\le&\ \bar{C}\Big[\|\mv_t\|^2_{L^2(\Omega)} +
\|\mv\|^2_{H^1(\Omega)} + \|\nabla \nabla_0 \mv\|^2_{L^2(\Omega_1')}
+ \|\mv\|^2_{H^2(\Gamma)}
+ \|\nabla_0^2 \mh^\eps1\|^2_{L^2(\Gamma)}\\
&\quad  + \|F\|^2_{H^1(\Omega)} + 1\Big] + \bar{C} \int_0^t
\|\mv\|^2_{H^4(\Gamma)} ds + \delta\|\mv\|^2_{H^3(\Omega)}
\end{align*}
for some constant $\bar{C}$ depending on $M$, $\kappa$, $\epsilon$ and $\delta$.
Integrating the inequality above in time from $0$ to $t$, by (\ref{L2H1estimate}) we find that
\begin{align}
&\ \|\nabla_0^2 \mv(t)\|^2_{L^2(\Omega_1)} + \int_0^t \Big[\|\nabla
\nabla_0^2 \mv\|^2_{L^2(\Omega_1)} + \kappa\|\mv\|^2_{H^4(\Gamma)}\Big]ds  \nonumber\\
\le&\ \bar{C}N_2(u_0,F) +
\bar{C}\int_0^t\Big[\|\mv_t\|^2_{L^2(\Omega)} + \|\nabla \nabla_0
\mv\|^2_{L^2(\Omega_1')} + \|\mv\|^2_{H^2(\Gamma)}
\Big]ds \label{L2H3temp}\\
& + \bar{C} \int_0^t \int_0^s \|\mv(r)\|^2_{H^4(\Gamma)} dr +
\delta\int_0^t \|\mv\|^2_{H^3(\Omega)}ds. \nonumber
\end{align}

\noindent
By using $\nabla_0 (\zeta_1^2 \nabla_0 \mv)$ as a testing function in (\ref{eps1weakform1}), similar computations leads to
\begin{align}
& \|\nabla_0 \mv(t)\|^2_{L^2(\Omega_1)} + \int_0^t \Big[\|\nabla \nabla_0 \mv\|^2_{L^2(\Omega_1)}
+\kappa\|\mv\|^2_{H^3(\Gamma)}\Big] ds \nonumber\\
\le&\ C(M) N_2(u_0,F) + C(M,\delta)
\int_0^t\Big[\|\mv_t\|^2_{L^2(\Omega)} + \kappa
\|\mv\|^2_{H^2(\Gamma)} \Big]ds
\label{L2H2temp} \\
& + C(M) \int_0^t \int_0^s \|\mv(r)\|^2_{H^4(\Gamma)} dr ds +
\delta\int_0^t \|\mv\|^2_{H^3(\Omega)}ds. \nonumber
\end{align}
%for some constant $\bar{C}$ depending on $M$, $\epsilon$ and $\delta$.\

\subsection{Energy estimates for $v_t$ - $L^2_tH^1_x$-estimates}\label{L2H1vtsec} In this section,
we time differentiate (\ref{eps1weakform1}) and then use $\mv_t$ as a test function to obtain
\begin{align*}
& \langle \mv_{tt},\mv_t\rangle + \nu\int_{\Omega}
\Big[\tta_\ell^k(D_\tteta \mv)_{\ell,k}^i\Big]_t \mv_t^i dx +
\sigma\int_{\Gamma} \Big[\bL^\eps1_\tth(\mh^\eps1)(-\nabla_0\tth\circ\tteta^\tau,1) \Big]_t \cdot \mv_t dS \\
& +\kappa \int_{\Gamma} |\Delta_0\mv_t|^2 dS - \int_{\Omega}
(\tta_k^\ell \mq)_t \mv_{t,\ell}^k dx = \langle F_t, \mv_t\rangle -
\sigma\int_{\Gamma}
\Big[\bG^\eps1_\tth(-\nabla_0\tth\circ\tteta^\tau,1) \Big]_t \cdot
\mv_t dS.
\end{align*}
By the chain rule,
\begin{align*}
& \int_{\Gamma} \Big[(\bL_\tth^\eps1(\mh^\eps1) +
\bG_\tth^\eps1)(-\nabla_0\tth\circ\tteta^\tau,1)\Big]_t\cdot
\mv_t dS \\
=& \int_{\Gamma} \ttTheta_t
\Big[L_\tth(\mh^\eps1)\Big]^\eps1\circ\tteta^\tau
(-\nabla_0\tth\circ\tteta^\tau,1)
\cdot \mv_t dS \\
& + \int_{\Gamma} \ttTheta
\tteta^\tau_t\cdot\Big[\nabla_0[L_\tth(\mh^\eps1)]^\eps1
(-\nabla_0\tth,1) \Big]\circ\tteta^\tau
\cdot \mv_t dS \\
& + \int_{\Gamma} \ttTheta \Big[[L_\tth(\mh^\eps1)]^\eps1
(\nabla_0\tth,-1)]\Big]_t \circ\tteta^\tau \cdot \mv_t dS.
\end{align*}
By using $H^2(\Gamma)$-$H^{-2}(\Gamma)$ duality pairing with
$\epsilon$-regularization on $\ttTheta$ and $\ttv$, it follows that
\begin{align*}
& \Big|\int_{\Gamma} \Big[(\bL_\tth^\eps1(\mh^\eps1) +
\bG_\tth^\eps1)(-\nabla_0\tth\circ\tteta^\tau,1)\Big]_t\cdot
\mv_t dS\Big| \\
\le&\ C(\epsilon) \Big[\|\nabla_0^3 \mh\|_{L^2(\Gamma)} +
\|\nabla_0^2 \mh_t\|_{L^2(\Gamma)} + 1\Big]
\|\mv_t\|_{H^2(\Gamma)}\\
\le&\ C(\epsilon,\delta_3) \Big[\int_0^t \|\mv\|^2_{H^4(\Gamma)} ds
+ \|\mv\|^2_{H^2(\Gamma)} + 1\Big]
+ \delta_3 \|\mv_t\|^2_{H^2(\Gamma)} \\
\le&\ \bar{C} \Big[\int_0^t \|\mv\|^2_{H^4(\Gamma)} ds +
\|\mv\|^2_{H^1(\Omega)} + 1\Big] + \delta \|\mv\|^2_{H^3(\Omega)} +
\delta_3 \|\mv_t\|^2_{H^2(\Gamma)}
\end{align*}
for some constant $\bar{C}$ depending on $M$, $\epsilon$, $\delta$ and $\delta_3$,
where we %use the interpolation inequality (\ref{interpolation5}) to
estimate $\|\mv\|^2_{H^2(\Gamma)}$ by interpolation.

Also by interpolation %(\ref{interpolation1a}) (or (\ref{interpolation1b}) is $n=2$),
\begin{align*}
\int_{\Omega} |D_\tteta \mv_t|^2 dx =&\ 2\int_{\Omega}\Big[\tta_i^k
D_\tteta (\mv)_i^j\Big]_t \mv^j_{t,k} dx
- 2\int_{\Omega}\Big[(\tta_i^k\tta_i^\ell)_t \mv^j_{,\ell} + (\tta_i^k\tta_j^\ell)_t \mv^i_{,\ell}\Big] \mv^j_{t,k}dx \\
\le&\ 2\int_{\Omega}\Big[\tta_i^k D_\tteta (\mv)_i^j\Big]_t
\mv^j_{t,k} dx
+ C(M)C(\delta,\delta_1)\|\nabla \mv\|^2_{L^2(\Omega)} \\
& -\int_{\Omega} (\tta_{k}^\ell \mq)_t \mv_{t,\ell}^k dx + \delta
\|\mv\|^2_{H^2(\Omega)} + \delta_1\|\mv_t\|^2_{H^1(\Omega)}.
\end{align*}
Note that
$$\langle F_t, \mv_t\rangle \le C\|F_t\|_{H^1(\Omega)'}\|\mv_t\|_{H^1(\Omega)} \le C(\delta_1)\|F_t\|^2_{H^1(\Omega)'}
+ \delta_1\|\mv_t\|^2_{H^1(\Omega)}.$$
Summing all the estimates above,% and choosing $\delta_1$ and $T$ small enough,
\begin{align}
&\ \frac{1}{2}\frac{d}{dt}\|\mv_t\|^2_{L^2(\Omega)} +
\frac{\nu}{4}\|\nabla \mv_t\|^2_{L^2(\Omega)}
+ \kappa\|\Delta_0 \mv_t\|^2_{L^2(\Gamma)} \nonumber\\
\le&\ \bar{C} \Big[\int_0^t \|\mv\|^2_{H^4(\Gamma)} ds +
\|\mv\|^2_{H^1(\Omega)} + 1\Big]
+ C(\delta_1) \|F_t\|^2_{H^1(\Omega)'} \label{L2H1vttemp1}\\
& + \delta \|\mv\|^2_{H^3(\Omega)} + \delta_1
\|\mv_t\|^2_{H^1(\Omega)} + \delta_3 \|\mv_t\|^2_{H^2(\Gamma)} +
\int_{\Omega} ({\tta_k^\ell} \mq)_t \mv^k_{t,\ell} dx\nonumber
\end{align}
for some constant $\bar{C}$ depending on $M$, $\kappa$, $\delta$ and $\delta_1$. As in \cite{CoSh2005} and \cite{CoSh2006},
the integral involving the pressure $q$ has the following estimate:
%By Appendix \ref{L2H1vtapp}.2,
\begin{align*}
\int_0^t \int_{\Omega} ({\tta_k^\ell} \mq)_t \mv^k_{t,\ell} dx ds
\le&\ C(M)C(\delta,\delta_1)N_3(u_0,F) + \delta \int_0^t \|\mv\|^2_{H^3(\Omega)} ds \\
& + \delta_1\int_0^t \|\mv_t\|^2_{H^1(\Omega)} ds
\end{align*}
where
\begin{align*}
N_3(u_0,F) :=&\ \|u_0\|^2_{H^{2.5}(\Omega)} + \|u_0\|^2_{H^{4.5}(\Gamma)} + \|F\|^2_{L^2(0,T;H^1(\Omega))} \\
& + \|F_t\|^2_{L^2(0,T;H^1(\Omega)')} + \|F(0)\|^2_{H^1(\Omega)} +
1.
\end{align*}

Integrating (\ref{L2H1vttemp1}) in time from $0$ to $t$ and choosing $\delta_1, \delta_3>0$ small enough,
(\ref{L2H1estimate}) and (\ref{L2L2vtineq}) imply that, for all $t\in [0,T]$,
\begin{align}
&\ \|\mv_t(t)\|^2_{L^2(\Omega)}
+ \int_0^t \Big[\|\nabla \mv_t\|^2_{L^2(\Omega)} + \kappa\|\mv_t\|^2_{H^2(\Gamma)} \Big] ds \nonumber\\
\le&\ \bar{C} N_3(u_0,F) + \bar{C}\int_0^t \int_0^s
\|\mv(r)\|^2_{H^4(\Gamma)} dr ds + \delta \int_0^t
\|\mv\|^2_{H^3(\Omega)} ds \label{L2H1vttemp}
\end{align}
for some constant $\bar{C}$ depending on $M$, $\kappa$, $\delta$ and
$\delta_2$. In (\ref{L2H1vttemp}), (\ref{zerovalue}) is used to
bound $\|v_t(0)\|^2_{L^2(\Omega)}$.

\subsection{$\eps1$-independent estimates}\label{eps1indep}
Integrating (\ref{regularity1temp}) in time from $0$ to $t$, (\ref{L2H1estimate}), (\ref{L2L2vtineq})
and (\ref{L2H2temp}) imply that
\begin{align}
& \int_0^t \Big[\|\mv\|^2_{H^2(\Omega)} + \|\mq\|^2_{H^1(\Omega)}\Big] ds \nonumber\\
\le&\ C(M) N_1(u_0,F) + \int_0^t \Big[\|\mv_t\|^2_{L^2(\Omega)} + \|\mv\|^2_{H^2(\Gamma)} \Big]ds \nonumber\\
\le&\ \bar{C} N_3(u_0,F) + \bar{C}\int_0^t \int_0^s
\|\mv(r)\|^2_{H^4(\Gamma)} dr ds + \delta \int_0^t
\|\mv\|^2_{H^3(\Omega)} ds \label{L2H2estimate1}
\end{align}
for some constant $\bar{C}$ depending on $M$, $\kappa$ and $\delta$.
Integrating (\ref{regularity2temp}) in time from $0$ to $t$, making use of (\ref{L2H3temp}), (\ref{L2H2temp}),
(\ref{L2H1vttemp}), (\ref{L2H2estimate1}), and then choosing $\delta>0$ small enough and $T$ even smaller, we find that
\begin{align}
\int_0^t\Big[\|\mv\|^2_{H^3(\Omega)} + \|\mq\|^2_{H^2(\Omega)}\Big]
ds \le \bar{C}N_3(u_0,F) + \bar{C} \int_0^t \int_0^s
\|\mv(r)\|^2_{H^4(\Gamma)} dr ds \label{L2H3temp0}
\end{align}
for some constant $\bar{C}$ depending on $M$, $\kappa$ and $\epsilon$.

Having (\ref{L2H3temp0}), by choosing $\delta_2>0$ small enough, the estimates (\ref{L2H3temp}) can be rewritten as
\begin{align}
&\ \|\nabla_0^2 \mv(t)\|^2_{L^2(\Omega_1)}
+ \int_0^t \Big[\|\nabla \nabla_0^2 \mv\|^2_{L^2(\Omega_1)} + \kappa\|\mv\|^2_{H^4(\Gamma)} \Big]ds \nonumber\\
\le&\ \bar{C}N_3(u_0,F) + \bar{C}\int_0^t \int_0^s
\|\mv(r)\|^2_{H^4(\Gamma)} dr ds \label{L2H3temp2}
\end{align}
for some constant $\bar{C}$ depending on $M$, $\kappa$ and $\epsilon$. Therefore,
\begin{align*}
X(t) \le \bar{C} \Big[\int_0^t X(s)ds + N_3(u_0,F)\Big]
\end{align*}
where
\begin{align*}
X(t) &= \int_0^t \|\mv\|^2_{H^4(\Gamma)} ds.
\end{align*}
By the Gronwall inequality,
\begin{align}
\int_0^t \int_0^s \|\mv(r)\|^2_{H^4(\Gamma)} dr ds \le&\ \bar{C}
N_3(u_0,F) \label{hestimate}
\end{align}
for all $t\in [0,T]$ for some constant $\bar{C}$ depending on $M$, $\kappa$, and $\epsilon$. Having (\ref{hestimate}), estimates
(\ref{L2L2vtineq}), (\ref{L2H1vttemp}), (\ref{L2H3temp0}) and (\ref{L2H3temp2}) along with the standard embedding theorem lead to
\begin{align}
& \sup_{0\le t\le T} \Big[\|\mv(t)\|^2_{H^2(\Omega)} +
\|\mv_t(t)\|^2_{L^2(\Omega)} \Big]
+ \|\mv\|^2_{{\mathcal V}^3(T)} + \|\mq\|^2_{L^2(0,T;H^2(\Omega))} \nonumber\\
& + \kappa \|\mv\|^2_{L^2(0,T;H^4(\Gamma))} \le \bar{C}N_3(u_0,F)
\label{mainestimate}
\end{align}
for some constant $\bar{C}$ depending on $M$, $\kappa$ and $\epsilon$.

\subsection{Weak limits as $\eps1 \rightarrow 0$}\label{eps1indp} Since
the estimate (\ref{mainestimate}) is independent of $\eps1$,
the weak limit as $\eps1 \rightarrow 0$ of the sequence
$(\mv, \mh, \mq)$ exists.
We will denote the weak limit of $(\mv,\mh,\mq)$ by $(\nv,\nh,\nq)$.
By lower semi-continuity, (\ref{zerovalue}) and thus (\ref{mainestimate})
hold for the weak limit $(\nv,\nh,\nq)$.
Furthermore,
\begin{align}
&\ \ \langle \nv_t,\varphi\rangle + \frac{\nu}{2}\int_\Omega
D_\tteta \nv : D_\tteta\varphi dx + \sigma\int_{\Gamma}
\ttTheta \Big[[\L_\tth(\nh) (-\nabla_0\tth,1)]\circ\tteta^\tau\Big] \cdot \varphi dS \nonumber\\
&\ + \kappa \int_{\Gamma} \Delta_0 \nv\cdot \Delta_0\varphi dS
- (\nq, \tta_k^\ell\varphi_{,\ell}^k)_{L^2(\Omega)} \label{kappaweakform}\\
&= \langle F, \varphi\rangle - \sigma\int_{\Gamma} \ttTheta
\Big[[\G(\tth) (-\nabla_0\tth,1)]\circ\tteta^\tau\Big] \cdot \varphi
dS \nonumber
\end{align}
for all $\varphi\in \H1H2$ and a.a. $t\in [0,T]$.

\section{Estimates independent of $\kappa$ and $\epsilon$}\label{energyestimatesection2}
\subsection{Energy estimates which are independent of $\kappa$}
Although (\ref{mainestimate}) doesn't imply that $\nh\in
H^4(\Gamma)$, $\nh$ is indeed in $H^4(\Gamma)$ by (\ref{hevol}).
Therefore, we have that $(\nv,\nh,\nq)$ satisfies
\begin{subequations}\label{NSequation3}
\begin{alignat}{2}
\nv^i_t - \nu [\tta_\ell^k D_\tteta(\nv)_\ell^i]_{,k} &= -(\tta_i^k \nq)_{,k} + \tF^i
&& \text{in}\ (0,T)\times\Omega\,,\label{NSequation3.a}\\
\tta_i^j \nv^i_{,j} &= 0
&& \text{in}\ (0,T)\times\Omega\,,\label{NSequation3.b}\\
[\nu D_\tteta(\nv)^j_i - \nq\delta^j_i]\tta_j^\ell N_\ell &= \sigma \ttTheta [\L_\tth(\nh)(-\nabla_0\tth,1)]\circ\tteta^\tau
&& \text{on}\ (0,T)\times\Gamma\,,\label{NSequation3.c}\\
& + \sigma \ttTheta [\G_\tth(-\nabla_0\tth,1)]\circ\tteta^\tau + \kappa\Delta_0^2 \nv\nonumber\\
h_t\circ\tteta^\tau &= [(\tth_{,\alpha})\circ\tteta^\tau] v_{\alpha} - v_z
&& \text{on}\ (0,T)\times\Gamma\,,\label{NSequation3.d}\\
v &= \tu_0
&& \text{on}\ \{t=0\}\times\Omega,\label{NSequation3.e}\\
h &= 0 && \text{on}\ \{t=0\}\times\Gamma \,.\label{NSequation3.f}
\end{alignat}
\end{subequations}

Having (\ref{NSequation3.c}), (\ref{H5kappaindep}) in Appendix
\ref{bcelliptic} implies that $\nh$ is in $H^5(\Gamma)$ for a.a.
$t\in [0,T]$ with estimate
\begin{align*}
\int_0^t \|\nabla_0^2 \nh\|^2_{H^3(\Gamma)} ds \le C(\epsilon)
\int_0^t \Big[\|\nabla_0^4 \nh\|^2_{L^2(\Gamma)} +
\|\nv\|^2_{H^3(\Omega)} + \|\nq\|^2_{H^2(\Omega)} + 1 \Big]ds,
\end{align*}
where the forcing $f$ in (\ref{H5kappaindep}) is given by
\begin{align*}
[\nu D_\tteta(\nv)^j_i - \nq\delta^j_i]\tta_j^\ell N_\ell - \sigma \ttTheta [\G_\tth(-\nabla_0\tth,1)]\circ\tteta^\tau.
\end{align*}
By the same argument, (\ref{eps1lagrange}) holds with all $\theta$ replaced by $\kappa$. Therefore, by (\ref{regularity2temp})
(which follows from (\ref{eps1lagrange})),
\begin{align}
\int_0^t \|\nabla_0^2 \nh\|^2_{H^3(\Gamma)} ds \le&\
C(\epsilon)\int_0^t \Big[
\|\nv_t\|^2_{H^1(\Omega)} + \|\nabla_0^4 \nh\|^2_{L^2(\Gamma)} + \|\nabla_0^2 \nv\|^2_{H^1(\Omega_1)} \Big]ds \nonumber\\
& + C(\epsilon) N_2(u_0,F). \label{H5estimate}
\end{align}

With this extra regularity of $\nh$, the energy estimate
(\ref{mainestimate}) can be made independent of $\kappa$. In
Appendix \ref{L2H3inequality}.2, we prove that
\begin{align*}
& \frac{\nu_1}{2} \|\nabla_0^4 \nh(t)\|^2_{L^2(\Gamma)} \le \int_0^t
\int_{\Gamma} \ttTheta
\Big[[L_\tth(\nh)(-\nabla_0\tth,1)]\circ\tteta^\tau\Big]
\cdot \nabla_0^2(\zeta_1^2 \nabla_0^2 \nv) dS ds \\
&\quad + C'\int_0^t \Big[1 + \|\tv\|^2_{H^3(\Omega)} +
\|\th_t\|^2_{H^{2.5}(\Gamma)}
+ \|\th\|^2_{H^5(\Gamma)}\Big] \|\nabla_0^4 \nh\|^2_{L^2(\Gamma)} ds \\
&\qquad + C'\int_0^t \Big[\|\th\|^2_{H^5(\Gamma)} + 1\Big] ds +
\delta \int_0^t \|\nv\|^2_{H^3(\Omega)} ds + \delta_1 \int_0^t
\|\nabla_0^2 \nh\|^2_{H^3(\Gamma)} ds
\end{align*}
for some constant $C'$ depending on $M$, $\epsilon$, $\delta$ and $\delta_1$. By (\ref{H5estimate}),
\begin{align}
& \frac{\nu_1}{2} \|\nabla_0^4 \nh(t)\|^2_{L^2(\Gamma)} \le \int_0^t
\int_{\Gamma} \ttTheta
\Big[[L_\tth(\nh)(-\nabla_0\tth,1)]\circ\tteta^\tau\Big]
\cdot \nabla_0^2(\zeta_1^2 \nabla_0^2 \nv) dS ds \nonumber\\
&\qquad + C' N_2(u_0,F) + C'\int_0^t \Big[ \|\nabla_0^2 \nv\|^2_{H^1(\Omega_1)}
+ K(s) \|\nabla_0^4 \nh\|^2_{L^2(\Gamma)}\Big] ds \label{L2H3kappaindep}\\
&\qquad\qquad + \delta \int_0^t \|\nv\|^2_{H^3(\Omega)} ds +
\delta_1 \int_0^t \|\nv_t\|^2_{H^1(\Omega)} ds \nonumber
\end{align}
where
\begin{align*}
K(s) := 1 + \|\tv\|^2_{H^3(\Omega)} + \|\th_t\|^2_{H^{2.5}(\Gamma)}
+ \|\th\|^2_{H^5(\Gamma)}.
\end{align*}
With (\ref{L2H3kappaindep}), (\ref{L2H3temp}) now is replaced by
\begin{align}
&\ \Big[\|\nabla_0^2 \nv(t)\|^2_{L^2(\Omega_1)} + \|\nabla_0^4
\nh(t)\|^2_{L^2(\Gamma)}\Big] + \int_0^t \Big[\|\nabla
\nabla_0^2 \nv\|^2_{L^2(\Omega_1)} + \kappa\|\nv\|^2_{H^4(\Gamma)}\Big]ds  \nonumber\\
\le&\ C'N_2(u_0,F) + C'\int_0^t \Big[ \|\nv_t\|^2_{L^2(\Omega)} +
\|\nabla_0^2 \nv\|^2_{H^1(\Omega_1)} + K(s)\|\nabla_0^4
\nh\|^2_{L^2(\Gamma)}
\Big]ds \nonumber\\
& + \delta\int_0^t \|\nv\|^2_{H^3(\Omega)}ds + \delta_1 \int_0^t
\|\nv_t\|^2_{H^1(\Omega)} ds \label{L2H3estimate}
\end{align}
for some $C'$ depending on $M$, $\epsilon$, $\delta$ and $\delta_1$,
where (\ref{H3kappa}) is applied to bound $\kappa
\|\nv\|^2_{H^3(\Gamma)}$ (this is where $\|\nv_t\|^2_{L^2(\Omega)}$
comes from). Similar computations leads to
\begin{align}
& \Big[\|\nabla_0 \nv(t)\|^2_{L^2(\Omega_1)} + \|\nabla_0^3
\nh(t)\|^2_{L^2(\Gamma)}\Big] + \int_0^t \Big[\|\nabla \nabla_0
\nv\|^2_{L^2(\Omega_1)} + \kappa \|\nv\|^2_{H^3(\Gamma)}\Big] ds \nonumber\\
\le&\ C N_2(u_0,F) + C \int_0^t \|\nabla_0^4\nh\|^2_{L^2(\Gamma)}ds
+ \delta\int_0^t \|\nv\|^2_{H^3(\Omega)}ds \label{L2H2estimate}
\end{align}
for some constant $C$ depending on $M$ and $\delta$.

In Appendix \ref{L2H1vtapp}, we establish the following $\kappa$- and $\epsilon$-independent
inequality for the time-differentiated problem:
\begin{align*}
& \int_0^t \|\nabla_0^2 \nh_t\|^2_{L^2(\Gamma)} ds
\le \int_0^t \int_{\Gamma} \Big[[L_\tth(\nh) (\nabla_0\tth,-1)]\circ\tteta^\tau\Big]_t \cdot \nv_t dS \\
& + C N_3(u_0,F) + C \int_0^t K(s) \Big[\|\nabla_0^4
\nh\|^2_{L^2(\Gamma)}
+ \|\nabla_0^2 \nh_t\|^2_{L^2(\Gamma)} \Big]ds \\
& + (\delta + C t^{1/2})\int_0^t \|\nv\|^2_{H^3(\Omega)} ds +
(\delta_1 + C t^{1/2}) \int_0^t \|\nv_t\|^2_{H^1(\Omega)} ds +
\delta_2 \|\nabla_0^4 \nh\|^2_{L^2(\Gamma)}
\end{align*}
for some constant $C$ depending on $M$, $\delta$, $\delta_1$ and $\delta_2$. Therefore,
(\ref{L2H1vttemp}) can be replaced by
the following estimate:
\begin{align}
&\ \ \Big[\|\nv_t\|^2_{L^2(\Omega)} + \|\nabla_0^2
\nh_t\|^2_{L^2(\Gamma)}\Big]
+ \int_0^t \Big[\|\nabla \nv_t\|^2_{L^2(\Omega)} + \kappa \|\Delta_0 \nv_t\|^2_{L^2(\Gamma)}\Big] ds\nonumber\\
&\le C N_3(u_0,F) + C \int_0^t K(s) \Big[\|\nabla_0^4
\nh\|^2_{L^2(\Gamma)}
+ \|\nabla_0^2 \nh_t\|^2_{L^2(\Gamma)} \Big]ds \label{L2H1vt1}\\
&\ + (\delta + C t^{1/2})\int_0^t \|\nv\|^2_{H^3(\Omega)} ds +
(\delta_1 + C t^{1/2})\int_0^t \|\nv_t\|^2_{H^1(\Omega)} ds +
\delta_2 \|\nabla_0^4 \nh\|^2_{L^2(\Gamma)}. \nonumber
\end{align}

\subsection{$\kappa$-independent estimates}\label{epsilonindependent}
Just as in Section \ref{eps1indep}, we find that
\begin{align}
& \int_0^t \Big[\|\nv\|^2_{H^3(\Omega)} + \|\nq\|^2_{H^2(\Omega)} \Big] ds \nonumber\\
\le&\ C(M) N_2(u_0,F) + C(M) \int_0^t\Big[\|\nv_t\|^2_{H^1(\Omega)}
+ \|\nabla_0^2 \nv\|^2_{H^1(\Omega_1)}\Big] ds. \label{L2H3ineq1}
\end{align}
By choosing $\delta=\delta_1=\delta_2 = 1/8$ and $T>0$ so that $C T^{1/2} < 1/8$ in (\ref{L2H1vt1}), we find that
\begin{align}
& \int_0^t \Big[\|\nv\|^2_{H^3(\Omega)} + \|\nq\|^2_{H^2(\Omega)}
\Big] ds
\le C N_3(u_0,F) + \frac{1}{8}\|\nabla_0^4 \nh\|^2_{L^2(\Gamma)} \nonumber\\
&\quad + C(M) \int_0^t\Big[\|\nabla_0^2 \nv\|^2_{H^1(\Omega_1)} +
K(s) \Big(\|\nabla_0^4 \nh\|^2_{L^2(\Gamma)} + \|\nabla_0^2
\nh_t\|^2_{L^2(\Gamma)} \Big)\Big]ds . \label{L2H1vt2}
\end{align}
Combining the estimates (\ref{L2H1estimate}), (\ref{L2L2vtineq}), (\ref{L2H3estimate}) and
(\ref{L2H2estimate}) with (\ref{L2H1vt1}),
\begin{align*}
& \Big[\|\nv\|^2_{H^1(\Omega)} + \|\nabla_0^2
\nv\|^2_{L^2(\Omega_1)} + \|\nabla_0^2 \nh\|^2_{H^2(\Gamma)}
+ \|\nv_t\|^2_{L^2(\Omega)} + \|\nabla_0^2 \nh_t\|^2_{L^2(\Gamma)} \Big](t) \\
& + \int_0^t \Big[\|\nabla \nv\|^2_{L^2(\Omega)} + \|\nabla \nabla_0
\nv\|^2_{L^2(\Omega_1)}
+ \|\nabla \nabla_0^2 \nv\|^2_{L^2(\Omega_1)} + \|\nv_t\|^2_{H^1(\Omega)} \Big]ds \\
\le&\ C'N_3(u_0,F) + C'\int_0^t \Big[\|\nv_t\|^2_{L^2(\Omega)} +
K(s)\Big(\|\nabla_0^4 \nh\|^2_{L^2(\Gamma)} + \|\nabla_0^2
\nh_t\|^2_{L^2(\Gamma)}\Big) \Big] ds
\end{align*}
for some constant $C'$ depending on $M$ and $\epsilon$.
By the Gronwall inequality and (\ref{regularity2temp}),
\begin{align*}
& \sup_{0\le t\le T}\Big[\|\nv\|^2_{H^2(\Omega)} +
\|\nv_t\|^2_{L^2(\Omega)}
+ \|\nabla_0^2 \nh_t\|^2_{L^2(\Gamma)} + \|\nabla_0^4 \nh\|^2_{L^2(\Gamma)} \\
& \qquad\quad + \|\nq\|^2_{H^1(\Omega)} \Big](t) +
\|\nv\|^2_{{\mathcal V}^3(T)} + \|\nq\|^2_{L^2(0,T;H^2(\Omega))} \le
C(\epsilon) N_3(u_0,F).
\end{align*}

\subsection{Weak limits as $\kappa \rightarrow 0$} Just as in Section \ref{eps1indp},
the weak limit $(\ov,\oh,\oq)$ of $(\nv,\nh,\nq)$ as $\kappa\to 0$
exists in $V(T)\times L^2(0,T;H^4(\Gamma))\times
L^2(0,T;H^2(\Omega))$ with estimate
\begin{align}
& \sup_{0\le t\le T}\Big[\|\ov(t)\|^2_{H^2(\Omega)} +
\|\ov_t(t)\|^2_{L^2(\Omega)}
+ \|\nabla_0^2 \oh_t(t)\|^2_{L^2(\Gamma)} + \|\nabla_0^4 \oh(t)\|^2_{L^2(\Gamma)} \nonumber\\
& \qquad\quad + \|\oq(t)\|^2_{H^1(\Omega)} \Big] +
\|\nv\|^2_{{\mathcal V}^3(T)} + \|\oq\|^2_{L^2(0,T;H^2(\Omega))} \le
C(\epsilon) N_3(u_0,F). \label{mainestimate1}
\end{align}
(\ref{mainestimate1}) implies that for a.a. $t\in [0,T]$,
$$\|\nv(t)\|_{H^{2.5}(\Gamma)} \le \bar{C}(t)$$
for some $\bar{C}(t)$ independent of $\kappa$, and therefore for a.a. $t\in [0,T]$,
\begin{align*}
\kappa \int_{\Gamma} \Delta_0 \nv\cdot \Delta_0 \varphi dS \to 0
\end{align*}
as $\kappa\to 0$. This observation with (\ref{kappaweakform}) shows that $(\ov,\oh,\oq)$ satisfies, for a.a. $t\in [0,T]$,
\begin{align}
&\ (\nv_t,\varphi)_{L^2(\Omega)} + \frac{\nu}{2} \int_{\Omega}
D_\tteta \nv:D_\tteta(\varphi) dx + \sigma \int_{\Gamma} \ttTheta
\L_\tth(\nh) \Big[-\tth_{,\sigma}\circ\tteta^\tau
\varphi^\sigma + \varphi^z\Big] dS \nonumber\\
& - (\tta^j_i \nq,\varphi^i_{,j})_{L^2(\Omega)} = \langle
\tF,\varphi\rangle + \sigma\langle \ttTheta
\G_\tth(-\nabla_0\tth\circ\tteta^\tau,1) ,\varphi \rangle_{\Gamma}.
\label{epsilonprimeweak}
\end{align}
for all $\varphi\in\H1H2$. Since (\ref{epsilonprimeweak}) also
defines a linear functional on $H^1(\Omega)$, by the density
argument, we have that (\ref{epsilonprimeweak}) holds for all
$\varphi\in H^1(\Omega)$. As $(\ov,\oh,\oq)$ are smooth enough, we
can integrate by parts and find that $(\ov,\oh,\oq)$ satisfies
(\ref{NSequation2}) with (\ref{NSequation2.c}) replaced by
\begin{equation}\label{NSequation4}
[\nu D_\tteta(\ov)^j_i -\oq\delta^j_i]\tta_j^\ell N_\ell = \sigma
\Big[\ttTheta [(\L_\tth(\oh) + \G(\tth))
(\nabla_0\tth,-1)]\circ\tteta^\tau \Big] \ \ \text{on}\ \
(0,T)\times\Gamma\,.
\end{equation}

\subsection{$H^{5.5}$-regularity of $\nh$}\label{H5}
By (\ref{NSequation4}), we have
\begin{lemma}\label{H55regularity}
For a.a. $t\in [0,T]$, $\oh(t)\in H^{5.5}(\Gamma)$ with
\begin{align}
\|\oh\|^2_{H^{5.5}(\Gamma)} \le&\ C(M)\Big[\|\ov_t\|^2_{H^1(\Omega)}
+ \|\nabla \ov\|^2_{L^2(\Omega)} + \|\nabla_0^2
\ov\|^2_{H^1(\Omega_1)} + \|\nabla_0^4 \oh\|^2_{L^2(\Gamma)} \nonumber\\
&\qquad\quad + \|F\|^2_{H^1(\Omega)} + 1\Big] \label{H55estimate}\,,
\end{align}
and hence
\begin{align}
\|\oh\|^2_{L^2(0,T;H^{5.5}(\Gamma))} \le C(M)e^{C(M)+T}N_3(u_0,F).
\label{L2H55estimate}
\end{align}
\end{lemma}
\begin{proof}
We write the boundary condition (\ref{NSequation4}) as
\begin{align}
\L_\tth(\oh) =&\ \frac{1}{\sigma}J_\tth^{-2}(-\nabla_0\tth,1)\cdot
\Big\{\ttTheta^{-1} \Big[[\nu D_\tteta (\ov)_i^j
- \oq\delta_i^j]\tta_j^\ell N_\ell \Big]\Big\}\circ\tteta^{-\tau} - \G(\tth). \label{hequation}
\end{align}
By Corollary \ref{ellipticconstant}, $\L_\tth$ is uniformly elliptic
with the elliptic constant $\nu_1$ which is independent of $M$ which
defines our convex subset $C_T(M)$. Since $\tth\in {\mathcal H}(T)$,
$\G(\tth)\in L^2(0,T;H^{2.5}(\Gamma))\cap
L^\infty(0,T;H^1(\Gamma))$, and hence by (\ref{mainestimate}), the
right-hand side of (\ref{hequation}) is bounded in
$H^{1.5}(\Gamma)$. The important point is that these bounds are
independent of $\epsilon$. Thus, elliptic regularity of $\L_{\tth}$
proves the estimate
\begin{align*}
\|\oh\|^2_{H^{5.5}(\Gamma)} \le
C(M)\Big[\|D_\tteta(\ov)\|^2_{H^{1.5}(\Gamma)} +
\|\oq\|^2_{H^{1.5}(\Gamma)} + 1\Big]\,,
\end{align*}
so that with (\ref{regularity2temp}), (\ref{H55estimate}) is proved.
\end{proof}

\subsection{Energy estimates which are independent of $\epsilon$}
Having estimate (\ref{H55estimate}), one can follow exactly the same procedure as in Section \ref{epsilonindependent} to show
that the constant $C'$ appearing in (\ref{mainestimate1}) is independent of $\epsilon$,
provided that we have an $\epsilon$-independent version of (\ref{L2H3estimate}).
By Appendix \ref{L2H3inequality}.2, we indeed
have such an estimate:
\begin{align*}
& \frac{\nu_1}{2} \|\nabla_0^4 \oh(t)\|^2_{L^2(\Gamma)} \le \int_0^t
\int_{\Gamma} \ttTheta
\Big[[L_\tth(\oh)(-\nabla_0\tth,1)]\circ\tteta^\tau\Big]
\cdot \nabla_0^2(\zeta_1^2 \nabla_0^2 \ov) dS ds\\
&\quad + C N_2(u_0,F) + C\int_0^t K(s) \|\nabla_0^4
\oh\|^2_{L^2(\Gamma)} ds
+ (\delta + C t^{1/2}) \int_0^t \|\ov\|^2_{H^3(\Omega)} ds\\
&\qquad  + (\delta_1 + C t^{1/2}) \int_0^t \|\ov_t\|^2_{H^1(\Omega)}
ds
\end{align*}
for some constant $C$ depending on $M$, $\delta$ and $\delta_1$. Therefore, we can conclude that
\begin{align}
& \sup_{0\le t\le T}\Big[\|\ov\|^2_{H^2(\Omega)} +
\|\ov_t\|^2_{L^2(\Omega)}
+ \|\nabla_0^2 \oh_t\|^2_{L^2(\Gamma)} + \|\nabla_0^4 \oh\|^2_{L^2(\Gamma)} \label{mainestimate2}\\
& \qquad\ \ + \|\oq\|^2_{H^1(\Omega)} \Big](t) +
\|\ov\|^2_{{\mathcal V}^3(T)} + \|\oq\|^2_{L^2(0,T;H^2(\Omega))} \le
C(M)e^{C(M)+T} N_3(u_0,F). \nonumber
\end{align}

\begin{remark}
Literally speaking, we cannot use $\nabla_0^2 (\zeta_1^2 \nabla_0^2
\ov)$ as a test function in (\ref{epsilonprimeweak}) since it is not
a function in $H^1(\Omega)$. However, since $\oh\in H^{5.5}(\Gamma)$
for a.a. $t\in [0,T]$, (\ref{epsilonprimeweak}) also holds for all
$\varphi\in H^1(\Omega)'\cap H^{-1.5}(\Gamma)$ and $\nabla_0^2
(\zeta_1^2 \nabla_0^2 \ov)$ is a function of this kind.
\end{remark}

\subsection{Weak limits as $\epsilon \rightarrow 0$} The same argument leads to that
weak limits of $(\ov, \oh, \oq)$ (denoted by $(v, h, q)$) as
$\epsilon\to 0$ exists and $(v,h,q)$ satisfies (\ref{NSequation0}).

\subsection{Uniqueness}
In this section, we show that for a given $(\tv,\th)\in Y_T$, the solution to (\ref{NSequation0}) is
unique in $Y_T$.
Suppose $(v_1,h_1)$ and $(v_2,h_2)$ are two solutions (in $Y_T$) to (\ref{weakform}). Let $w=v_1-v_2$ and
$g=h_1-h_2$, then $w$ and $g$ satisfy
\begin{align}
&\ \langle w_t,\varphi\rangle + \frac{\nu}{2}\int_\Omega D_\teta w :
D_\teta\varphi dx + \sigma\int_{\Gamma}
\tTheta\Big[\tilde{L}_\th\Big(\int_0^t (\th_{,\alpha}w_{\alpha} - w_z) ds\Big)\Big]\circ\teta^\tau\times \nonumber\\
& \times(-\th_{,\alpha}\circ\teta^\tau \varphi^\alpha + \varphi^z) dS = 0 \label{uniqueness1}
\end{align}
for all $\varphi \in {\mathcal V}_v(T)$ with $w(0)=0$, where $\tilde{L}$ equals $L$ except $L_1 = L_2 =0$.
Since $w$ is in ${\mathcal V}_v(T)$, letting $w=\varphi$ in (\ref{uniqueness1}) leads to
\begin{align*}
& \Big[\|v\|^2_{H^1(\Omega)} + \|\nabla_0^2 v\|^2_{L^2(\Omega_1)} +
\|\nabla_0^4 h\|^2_{L^2(\Gamma)}
+ \|v_t\|^2_{L^2(\Omega)} + \|\nabla_0^2 h_t\|^2_{L^2(\Gamma)} \Big](t) \\
& + \int_0^t \Big[\|\nabla v\|^2_{L^2(\Omega)} + \|\nabla \nabla_0
v\|^2_{L^2(\Omega_1)}
+ \|\nabla \nabla_0^2 v\|^2_{L^2(\Omega_1)} + \|v_t\|^2_{H^1(\Omega)} \Big]ds \\
\le&\ C(M)\int_0^t K(s)\Big[\|\nabla_0^4 h\|^2_{L^2(\Gamma)} +
\|\nabla_0^2 h_t\|^2_{L^2(\Gamma)}\Big] ds.
\end{align*}

Therefore, by the Gronwall inequality and the zero initial condition ($w(0)=0$),
we have that $w$ (and hence $g$) is identical to zero.

\section{Fixed-Point argument}\label{fixedpointuniqueness}
>From previous sections, we establish a map $\Theta_T$ from $Y_T$ into $Y_T$, i.e., given $(\tv,\th) \in C_T(M)$, there
exists a unique $\Theta_T(\tv,\th) = (v,h)$ satisfying (\ref{NSequation0}). Theorem \ref{maintheorem} is then proved if
this mapping $\Theta_T$ has a fixed point. We shall make use of the Tychonoff
Fixed-Point Theorem which states as follows:
\begin{theorem}\label{Tychonoff} For a reflexive Banach space $X$, and $C\subset X$ a closed, convex, bounded subset, if
$F:C\to C$ is weakly sequentially continuous into $X$, then $F$ has at least one fixed-point.
\end{theorem}

In order to apply the Tychonoff Fixed-Point Theorem, we need to show that $\Theta(\tv,\th)\in C_T(M)$ and this is the case if $T$ is small enough. In the following
discussion, we will always assume $T$ is smaller than a fixed constant (for example, say $T\le 1$) so that the right-hand side of (\ref{mainestimate2}) can be
written as $C(M)N_3(u_0,F)$.

\begin{remark} The space $Y_T$ is not reflexive. We will treat $C_T(M)$ as a convex subset of $X_T$ and applied the Tychonoff
Fixed-Point Theorem on the space $X_T$.
\end{remark}

Before proceeding the fixed-point proof, we note that lemma \ref{a} implies that for a short time, the constant $C(M)$
in (\ref{lagrange1}) and (\ref{regularity2temp}) can be chosen to be independent of $M$. To be more precise, for almost
all $0<t\le T_1$,
\begin{align}
\|q\|^2_{L^2(\Omega)} \le C \Big[\|v_t\|^2_{L^2(\Omega)} + \|\nabla
v\|^2_{L^2(\Omega)} + \|\nabla_0^4 h\|^2_{L^2(\Gamma)} +
\|F\|^2_{L^2(\Omega)} + 1\Big], \label{lagrange}
\end{align}
\begin{align}
\|v\|^2_{H^3(\Omega)} + \|q\|^2_{H^2(\Omega)} \le&\ C
\Big[\|v_t\|^2_{H^1(\Omega)}
+ \|\nabla v\|^2_{H^1(\Omega) } + \|\nabla_0 v\|^2_{H^1(\Omega_1)} \label{fullregularity}\\
&\quad + \|\nabla_0^2 v\|^2_{H^1(\Omega)} + \|F\|^2_{H^1(\Omega)} +
1 \Big], \nonumber
\end{align}
and
\begin{align}
\|h\|^2_{H^{5.5}(\Gamma)} \le&\ C\Big[\|v_t\|^2_{H^1(\Omega)} +
\|\nabla v\|^2_{L^2(\Omega)} + \|\nabla_0^2
v\|^2_{H^1(\Omega_1)} + \|\nabla_0^4 h\|^2_{L^2(\Gamma)} \nonumber\\
&\quad + \|F\|^2_{H^1(\Omega)} + 1\Big] \label{H55ineq}
\end{align}
for some constant $C$ independent of $M$.

\subsection{Continuity in time of $h$}
By the evolution equation (\ref{NSequation0.d}) and the fact that
$v\in {\mathcal V}^3(T_1)$, $h_t\in L^2(0,T_1;H^{2.5}(\Gamma))$.
Since $h\in L^2(0,T_1;H^{5.5}(\Gamma))$, we have that $h\in
{\mathcal C}^0([0,T_1];H^4(\Gamma))$ by standard interpolation
theorem. Although there is no uniform rate that $h$ converges to
zero in $H^4(\Gamma)$, we have the following.
\begin{lemma}\label{rateofconvergence}
Let $(v,h)=\Theta_{T_1}(\tv,\th)$. Then $\|h(t)\|_{H^{2.5}(\Gamma)}$
converges to zero as $t \to 0$, uniformly for all $(\tv,\th)\in
C_{T_1}(M)$.
\end{lemma}
\begin{proof}
By the evolution equation (\ref{NSequation0.d}),
$$\|h(t)\|_{H^{2.5}(\Gamma)} \le \int_0^t \|\th_{,\alpha} v_{\alpha} - v_z\|_{H^{2.5}(\Gamma)}dS \le C(M)N_3(u_0,F)^{1/2}t^{1/2}.$$
The lemma follows directly from the inequality.
\end{proof}

By lemma \ref{rateofconvergence} and the interpolation inequality, we also have
\begin{lemma}
$\|\nabla_0^2 h(t)\|_{H^{1.5}(\Gamma)}$ converges to zero as $t \to
0$, uniformly for all $\th\in C_{T_1}(M)$ with estimate
\begin{align}
\|\nabla_0^2 h(t)\|_{H^{1.5}(\Gamma)} \le C(M)N_3(u_0,F) t^{1/4}
\label{rateofconvergence1}
\end{align}
for all $0<t\le T_1$.
\end{lemma}

\subsection{Improved energy estimates}
In order to apply the fixed-point theorem, we have to use the fact
that the forcing $F$ is in ${\mathcal V}^2(T)$. We also define a new
constant
\begin{align*}
N(u_0,F) := \|u_0\|^2_{H^{2.5}(\Omega)} + \|F\|^2_{{\mathcal
V}^2(T_1)} + \|F\|^2_{L^\infty(0,T_1;L^2(\Omega))} +
\|F(0)\|^2_{H^1(\Omega)} + 1.
\end{align*}
Note that $N_3(u_0,F) \le N(u_0,F)$.

\begin{remark}
For the linearized problem (\ref{NSequation0}), we only need $F\in
{\mathcal V}^1(T)$ to obtain a unique solution $(v,h)\in Y_T$.
\end{remark}

\subsubsection{Estimates for $\nabla_0^2 v$ near the boundary} \label{L2H3existsec}
Note that
\begin{align*}
&\ \frac{1}{2}\frac{d}{dt}\Big[\|\zeta_1\nabla_0^2
v\|^2_{L^2(\Omega)} + \sigma \int_{\Gamma} \tTheta B
\tA^{\alpha\beta\gamma\delta} \nabla_0^2 h_{,\alpha\beta} \nabla_0^2
h_{,\gamma\delta} dS \Big] + \frac{\nu}{2}\|\zeta_1 D_\teta
(\nabla_0^2
v)\|^2_{L^2(\Omega)} \\
&= \langle F, \nabla^2_0(\zeta_1^2 \nabla_0^2 v)\rangle -
\frac{\nu}{4}\int_{\Omega}
\Big[\nabla_0^2(\ta_i^k\ta_i^\ell)v_{,\ell}^j +
\nabla_0^2(\ta_i^k\ta_j^\ell)v_{,\ell}^i\Big]
(\zeta_1^2 \nabla_0^2 v^j)_{,k} dx \\
& -
\frac{\nu}{2}\int_{\Omega}\Big[\nabla_0(\ta_i^k\ta_i^\ell)\nabla_0v_{,\ell}^j
+ \nabla_0(\ta_i^k\ta_j^\ell)\nabla_0 v_{,\ell}^i\Big](\zeta_1^2 \nabla_0^2 v^j)_{,k} dx \\
& - \frac{\nu}{2}\int_{\Omega} D_\teta (\nabla_0^2 v)_i^j \ta_i^k
\zeta_1\zeta_{1,k}\nabla_0^2 v^j dx + \int_{\Omega} q \ta_k^\ell
[\nabla_0^2 (\zeta_1^2 \nabla_0^2 v^k)]_{,\ell} dx - \sigma
(\sum_{k=1}^3 I_k + \sum_{k=1}^8 J_k)
\end{align*}
where $I_k$'s and $J_k$'s are defined in Appendix \ref{L2H3inequality}.1 (with $\bar{\ }$ replaced by $\tilde{\ }$, and no
$\epsilon$ and $\eps1$).

As in \cite{CoSh2005} and \cite{CoSh2006}, we study the time
integral of the right-hand side of the identity above in order to
prove the validity of the requirement of applying Tychonoff
Fixed-Point Theorem. By interpolation and (\ref{mainestimate1}),
\begin{align*}
& \int_0^t \int_{\Omega}
\Big[\nabla_0^2(\ta_i^k\ta_i^\ell)v_{,\ell}^j
+ \nabla_0^2(\ta_i^k\ta_j^\ell)v_{,\ell}^i\Big] (\zeta_1^2 \nabla_0^2 v^j)_{,k} dx ds \\
\le&\ C \int_0^t \|\ta \ta\|_{H^2(\Omega)} \|\nabla v\|_{L^\infty(\Omega)} \|v\|_{H^3(\Omega)}ds \\
\le&\ C(M)C(\delta)\int_0^t \|v\|^{1/2}_{H^3(\Omega)}
\|v\|^{1/2}_{H^1(\Omega)}ds
+ \delta \|v\|^2_{L^2(0,T;H^3(\Omega))} \\
\le&\ C(M)C(\delta)N(u_0,F)^{1/2} \int_0^t \|v\|^{1/2}_{H^3(\Omega)}ds + \delta C(M)N(u_0,F) \\
\le&\ C(M)N(u_0,F)\Big[C(\delta) t^{3/4} + \delta\Big].
\end{align*}
Similarly,
\begin{align*}
&\ \int_0^t
\int_{\Omega}\Big[\nabla_0(\ta_i^k\ta_i^\ell)\nabla_0v_{,\ell}^j
+ \nabla_0(\ta_i^k\ta_j^\ell)\nabla_0 v_{,\ell}^i\Big](\zeta_1^2 \nabla_0^2 v^j)_{,k} dx ds \\
&+ \int_0^t \int_{\Omega} D_\teta (\nabla_0^2 v)_i^j \ta_i^k
\zeta_1\zeta_{1,k}\nabla_0^2 v^jdx ds \le C(M)N(u_0,F)\Big[t^{1/2} +
C(\delta)t + \delta\Big].
\end{align*}
For the pressure term, by interpolation and (\ref{qidentity}),
\begin{align*}
&\int_0^t \int_{\Omega} q \ta_k^\ell [\nabla_0^2 (\zeta_1^2 \nabla_0^2 v^k)]_{,\ell} dx ds \\
\le&\ C(M)\int_0^t \Big[\|q\|_{L^\infty(\Omega)} +
\|q\|_{W^{1,4}(\Omega)} + \|q\|_{H^1(\Omega)}\Big]
\|v\|_{H^3(\Omega)} ds \\
\le&\ C(M)C(\delta)\int_0^t \|q\|^2_{H^1(\Omega)}ds + \delta
\Big[\|v\|^2_{L^2(0,T;H^3(\Omega)}
+ \|q\|^2_{L^2(0,T;H^2(\Omega))} \Big] \\
\le&\ C(M)N(u_0,F)\Big[C(\delta)t^{1/2} + \delta\Big].
\end{align*}
By the estimates already established in Appendix \ref{L2H3inequality}, with the help of (\ref{hconvergence}), it is also easy
to see that
\begin{align*}
\int_0^t (\sum_{k=1}^3 I_k + \sum_{k=1}^8 J_k) ds
\le C(M)N(u_0,F) \Big[t^{1/4} + t^{1/2} + C(\delta)t^{2/3} + \delta\Big].
\end{align*}
Finally, for the forcing term, by the extra regularity we assume on $F$,
\begin{align*}
\int_0^t \langle F,\nabla_0^2 (\zeta_1^2 \nabla_0^2 v)\rangle ds
\le&\ \int_0^t \|F\|_{H^2(\Omega)}\|v\|_{H^2(\Omega)}ds
\le N(u_0,F) + \int_0^t \|v\|^2_{H^2(\Omega)} ds \\
\le&\ N(u_0,F) + C(M)N(u_0,F)t.
\end{align*}
Therefore,
\begin{align*}
&\ \Big[\|\nabla_0^2 v(t)\|^2_{L^2(\Omega_1)} + \sigma E_\th(\nabla_0^2 h)\Big]
+ \nu\int_0^t \|D_\teta (\nabla_0^2
v)\|^2_{L^2(\Omega_1)} ds \\
\le&\ \|u_0\|^2_{H^2(\Omega)} + C N(u_0,F) +
C(M)N(u_0,F)\Big[C(\delta)(t^{3/4} + t^{2/3} + t^{1/2} + t) +
\delta\Big].
\end{align*}
By Corollary \ref{ellipticconstant},
\begin{align}
&\ \Big[\|\nabla_0^2 v(t)\|^2_{L^2(\Omega_1)} + \|\nabla_0^4
h(t)\|^2_{L^2(\Gamma)} \Big]
+ \int_0^t \|\nabla_0^2 v\|^2_{H^1(\Omega_1)} ds \nonumber\\
\le&\ C N(u_0,F) + C(M)N(u_0,F)\Big[C(\delta){\mathcal O}(t) + \delta\Big] \quad\text{as}\quad t\to 0 \label{A}
\end{align}
where $C$ depends on $\nu$, $\sigma$, $\nu_1$ and the geometry of
$\Gamma$.

By similar computations, we can also conclude (the (\ref{L2H1estimate}), (\ref{L2L2vtineq}) and (\ref{L2H2estimate}) variants)
that
\begin{align}
&\ \Big[\|v(t)\|^2_{L^2(\Omega)} + \|\nabla_0^2
h(t)\|^2_{L^2(\Gamma)} \Big]
+ \int_0^t \|v\|_{H^1(\Omega)}^2 ds \nonumber\\
\le&\ C N(u_0,F) + C(M)N(u_0,F){\mathcal O}(t) \qquad\text{as}\quad t\to 0 \label{A1}\ ; \\
&\ \Big[\|\nabla_0 v(t)\|^2_{L^2(\Omega_1)} + \|\nabla_0^3
h(t)|^2_{L^2(\Gamma)} \Big]
+ \int_0^t \|\nabla_0 v\|_{H^1(\Omega_1)}^2 ds \nonumber\\
\le&\ CN(u_0,F) + C(M)N(u_0,F){\mathcal O}(t) \qquad\text{as}\quad t\to 0 \label{A2}\ ; \\
&\ \|\nabla v(t)\|^2_{L^2(\Omega)} + \int_0^t \|v_t\|^2_{L^2(\Omega)}ds \nonumber\\
\le&\ CN(u_0,F) + C(M)N(u_0,F){\mathcal O}(t)
\qquad\text{as}\quad t\to 0 \label{A3}
\end{align}
where $C$ depends on $\nu$, $\sigma$, $\nu_1$ and the geometry of
$\Gamma$.

\subsubsection{$L^2_tH^1_x$-estimate for $v_t$}
For the time-differentiated problem, we are not able to use estimates such
as those in
sections \ref{L2H1vtsec} and
\ref{L2H3existsec}, since no $\epsilon$-regularization is present;
nevertheless, we can obtain estimates at the
$\epsilon$-regularization level and then pass $\epsilon$ to the limit once the
estimate is found to be  $\epsilon$-independent. We have that
\begin{align*}
&\ \frac{1}{2}\frac{d}{dt}\|v_t\|^2_{L^2(\Omega)} +
\frac{\nu}{2}\|D_\tteta v_t\|^2_{L^2(\Omega)} +
\frac{\sigma}{2}\frac{d}{dt}\int_{\Gamma} \ttTheta \ttA^{\alpha\beta\gamma\delta} h_{t,\alpha\beta} h_{t,\gamma\delta} dS \\
=&\ \langle F_t,v_t\rangle -\nu \int_{\Omega}\Big[(\tta_i^k
\tta_j^\ell)_t v_{,\ell}^j + (\tta_i^k \tta_j^\ell)_t v_{,\ell}^i
\Big]v^j_{t,k} dx + \int_{\Omega} q_t \tta_{kt}^\ell v^k_{,\ell} dx \\
& + \frac{1}{2}\int_{\Gamma} (\ttTheta
\ttA^{\alpha\beta\gamma\delta})_t h_{t,\alpha\beta}
h_{t,\gamma\delta} dS - \int_{\Gamma}
\frac{\ttTheta}{\sqrt{\det(g_0)}}
\Big[\sqrt{\det(g_0)}(\ttA^{\alpha\beta\gamma\delta})_t
h_{,\alpha\beta}\Big]_{,\gamma\delta}
h_{tt} dS \\
& - 2\int_{\Gamma} \ttTheta_{,\gamma}\ttA^{\alpha\beta\gamma\delta}
h_{t,\alpha\beta} h_{tt,\delta} dS
- \int_{\Gamma} \ttTheta_{,\gamma\delta} \ttA^{\alpha\beta\gamma\delta} h_{t,\alpha\beta} h_{tt} dS \\
& - \int_{\Gamma} \ttTheta\Big[L_1^{\alpha\beta\gamma}
\tth_{,\alpha\beta\gamma}\Big]_t h_{tt} dS - \int_{\Gamma} \ttTheta
(L_2)_t h_{tt} dS + K_1 + K_3 + K_4 + K_5 + K_6
\end{align*}
where $K_i's$ are defined in Appendix \ref{L2H1vtapp} (without$\eps1$).

As in the previous section, the time integral of the right-hand side of the identity above is studied.
It is easy to see that
\begin{align*}
&\int_0^t \Big[ \langle F_t,v_t\rangle -\nu \Big((\tta_i^k \tta_j^\ell)_t v_{,\ell}^j
+ (\tta_i^k \tta_j^\ell)_t v_{,\ell}^i \Big)v^j_{t,k} + K_1 + K_5 + K_6\Big] ds\\
\le&\ C(M)N(u_0,F)\Big[t^{1/4} + t^{1/2} + C(\delta)(t^{1/2} + t) + \delta\Big]
\end{align*}
and by Appendix \ref{L2H1vtapp}, particularly Remark \ref{pqr},
\begin{align*}
& \int_0^t \int_{\Gamma} \Big[\frac{1}{2} (\ttTheta
\ttA^{\alpha\beta\gamma\delta})_t h_{t,\alpha\beta}
h_{t,\gamma\delta}
- \frac{\ttTheta}{\sqrt{\det(g_0)}}\Big[\sqrt{\det(g_0)}(\ttA^{\alpha\beta\gamma\delta})_t h_{,\alpha\beta}\Big]_{,\gamma\delta} h_{tt} \\
&\qquad\qquad - 2 \ttTheta_{,\gamma}\ttA^{\alpha\beta\gamma\delta} h_{t,\alpha\beta} h_{tt,\delta}
- \ttTheta_{,\gamma\delta}\ttA^{\alpha\beta\gamma\delta} h_{t,\alpha\beta} h_{tt} \Big] dS ds \\
\le&\ C(M) N(u_0,F) t^{1/2}.
\end{align*}
Special treatment needed to be done for the rest terms, and we break this procedure into several steps.

\noindent {\bf Step 1.} Let $\displaystyle{B_1 = \int_0^t
\int_{\Omega} (q \tta_k^\ell)_t v^k_{t,\ell} dx ds}$.
%As in Appendix \ref{L2H1vtapp}.2,
By the ``divergence free'' condition (\ref{NSequation2.b}),
\begin{align*}
B_1 = \int_0^t \int_{\Omega} \tta_{kt}^\ell q v^k_{t,\ell} dx ds -
\int_0^t \int_{\Omega} \tta_{kt}^\ell q_t v^k_{,\ell} dx ds.
\end{align*}
By interpolation
%(\ref{interpolation1a}) (or (\ref{interpolation1b}) if $n=2$)
and (\ref{lagrange1}),
\begin{align*}
& \Big|\int_0^t \int_{\Omega} \tta_{kt}^\ell q v^k_{t,\ell} dx ds\Big| \\
\le&\ C(M)C(\delta)\int_0^t \|q\|^2_{L^2(\Omega)} ds + \delta
\Big[\|q\|^2_{L^2(0,T;H^1(\Omega))} +
\|v_t\|^2_{L^2(0,T;H^1(\Omega))}\Big] \\
\le&\ C(M)N(u_0,F)\Big[C(\delta) t + \delta\Big].
\end{align*}
For the second integral, we have the following identity:
\begin{align*}
\int_0^t \int_{\Omega} \tta_{kt}^\ell q_t v^k_{,\ell} dx ds
= &\ \int_{\Omega}(\tta_{kt}^\ell q v^k_{,\ell})(t)dx - \int_{\Omega} \tta_{kt}^\ell(0)q(0) u_{0,\ell}^k dx \\
& -\int_0^t \int_{\Omega} (\tta_{kt}^\ell v^k_{,\ell})_t q dx ds.
\end{align*}
By the identity $\tta_{kt}^\ell = - \tta_k^i \ttv_{,i}^j \tta_j^\ell$,
\begin{align*}
\Big|\int_0^t \int_{\Omega} (\tta_{kt}^\ell v^k_{,\ell})_t q dx
ds\Big|
\le&\ \int_0^t \int_{\Omega} \Big|\Big[\tta_{ktt}^\ell v_{,\ell}^k + \tta_{kt}^\ell v_{t,\ell}^k\Big] q\Big| dx ds \\
\le&\ C(M)\int_0^t (1+\|\ttv_t\|_{H^1(\Omega)})\|\nabla
v\|_{L^4(\Omega)} \|q\|_{L^4(\Omega)} ds.
\end{align*}
Therefore,
\begin{align*}
& \Big|\int_0^t \int_{\Omega} (\tta_{kt}^\ell v^k_{,\ell})_t q dx ds\Big| \\
\le&\ C(M)C(\delta)N(u_0,F) \int_0^t \|q\|^{2\alpha}_{H^1(\Omega)}
\|q\|^{2(1-\alpha)}_{L^2(\Omega)} ds
+ \delta \int_0^t (1+\|\ttv_t\|^2_{H^1(\Omega)}) ds\\
\le&\ C(M)N(u_0,F)^2\Big[C(\delta) (t + t^{\frac{1-\alpha}{2}}) + \delta \Big]
\end{align*}
where $\displaystyle{\alpha=\frac{3}{4}}$ if $n=3$ and $\displaystyle{\alpha=\frac{1}{2}}$ if $n=2$.\\
The second integral equals $\displaystyle{\int_{\Omega} \nabla u_0 :
(\nabla u_0)^T q(0) dx}$ which is bounded by $CN(u_0,F)$. It remains
to estimate the first integral. By adding and subtracting
$\displaystyle{\int_{\Omega} \tta_{kt}^\ell(0) q v^k_{,\ell} dx}$,
we find, by $\tta_t(0)\in H^2(\Omega)$, that
\begin{align*}
\Big|\int_{\Omega}(\tta_{kt}^\ell q v^k_{,\ell})(t)dx\Big| \le&\
\int_{\Omega}\Big|(\tta_{kt}^\ell - \tta_{kt}^\ell(0))
(q v^k_{,\ell})(t)\Big| dx + \int_{\Omega} \Big|\tta_{kt}^\ell(0) q v^k_{,\ell}\Big| dx \\
\le&\ C\|\tta_t(t)-\tta_t(0)\|_{L^4(\Omega)} \|q\|_{L^2(\Omega)} \|\nabla v\|_{L^4(\Omega)} \\
& + C(\delta_1) \|\nabla v\|^2_{L^2(\Omega)} + \delta_1
\|q\|^2_{L^2(\Omega)}.
\end{align*}
Noting that
\begin{align*}
\|\nabla v\|^2_{L^2(\Omega)} =&\ \|\nabla u_0 + \int_0^t \nabla v_t
ds\|^2_{L^2(\Omega)}
\le \Big[\|\nabla u_0\|_{L^2(\Omega)} + \int_0^t \|\nabla v_t\|_{L^2(\Omega)} ds \Big]^2 \\
\le&\ 2\Big[\|u_0\|^2_{H^1(\Omega)} + C(M)N(u_0,F)t\Big],
\end{align*}
(\ref{mainestimate1}), (\ref{aestimate}c) and (\ref{lagrange}) imply
\begin{align*}
\Big|\int_{\Omega} (\tta_{kt}^\ell q v^k_{,\ell}(t) dx\Big|
\le&\ C(M)N(u_0,F) t^{1/2} + C(\delta_1)N(u_0,F) \\
& + \delta_1\Big[\|v_t\|^2_{L^2(\Omega)} + \|\nabla_0^4
h\|^2_{L^2(\Gamma)}\Big].
\end{align*}
Summing all the estimates above, we find that
\begin{align*}
|B_1| \le&\ C(\delta_1)N(u_0,F) + C(M)N(u_0,F)^2\Big[C(\delta)(t+t^\frac{1-\alpha}{2}) + \delta \Big] \\
& + \delta_1\Big[\|v_t\|^2_{L^2(\Omega)} + \|\nabla_0^4
h\|^2_{L^2(\Gamma)}\Big].
\end{align*}

\begin{remark}
It may be tempting to use an interpolation inequality to show that
$q\in {\mathcal C}([0,T];X)$ for some Banach space $X$ by analyzing
$q_t$  via Laplace's equation. The problem, however, is that the
boundary condition for $q_t$ has low regularity
$L^2(0,T;H^{-1.5}(\Gamma))$ (by the fact that $h_t\in
L^2(0,T;H^{2.5}(\Gamma))$), and thus standard elliptic estimates do
not provide the desired conclusion that $q_t\in
L^2(0,T;H^1(\Omega)')$ (and hence by interpolation, $q\in {\mathcal
C}([0,T];H^{0.5}(\Omega))$. However, suppose that $q_t\in
L^2(0,T;H^1(\Omega)')$; then we can estimate $\displaystyle{\int_0^t
\int_{\Omega} \tta_{kt}^\ell q_t v^k_{,\ell} dx ds}$ by the
following method:
\begin{align*}
\Big|\int_0^t \int_{\Omega} \tta_{kt}^\ell q_t v^k_{,\ell} dx
ds\Big|
\le& \int_0^t \|\tta_k^i \ttv_{,i}^j \tta_j^\ell v^k_{,\ell}\|_{H^1(\Omega)} \|q_t\|_{H^1(\Omega)'} ds \\
\le&\ C(M)N(u_0,F)\Big[t + t^{5/8}\Big].
\end{align*}
\end{remark}

\noindent {\bf Step 2.} Let $\displaystyle{B_2 = \int_0^t
\int_{\Gamma} \tTheta\Big[[L_1^{\alpha\beta\gamma}
\th_{,\alpha\beta\gamma}]_t h_{tt} + (L_2)_t h_{tt}\Big] dS ds}$. It
is easy to see that
\begin{align*}
\Big|\int_0^t \int_{\Gamma} \ttTheta (L_2)_t h_{tt} dS ds\Big| \le&\
C(M)\int_0^t \Big[\|v\|_{L^\infty(\Gamma)} +
\|v_t\|_{L^2(\Gamma)}\Big]ds \\
\le&\ C(M) N(u_0,F)^{1/2} (t + t^{3/4}).
\end{align*}
For parts involving $L_1$, we have
\begin{align*}
\int_0^t \int_{\Gamma} \tTheta\Big[[L_1^{\alpha\beta\gamma}
\th_{,\alpha\beta\gamma}]_t h_{tt} dS ds
=& \int_0^t \int_{\Gamma} \tTheta \Big[L_1^{\alpha\beta\gamma}\Big]_t \tth_{,\alpha\beta\gamma} h_{tt} dS ds \quad(\equiv B_2^1)\\
& + \int_0^t \int_{\Gamma} \tTheta L_1^{\alpha\beta\gamma}
\tth_{t,\alpha\beta\gamma} h_{tt} dS ds. \quad(\equiv B_2^2)
\end{align*}
By interpolation, %(\ref{interpolation1b}),
\begin{align*}
|B_2^1| \le&\ C(M)\int_0^t \|\ttTheta\|_{L^\infty(\Gamma)} \|\th\|_{W^{1,4}(\Gamma)} \|h_{tt}\|_{L^4(\Gamma)} dS ds \\
\le&\ C(M)\int_0^t \Big[\|v\|_{H^2(\Omega)} + \|v_t\|_{H^1(\Omega)} \Big]ds \\
\le&\ C(M)N(u_0,F)^{1/2} t^{1/2}
\end{align*}
while by (\ref{hconvergence}) and Corollary \ref{Lconvergence},
\begin{align*}
|B_2^2| \le&\ \int_0^t \|\ttTheta\|_{H^{1.5}(\Gamma)}
\|\th_t\|_{H^{2.5}(\Gamma)}
\|L_1^{\alpha\beta\gamma}\|_{H^{1.5}(\Gamma)} \|h_{tt}\|_{H^{0.5}(\Gamma)} ds \\
\le&\ C(M)\|L^{\alpha\beta\gamma}_1\|_{H^{1.5}(\Gamma)}\int_0^t
\|\th\|_{H^{2.5}(\Gamma)}\Big[\|v\|_{H^2(\Omega)} +
\|v_t\|_{H^1(\Omega)}\Big]ds \\
\le&\ C(M) N(u_0,F) t^{1/4}.
\end{align*}
Therefore,
\begin{align*}
|B_2|\le C(M) N(u_0,F) (t + t^{3/4} + t^{1/4}).
\end{align*}

\noindent {\bf Step 3.} Let $\displaystyle{B_3 = \int_0^t K_3 ds =
\int_0^t \int_{\Gamma} \ttTheta [L_\tth(h)]_t
[(\ttv\circ\tteta^{-\tau})\cdot (\nabla_0 h_t)] dS ds}$. The $L_1$
and $L_2$ part of $B_3$ is bounded by
\begin{align*}
C(M) \int_0^t \|\ttTheta\|_{H^{1.5}(\Gamma)}
\|\ttv\|_{H^{1.5}(\Gamma)} \|\tth\|_{H^{3.5}(\Gamma)}
\|\tth_t\|_{H^2(\Gamma)} \|h_t\|_{H^2(\Omega)} ds
\end{align*}
and hence
\begin{align*}
\Big|\int_0^t \ttTheta \Big[L_1^{\alpha\beta\gamma} \tth_{,\alpha\beta\gamma} + L_2\Big]_t  [(\ttv\circ\tteta^{-\tau})\cdot
(\nabla_0 h_t)] dS ds\Big| \le C(M) N(u_0,F) t^{1/4}.
\end{align*}
By the integration by parts formula, the highest order part of $B_3$ can be expressed as
\begin{align*}
& \int_0^t \int_{\Gamma} \frac{\ttTheta
(\ttv\circ\tteta^{-\tau})}{\sqrt{\det(g_0)}} \Big[\sqrt{\det(g_0)}
(\ttA^{\alpha\beta\gamma\delta})_t
h_{,\alpha\beta}\Big]_{,\gamma\delta} \nabla_0 h_t dS ds \quad(\equiv B_3^1)\\
+& \int_0^t \int_{\Gamma} \ttTheta (\ttv\circ\tteta^{-\tau})
\ttA^{\alpha\beta\gamma\delta} h_{t,\alpha\beta} \nabla_0
h_{t,\gamma\delta} dS ds \quad(\equiv B_3^2)\\
+&\ 2\int_0^t \int_{\Gamma} [\ttTheta
(\ttv\circ\tteta^{-\tau})]_{,\gamma}
\ttA^{\alpha\beta\gamma\delta} h_{t,\alpha\beta} \nabla_0 h_{t,\delta} dS ds \quad(\equiv B_3^3)\\
+& \int_0^t \int_{\Gamma} [\ttTheta
(\ttv\circ\tteta^{-\tau})]_{,\gamma\delta}
\ttA^{\alpha\beta\gamma\delta} h_{t,\alpha\beta} \nabla_0 h_t dS
ds.\quad(\equiv B_3^4)
\end{align*}
It is easy to see that
\begin{align*}
|B_3^1| \le&\ C(M) \int_0^t \|\ttTheta
\ttv\circ\tteta^{-\tau}\|_{H^{1.5}(\Gamma)}
\|\tth_t\|_{H^2(\Gamma)} \|h\|_{H^4(\Gamma)} \|h_t\|_{H^2(\Gamma)} dS \\
\le&\ C(M)N(u_0,F)t
\end{align*}
and
\begin{align*}
|B_3^3| \le&\ C(M) \int_0^t \|\ttTheta
\ttv\circ\tteta^{-\tau}\|_{W^{1,4}(\Gamma)}
\|\ttA\|_{L^\infty(\Gamma)} \|h_t\|_{H^2(\Gamma)} \|h_t\|_{W^{2,4}(\Gamma)} dS \\
\le&\ C(M) N(u_0,F) t^{1/2}.
\end{align*}
For $B_3^2$, by the integration by parts formula,
\begin{align*}
B_3^2 =&\ \frac{1}{2}\int_0^t \int_{\Gamma} \ttTheta
(\ttv\circ\tteta^{-\tau}) \ttA^{\alpha\beta\gamma\delta} \nabla_0
\Big[h_{t,\alpha\beta} h_{t,\gamma\delta}\Big] dS ds \\
=& -\frac{1}{2}\int_0^t \int_{\Gamma}
\frac{1}{\sqrt{\det(g_0)}}\nabla_0 \Big[\sqrt{\det(g_0)} \ttTheta
(\ttv\circ\tteta^{-\tau}) \ttA^{\alpha\beta\gamma\delta}
\Big]h_{t,\alpha\beta} h_{t,\gamma\delta} dS ds
\end{align*}
and hence
\begin{align*}
|B_3^2| \le&\ \int_0^t \Big[ \|\nabla_0 \ttTheta\|_{L^4(\Gamma)}
\|\ttv \ttA\|_{L^\infty(\Gamma)} +
\|\ttTheta\|_{L^\infty(\Gamma)} \|\ttv \ttA\|_{W^{1,4}(\Gamma)} \Big] \\
&\qquad \times \|h_t\|_{W^{2,4}(\Gamma)} \|h_t\|_{H^2(\Gamma)} ds \\
\le&\ C(M) N(u_0,F)^{1/2} \int_0^t \|v\|_{H^3(\Omega)} ds \\
\le&\ C(M) N(u_0,F) t^{1/2}.
\end{align*}
For $B_3^4$, noting that
\begin{align*}
\ttTheta_{,\gamma\delta} =&\ \det(\nabla_0\tteta^\tau)_{,\gamma\delta} \sqrt{\det(G_\tth)\circ\tteta^\tau}
+ \det(\nabla_0\tteta^\tau)_{,\gamma} \sqrt{\det(G_\tth)\circ\tteta^\tau}_{,\delta} \\
& + \det(\nabla_0\tteta^\tau)_{,\delta} \sqrt{\det(G_\tth)\circ\tteta^\tau}_{,\gamma}
+ \det(\nabla_0\tteta^\tau) \sqrt{\det(G_\tth)\circ\tteta^\tau}_{,\gamma\delta}
\end{align*}
and $\|\nabla_0 \det(\nabla_0 \tteta^\tau)\|_{H^{0.5}(\Gamma)} \le
C(M) t^{1/2}$, we find that
\begin{align*}
|B_3^4| \le&\ C(M)\int_0^t \|\nabla_0 \det(\nabla_0
\tteta^\tau)\|_{H^{0.5}(\Gamma)} \|\nabla_0^2
h_t\|_{H^{0.5}(\Gamma)} \|\nabla_0 h_t\|_{H^{1.5}(\Gamma)} ds \\
%& + \int_0^t \|\nabla_0 \ttTheta\|_{L^4(\Gamma)} \|\nabla_0 \ttv\|_{L^\infty(\Gamma)} \|\ttA\|_{L^\infty(\Gamma)} \|\nabla_0^2 h_t\|_{L^2(\Gamma)} \|\nabla_0
%h_t\|_{L^4(\Gamma)} ds \\
& + C(M)\int_0^t \|\det(\nabla_0\tteta^\tau)\|_{L^\infty(\Gamma)}
\|\nabla_0\tteta^\tau\|^2_{L^\infty(\Gamma)}
\|\nabla_0^2 h_t\|_{L^2(\Gamma)} \|\nabla_0 h_t\|_{L^2(\Gamma)} ds \\
\le&\ C(M)N(u_0,F) t^{1/2} + C(M) N(u_0,F)^{3/4} \int_0^t \|v\|_{H^3(\Omega)}^{1/2} ds \\
\le&\ C(M)N(u_0,F) (t^{1/2}+t^{3/4}).
\end{align*}
Combining all the estimates, we find that
\begin{align*}
|B_3| \le C(M)N(u_0,F) (t+t^{1/2}+ t^{3/4}).
\end{align*}

\noindent {\bf Step 4.} Let $\displaystyle{B_4 = \int_0^t K_4 ds =
\int_0^t \int_{\Gamma} \ttTheta
\Big[L_\tth(h)\Big]_t[(\nabla_0\tth,-1)_t\cdot
(v\circ\tteta^{-\tau})] dS ds}$. Integrating by parts,
\begin{align*}
B_4 =& - \int_0^t \int_{\Gamma} L_\tth(h)\Big[\ttTheta_t
(\nabla_0\tth,-1)_t \cdot (v\circ\tteta^{-\tau}) + \ttTheta
(\nabla_0 \tth,-1)_t \cdot
(v\circ\tteta^{-\tau})_t \\
& + \ttTheta (\nabla_0\tth,-1)_{tt}\cdot (v\circ\tteta^{-\tau})\Big]
dS ds + \int_{\Gamma} \ttTheta L_\th(h) [(\nabla_0\th,-1)_t\cdot
(v\circ\tteta^{-\tau})] dS.
\end{align*}
For the first integral, (\ref{htconvergence}) implies
\begin{align*}
&\ \Big|\int_{\Gamma} \ttTheta L_\th(h) [(\nabla_0\th,-1)_t\cdot (v\circ\tteta^{-\tau})] dS\Big| \\
\le&\ \|\ttTheta\|_{L^\infty(\Gamma)}
\|L_\th(h)\|_{L^2(\Gamma)} \|\nabla_0\th_t\|_{L^4(\Gamma)} \|v\circ\tteta^{-\tau}\|_{L^4(\Gamma)} \\
\le&\ C(M)N(u_0,F) \|\th_t\|_{H^{1.5}(\Gamma)} \\
\le&\ C(M)N(u_0,F) t^{1/8}.
\end{align*}
It is also easy to see that
\begin{align*}
& \Big|\int_0^t \int_{\Gamma} L_\tth(h)\Big[\ttTheta_t
(\nabla_0\tth,-1)_t \cdot (v\circ\tteta^{-\tau})
+ \ttTheta (\nabla_0\th,-1)_t\cdot (v\circ\tteta^{-\tau})_t \Big]dS ds\Big| \\
\le&\ C(M)\int_0^t \Big[\|v\|_{L^\infty(\Gamma)} +
\|v_t\|_{L^4(\Gamma)}\Big]
\|L_\th(h)\|_{L^2(\Gamma)} \|\nabla_0 \th_t\|_{L^4(\Gamma)}  ds \\
\le&\ C(M)N(u_0,F)^{1/2} \int_0^t \Big[\|v\|_{H^3(\Omega)} + \|v_t\|_{H^1(\Omega)}\Big] ds \\
\le&\ C(M)N(u_0,F) t^{1/2}.
\end{align*}
For the remaining terms, $H^{0.5}(\Gamma)$-$H^{-0.5}(\Gamma)$
duality pairing leads to
\begin{align*}
& \Big|\int_0^t \int_{\Gamma} \ttTheta L_\th(h) (\nabla_0\th,-1)_{tt}\cdot v dS ds\Big| \\
\le& \int_0^t
\|\ttTheta\|_{H^{1.5}(\Gamma)}\|L_\th(h)\|_{H^{0.5}(\Gamma)}
\|v\|_{H^{1.5}(\Gamma)} \|\th_{tt}\|_{H^{0.5}(\Gamma)} ds.
\end{align*}
By interpolation, %(\ref{interpolation4}) and (\ref{interpolation3}),
\begin{align*}
\|L_\th(h)\|_{H^{0.5}(\Gamma)} \le
C(M)\Big[\|h\|^{1/2}_{H^{5.5}(\Gamma)} \|h\|^{1/2}_{H^{3.5}(\Gamma)}
+ 1\Big]
\end{align*}
and hence
\begin{align*}
& \Big|\int_0^t\int_{\Gamma} \tTheta L_\th(h)(\nabla_0\th,-1)_{tt}\cdot (v\circ\tteta^{-\tau}) dS ds\Big| \\
\le&\ C(M)N(u_0,F)\int_0^t \|\th_{tt}\|_{H^{0.5}(\Gamma)}\Big[
\|\nabla_0^5 h\|^{1/2}_{L^2(\Gamma)} + 1\Big] ds \\
\le&\ C(M)C(\delta)N(u_0,F)\int_0^t \Big[\|\nabla_0^5 h\|_{L^2(\Gamma)}+1\Big] ds + \delta C(M)N(u_0,F) \\
\le&\ C(M)N(u_0,F)\Big[C(\delta) (t^{1/2}+t) + \delta\Big].
\end{align*}
All the inequalities above give us
\begin{align*}
|B_4| \le C(M)N(u_0,F)\Big[C(\delta)(t^{1/2}+t) + t^{1/8} + \delta\Big].
\end{align*}

\noindent
Summing all the estimates above, we find that
\begin{align*}
&\ \Big[\|v_t\|^2_{L^2(\Omega)} + \sigma \int_{\Gamma} \ttTheta
\ttA^{\alpha\beta\gamma\delta} h_{t,\alpha\beta}
h_{t,\gamma\delta}|^2 dS\Big](t)
+ \nu \int_0^t \|D_\teta v_t\|^2_{L^2(\Omega)} ds \nonumber\\
\le&\ \|v_t(0)\|^2_{L^2(\Omega)} + \sigma \int_{\Gamma}
|G_0^{\alpha\beta} h_{t,\alpha\beta}(0)|^2 dS
+ (C + C(\delta_1))N(u_0,F) \\
&  + C(M)N(u_0,F)\Big[C(\delta)(t + t^{3/4} + t^{1/2} + t^{1/4} + t^{1/8} + t^\frac{1-\alpha}{2})
+ \delta \Big] \nonumber\\
& + \delta_1\Big[\|v_t\|^2_{L^2(\Omega)} + \|\nabla_0^4
h\|^2_{L^2(\Gamma)}\Big]
\end{align*}
and by Corollary \ref{ellipticconstant},
\begin{align}
&\ \Big[\|v_t(t)\|^2_{L^2(\Omega)} + \|\nabla_0^2
h_t(t)\|^2_{L^2(\Gamma)} \Big]
+ \int_0^t \|v_t\|^2_{H^1(\Omega)} ds \nonumber\\
\le&\ (C + C(\delta_1))N(u_0,F) + C(M)N(u_0,F)\Big[C(\delta){\mathcal O}(t) + \delta \Big] \label{B}\\
& + \delta_1\Big[\|v_t\|^2_{L^2(\Omega)} + \|\nabla_0^4
h\|^2_{L^2(\Gamma)}\Big] \nonumber
\end{align}
where $C$ depends on $\nu$, $\sigma$, $\nu_1$ and the geometry of
$\Gamma$. Since this estimate is independent of $\epsilon$, we pass
$\epsilon$ to zero and conclude that the solution $(v,h)$ to
(\ref{NSequation0}) also satisfies (\ref{B}).

\subsection{Mapping from $C_T(M)$ into $C_T(M)$}
In this section, we are going to choose $M$ so that $\Theta(\tv,\th)\in C_T(M)$ if $(\tv,\th)\in C_T(M)$.

Summing (\ref{A}), (\ref{A1}), (\ref{A2}), (\ref{A3}) and (\ref{B}), by (\ref{aestimate}) we find that
\begin{align*}
&\ \ \Big[\|v(t)\|^2_{L^2(\Omega)} + \|\nabla_0
v(t)\|^2_{L^2(\Omega_1)}
+ \|\nabla_0^2 v(t)\|^2_{L^2(\Omega_1)} + \|v_t(t)\|^2_{L^2(\Omega)} \\
&\ \ + \|\nabla_0^2 h(t)\|^2_{L^2(\Gamma)} + \|\nabla_0^3
h(t)\|^2_{L^2(\Gamma)} + \|\nabla_0^4 h(t)\|^2_{L^2(\Gamma)}
+ \|\nabla_0^2 h_t(t)\|^2_{L^2(\Gamma)} \Big] \\
&\ \ + \int_0^t \Big[\|v\|^2_{H^1(\Omega)} + \|\nabla_0
v\|^2_{H^1(\Omega_1)} + \|\nabla_0^2 v\|^2_{H^1(\Omega_1)}
+ \|v_t\|^2_{H^1(\Omega)} \Big] ds \\
&\le (C+C(\delta_1)) N(u_0,F) + C(M)N(u_0,F) \Big[C(\delta){\mathcal O}(t) + \delta\Big] \\
&\ \ + \delta_1\Big[\|v_t\|^2_{L^2(\Omega)} + \|\nabla_0^4
h\|^2_{L^2(\Gamma)}\Big]
\end{align*}
where $C$ depends on $\nu$, $\sigma$, $\nu_1$ and the geometry of
$\Gamma$. Choose $\displaystyle{\delta_1 = \frac{1}{2}}$,
\begin{align*}
&\ \ \Big[\|v(t)\|^2_{L^2(\Omega)} + \|\nabla_0
v(t)\|^2_{L^2(\Omega_1)}
+ \|\nabla_0^2 v(t)\|^2_{L^2(\Omega_1)} + \|v_t(t)\|^2_{L^2(\Omega)} \\
&\ \ + \|\nabla_0^2 h(t)\|^2_{L^2(\Gamma)} + \|\nabla_0^3
h(t)\|^2_{L^2(\Gamma)} + \|\nabla_0^4 h(t)\|^2_{L^2(\Gamma)}
+ \|\nabla_0^2 h_t(t)\|^2_{L^2(\Gamma)} \Big] \\
&\ \ + \int_0^t \Big[\|v\|^2_{H^1(\Omega)} + \|\nabla_0
v\|^2_{H^1(\Omega_1)} + \|\nabla_0^2 v\|^2_{H^1(\Omega_1)}
+ \|v_t\|^2_{H^1(\Omega)} \Big] ds \\
&\le C_1 N(u_0,F) + C(M)N(u_0,F)^2 \Big[C(\delta){\mathcal O}(t) + \delta\Big]
\end{align*}
where $C_1$ depends on $\nu$, $\sigma$ $\mu$ and the geometry of
$\Gamma$. Similar to Section \ref{eps1indep}, for almost all $0<t\le
T$,
\begin{align}
&\ \ \Big[\|v(t)\|^2_{H^2(\Omega)} + \|v_t(t)\|^2_{L^2(\Omega)} +
\|\nabla_0^2 h(t)\|^2_{H^2(\Gamma)}
+ \|\nabla_0^2 h_t(t)\|^2_{L^2(\Gamma)} \Big] \nonumber\\
&\ \ + \int_0^t \Big[\|v\|^2_{H^3(\Omega)} + \|v_t\|^2_{H^1(\Omega)} + \|q\|^2_{H^2(\Omega)}\Big] ds \label{fixpoint1}\\
&\le C_2 N(u_0,F) + C(M)N(u_0,F)^2 \Big[C(\delta){\mathcal O}(t) + \delta\Big] \nonumber
\end{align}
for some constant $C_2$ depending on $C_1$.

By (\ref{hconvergence}), (\ref{htconvergence}) and (\ref{NSequation0.d}),
\begin{align}
\int_0^t \|h_t\|^2_{H^{2.5}(\Gamma)} ds &\le \int_0^t
\Big[1+\|\th\|^2_{H^{3.5}(\Gamma)}\Big]
\|v\|^2_{H^{2.5}(\Gamma)} ds \nonumber\\
&\le C(M)N(u_0,F)t^{1/4} \label{fixpoint2}
\end{align}
and
\begin{align}
\int_0^t \|h_{tt}\|^2_{H^{0.5}(\Gamma)} ds \le&\ C(M)\int_0^t
\Big[\|\th_t\|^2_{H^{1.5}(\Gamma)} \|v\|^2_{H^2(\Omega)} +
\|\th\|^2_{H^{2.5}(\Gamma)} \|v_t\|^2_{H^1(\Omega)}\Big] ds \nonumber\\
\le&\ C(M)N(u_0,F)\Big[t^{1/4} + t^{1/2}\Big]. \label{fixpoint3}
\end{align}
Also, by (\ref{H55ineq}) and (\ref{fixpoint1}),
\begin{align}
\int_0^t \|h\|^2_{H^{5.5}(\Gamma)} ds \le&\ C \int_0^t
\Big[\|v_t\|^2_{H^1(\Omega)} + \|\nabla v\|^2_{L^2(\Omega)} +
\|\nabla_0^2
v\|^2_{H^1(\Omega_1)} + \|\nabla_0^4 h\|^2_{L^2(\Gamma)} \nonumber\\
&\qquad\quad + \|F\|^2_{H^1(\Omega)} + 1\Big]ds \nonumber\\
\le&\ C_3 N(u_0,F) + C(M)N(u_0,F)^2\Big[C(\delta){\mathcal O}(t) + \delta\Big] \label{fixpoint4}
\end{align}
for some constant $C_3$ depending on $C_2$.

Combining (\ref{fixpoint1}), (\ref{fixpoint2}), (\ref{fixpoint3}) and (\ref{fixpoint4}), we have the following inequality:
\begin{align*}
&\ \ \Big[\|v(t)\|^2_{H^2(\Omega)} + \|v_t(t)\|^2_{L^2(\Omega)} +
\|h(t)\|^2_{H^4(\Gamma)} +
\|h_t(t)\|^2_{H^2(\Gamma)}\Big] \\
&\ \ + \int_0^t \Big[ \|v\|^2_{H^3(\Omega)} +
\|v_t\|^2_{H^1(\Omega)} + \|h\|^2_{H^{5.5}(\Gamma)}
+ \|h_t\|^2_{H^{2.5}(\Gamma)} + \|h_{tt}\|^2_{H^{0.5}(\Gamma)} \Big] ds \\
&\le (C_2+C_3)N(u_0,F) + C(M)N(u_0,F)^2 \Big[C(\delta){\mathcal O}(t) + \delta\Big].
\end{align*}

Let $M=2(C_2+C_3)N(u_0,F) + 1$ (and hence corresponding $T_0$ and $T$ in Lemma \ref{a} and Corollary \ref{ellipticconstant} are fixed).
Choose $\delta>0$ small enough (but fixed one such $\delta$) so that
$$C(M)N(u_0,F)^2 \delta \le \frac{1}{4}$$
and then choose $T>0$ small enough so that
$$C(M)N(u_0,F)^2 C(\delta) T \le \frac{1}{4}.$$
Then for almost all $0< t\le T$,
\begin{align*}
&\ \ \Big[\|v(t)\|^2_{H^2(\Omega)} + \|v_t(t)\|^2_{L^2(\Omega)} +
\|h(t)\|^2_{H^4(\Gamma)} +
\|h_t(t)\|^2_{H^2(\Gamma)}\Big] \\
&\ \ + \int_0^t \Big[ \|v\|^2_{H^3(\Omega)} +
\|v_t\|^2_{H^1(\Omega)} + \|h_t\|^2_{H^{2.5}(\Gamma)} +
\|h_{tt}\|^2_{H^{0.5}(\Gamma)} \Big] ds \\
&\le C_2N(u_0,F) + \frac{1}{2}
\end{align*}
and therefore
\begin{align}
&\sup_{0\le t\le T}\Big[\|v(t)\|^2_{H^2(\Omega)} +
\|v_t(t)\|^2_{L^2(\Omega)} + \|h(t)\|^2_{H^4(\Gamma)} +
\|h_t(t)\|^2_{H^2(\Gamma)}\Big] \nonumber\\
&+ \|v\|^2_{{\mathcal V}^3(T)} + \|h\|^2_{{\mathcal H}(T)} \le
2C_2N(u_0,F) + 1, \label{fixpoint}
\end{align}
or in other words,
\begin{align*}
\|(v,h)\|^2_{Y(T)} \le 2C_2N(u_0,F) + 1.
\end{align*}

\begin{remark}
(\ref{fixpoint}) implies that for $(\tv,\th)\in C_T(M)$ (with $M$
and $T$ chosen as above), the corresponding solution to the linear
problem (\ref{NSequation0}) $(v,h)=\Theta_T(\tv,\th)$ is also in
$C_T(M)$.
\end{remark}

\subsection{Weak continuity of the mapping $\Theta_T$}
\begin{lemma} The mapping $\Theta_T$ is weakly sequentially continuous from $C_T(M)$ into $C_T(M)$ (endowed with the norm of
$X_T$).
\end{lemma}
\begin{proof}
Let $(v_p,h_p)_{p\in{\mathbb N}}$ be a given sequence of elements of $C_T(M)$ weakly convergent (in $Y_T$) toward a given element
$(v,h)\in C_T(M)$ ($C_T(M)$ is sequentially weakly closed as a closed convex set) and let
$(v_{\sigma(p)},h_{\sigma(p)})_{p\in{\mathbb N}}$ be any subsequence of this sequence.

Since ${\mathcal V}^3(T)$ is compactly embedded into
$L^2(0,T;H^2(\Omega))$, we deduce the following strong convergence
results in $L^2(0,T;L^2(\Omega))$ as $p\to\infty$:\
\begin{subequations}\label{weakconvergence}
\begin{align}
(a_\ell^j)_p(a_\ell^k)_p \to a_\ell^j a_\ell^k\quad &\text{and}\quad(a_\ell^j)_p(a^\ell_k)_p \to a_\ell^j a^\ell_k, \\
[(a_\ell^j)_p(a_\ell^k)_p]_{,j} \to (a_\ell^j a_\ell^k)_{,j}\quad &\text{and}\quad [(a_\ell^j)_p(a^\ell_k)_p]_{,j}
\to (a_\ell^j a^\ell_k)_{,j}\ , \\
(a_i^k)_p &\to a_i^k.
\end{align}
\end{subequations}
Now, let $(w_p,g_p)=\Theta_T(v_p,h_p)$ and let $q_p$ be the
associated pressure, so that $(q_p)_{p\in{\mathbb N}}$ is in a
bounded set of ${\mathcal V}^2(T)$. Since $X_T$ is a reflexive
Hilbert space, let
$(w_{\sigma(p)},g_{\sigma(p)},q_{\sigma(p)})_{p\in{\mathbb N}}$ be a
subsequence weakly converging in $X_T\times {\mathcal V}^2(T)$
toward an element $(w,g,q)\in X_T\times {\mathcal V}^2(T)$. Since
$C_T(M)$ is weakly closed in $X_T$, we also have $(w,g)\in C_T(M)$.

For each $\phi\in L^2(0,T;H^1(\Omega))$, we deduce from (\ref{weakform}) (and Remark \ref{simpleht}) that
\begin{align*}
& \int_0^T \Big[(w_t, \phi)_{L^2(\Omega)}  + \frac{\mu}{2}\int_\Omega D_\eta w: D_\eta \phi dx +
\sigma\int_{\Gamma} L_h(g) (g_{,\alpha} \phi_{\alpha} - \phi_z) dS \\
&\qquad + \int_{\Omega} q a_i^j \phi_{,j}^i dx\Big] dt = \int_0^T
\langle F, \phi\rangle  dt
\end{align*}
which with the fact that, from (\ref{weakconvergence}), for all $t\in[0,T]$, $w\in {\mathcal V}_v$, provides that $(w,g)$
is a solution of (\ref{NSequation1}) in $C_T(M)$, i.e., $(w,g)=\Theta_T(v,h)$.

Therefore, we deduce that the whole sequence $(\Theta_T(v_n,h_n))_{n\in{\mathbb N}}$ weakly converges in $C_T(M)$ toward
$\Theta_T(v,h)$, which concludes the lemma.
\end{proof}

\subsection{Uniqueness}\label{uniqueness}
For the uniqueness result, we assume that $u_0$, $F$ and $\Gamma$
are smooth enough (e.g. $u_0\in H^{5.5}(\Omega)$, $F\in {\mathcal
V}^4(T)$, $\Gamma$ is a $H^{8.5}$ surface) so that $u_0$ and the
associated $u_1$, $q_0$ satisfy compatibility condition
(\ref{compatibility1}). Therefore, the solution $(v,h,q)$ are such
that $v\in {\mathcal V}^6(T)$, $q\in L^2(0,T;H^5(\Omega))$ and $h\in
L^\infty(0,T;H^7(\Gamma))\cap L^2(0,T;H^{8.5}(\Gamma))$, $h_t\in
L^\infty(0,T;H^5(\Gamma))\cap L^2(0,T;H^{5.5}(\Gamma))$, $h_{tt}\in
L^\infty(0,T;H^2(\Gamma))\cap L^2(0,T;H^{3.5}(\Gamma))$. This
implies $a\in L^\infty(0,T;H^5(\Omega))$ and hence by studying the
elliptic equation
\begin{alignat*}{2}
(a_i^\ell a_i^k q_{t,k})_{,\ell} &= \Big[\nu a_i^\ell(a_p^k a_p^j v_{,j}^i)_{,k\ell}
+ a_{it}^\ell v_{,\ell}^i + a_i^\ell F_{,\ell}\Big]_t - [(a_i^\ell a_i^k)_t q_{,k}]_{,\ell} &&\quad\text{in}\quad\Omega, \\
q_t &= J_h^{-2} \Big[\Big(\sigma L_h(h) N_i -\nu D_\eta(v)_i^\ell
a_i^j N_j\Big)_t - (a_i^j N_j)_t q \Big] a_i^\ell N_\ell
&&\quad\text{on}\quad\Gamma,
\end{alignat*}
we find that $q_t\in L^2(0,T;H^2(\Omega))$ and this implies
$v_{tt}\in L^2(0,T;H^1(\Omega))$. By the interpolation theorem, we
also conclude that $v_t\in {\mathcal C}^0([0,T];H^{2.5}(\Omega))$.

Suppose $(v,h,q)$ and $(\tv,\th,\tq)$ are two set of solutions of (\ref{NSequation}). Then
\begin{subequations}\label{differenceequation}
\begin{align}
(v-\tv)_t - \nu [a_\ell^k D_\eta (v-\tv)_\ell^i]_{,k} =&\ -a_i^k (q-\tq)_{,k} + \delta F \\
a_i^j (v - \tv)_{,j}^i =&\ \delta a \\
\Big[\nu[D_\eta (v-\tv)]_i^\ell - (q-\tq)\delta_i^\ell\Big]a_\ell^j N_j =&\ \sigma\Theta
\Big[L_h (h-\th)(-\nabla_0 h ,1)\Big]\circ\eta^\tau \\
& + \delta L_1 + \delta L_2 + \delta L_3 \nonumber\\
(h -\th)_t\circ\eta^\tau =&\ [h_{,\alpha}\circ\eta^\tau](v_\alpha-\tv_\alpha) - (v_z -\tv_z) \\
& + \delta h_1 + \delta h_2 + \delta h_3 \nonumber\\
(v-\tv)(0) =&\ 0 \\
(h-\th)(0) =&\ 0
\end{align}
\end{subequations}
where
\begin{subequations}
\begin{align}
\delta F =&\ f\circ\eta - f\circ\teta + \nu [(a_\ell^k a_\ell^j - \ta_\ell^k\ta_\ell^j) \tv_{,j}^i]_{,k}
+ \nu [(a_\ell^k a_i^j - \ta_\ell^k\ta_i^j) \tv_{,j}^\ell]_{,k} \\
& - (a_i^k - \ta_i^k)\tq_{,k} \nonumber\\
\delta a =&\ (a_i^j -\ta_i^j)\tv_{,j}^i \\
\delta L_1 = &\ \sigma \Theta \Big[L_h(\th)(\nabla_0 h - \nabla_0\th,0)\Big]\circ\eta^\tau
- \nu (a_i^k a_\ell^j - \ta_i^k \ta_\ell^j) \tv_{,k}^\ell N_j \\
& - \nu (a_\ell^k a_\ell^j - \ta_\ell^k \ta_\ell^j) \tv_{,k}^i N_j + (a_i^j - \ta_i^j)\tq N_j \nonumber\\
\delta L_2 = &\ \tTheta [L_\th(\th)\circ\eta^\tau] (\nabla_0\th\circ\eta^\tau - \nabla_0\th\circ\teta^\tau,0) \\
& + \Big[\Theta L_h(\th)\circ\eta^\tau - \tTheta L_h(\th)\circ\teta^\tau\Big] (\nabla_0\th\circ\teta^\tau,-1) \nonumber\\
\delta L_3 = &\ \Big[[L_h(\th) - L_\th(\th)](\nabla_0 \th,-1)\Big]\circ\teta^\tau \\
\delta h_1 =&\ (h_{,\alpha}\circ\eta^\tau - h_{,\alpha}\circ\teta^\tau )\tv_\alpha \\
\delta h_2 =&\ \Big[(h_{,\alpha} - \th_{,\alpha})\circ\teta^\tau\Big] \tv_\alpha \\
\delta h_3 =&\ - (\th_t\circ\eta^\tau - \th_t\circ\teta^\tau)
\end{align}
\end{subequations}
We will also use $\delta L$ and $\delta h$ to denote $\sum\limits_{k=1}^3 L_k$ and $\sum\limits_{k=1}^3 \delta h_k$ respectively.

Similar to (11.3) in \cite{CoSh2006}, we also have the following estimates.
\begin{lemma}
For $f\in H^2(\Omega)$ and $g\in H^{1.5}(\Gamma)$,
\begin{align}
\|f\circ\eta - f\circ\teta\|_{L^2(\Omega)} \le&\ C \sqrt{t}
\|f\|_{H^2(\Omega)} \Big[\int_0^t \|v -
\tv\|^2_{H^1(\Omega)} ds\Big]^{1/2}, \label{compdiff1}\\
\|g\circ\eta^\tau - g\circ\teta^\tau\|_{L^2(\Gamma)} \le&\ C
\sqrt{t} \|g\|_{H^{1.5}(\Gamma)} \Big[\int_0^t \|v
-\tv\|^2_{H^1(\Omega)} ds\Big]^{1/2}. \label{compdiff2}
\end{align}
for some constant $C$.
\end{lemma}

\begin{remark} Assuming the regularity of $h$, $h_t$ and $h_{tt}$ given in the beginning of this section, we have
\begin{align}
\|\delta L_2\|_{H^2(\Gamma)} + \|\delta h_1 + \delta
h_3\|_{H^{2.5}(\Gamma)} \le&\ C \sqrt{t}\Big[\int_0^t \|v -
\tv\|^2_{H^3(\Omega)} ds\Big]^{1/2} \label{tempineq1}
\end{align}
and
\begin{align}
&\ \|(\delta L_2)_t\|_{L^2(\Gamma)} + \|(\delta h_1 + \delta h_3)_t\|_{H^1(\Gamma)} \label{tempineq2}\\
\le&\ C \Big[\|v-\tv\|_{H^1(\Omega)} + \sqrt{t} \Big(\int_0^t \|v -
\tv\|^2_{H^2(\Omega)} ds\Big)^{1/2} \Big] \nonumber
\end{align}
and
\begin{align}
\|\nabla_0^2 (\delta h_3)_t\|_{L^2(\Gamma)} \le&\ C\Big[ \|v-\tv\|_{H^1(\Omega)} + \|v-\tv\|_{H^3(\Omega)} \nonumber\\
&\quad + \sqrt{t} \|\th_{tt}\|_{H^{3.5}(\Gamma)} \Big(\int_0^t \|v -
\tv\|^2_{H^3(\Omega)} ds\Big)^{1/2}\Big]. \label{tempineq3}
\end{align}
\end{remark}
\vspace{.1 in}

\noindent By using (\ref{compdiff1}) to estimate $\|\delta
F\|_{L^2(\Omega)}$, we find that
\begin{align}
&\ \|\nabla (v-\tv)(t)\|^2_{L^2(\Omega)} + \int_0^t \|(v-\tv)_t\|^2_{L^2(\Omega)} ds \nonumber\\
\le&\ C(\delta)\int_0^t \Big[\|v-\tv\|^2_{H^1(\Omega)}
+ \|h-\th\|^2_{H^4(\Gamma)}\Big] ds + (C(\delta)t^2 + \delta) \int_0^t \|v-\tv\|^2_{H^2(\Omega)} ds \nonumber\\
& + \delta \int_0^t \Big[\|(v-\tv)_t\|^2_{H^1(\Omega)} +
\|q-\tq\|^2_{H^1(\Omega)} \Big]ds. \label{uniqueL2L2vt}
\end{align}

\vspace{.1 in}
\noindent
For the $L^2_tH^3_x$ estimate for $v-\tv$ and the $L^2_tH^1_x$ estimate for $(v-\tv)_t$, we have
\begin{align*}
&\ \frac{1}{2}\frac{d}{dt}\Big[\|\zeta_1 \nabla_0^2
(v-\tv)\|^2_{L^2(\Omega)}
+ 2\sigma E_h(\nabla_0^2 (h-\th))\Big] + \frac{\nu}{4} \|\zeta_1 D_\tteta \nabla_0^2 (v-\tv)\|^2_{L^2(\Omega)} \\
\le&\ C\Big[\|\delta F\|^2_{H^1(\Omega)} + \|(v -
\tv)_t\|^2_{L^2(\Omega)}
+ \|\nabla (v-\tv)\|^2_{L^2(\Omega)} + \|\nabla \nabla_0 (v-\tv)\|^2_{L^2(\Omega_1')} \\
&\quad + \|\nabla_0^4 (h-\th)\|^2_{L^2(\Gamma)}\Big] +
\delta\|v-\tv\|^2_{H^3(\Omega)} + D_1 + D_2 + D_3
\end{align*}
and
\begin{align*}
&\ \frac{1}{2}\frac{d}{dt}\Big[\|(v-\tv)_t\|^2_{L^2(\Omega)} +
2\sigma E_h((h-\th)_t)\Big]
+ \frac{\nu}{4}\|\nabla (v-\tv)_t\|^2_{L^2(\Omega)} \\
\le&\ C\Big[(\|\nabla_0^4 (h -\th)\|^2_{L^2(\Gamma)} + \|\nabla_0^2
(h - \th)_t\|^2_{L^2(\Gamma)}) + \|\delta F_t\|^2_{H^1(\Omega)'}
\Big]
+ \delta\|v-\tv\|^2_{H^3(\Omega)} \\
& + E_1 + E_2 + E_3.
\end{align*}
where
\begin{alignat*}{2}
D_1 &:= \int_{\Omega} \zeta_1^2 \nabla_0^2 (q - \tq) \nabla_0^2
\delta a dx, \quad
&& D_2 := \int_{\Gamma} \Theta \Big[[L_h(h-\th)]\circ\eta^\tau\Big](\nabla_0^4 \delta h) dS, \\
D_3 &:= \int_{\Gamma} \delta L \cdot \nabla_0^4 (v - \tv) dS,
\end{alignat*}
and
\begin{alignat*}{2}
E_1 &:= \int_{\Omega} (q-\tq)_t (\delta a)_t dx,\qquad
&& E_2 := \int_{\Gamma} \Big[\Theta [L_h(h-\th)]\circ\eta^\tau\Big]_t(\delta h)_t dS, \\
E_3 &:= \int_{\Gamma} (\delta L)_t \cdot (v - \tv)_t dS.
\end{alignat*}
By using (\ref{tempineq1}) to estimate $D_i$ and (\ref{tempineq2}), (\ref{tempineq3}) to
estimate $E_i$, we obtain
\begin{align}
&\ \Big[\|\nabla_0^2 (v-\tv)(t)\|^2_{L^2(\Omega_1)} + \|\nabla_0^4
(h-\th)(t)\|^2_{L^2(\Gamma)}\Big] + \int_0^t \|\nabla
\nabla_0^2 (v-\tv)\|^2_{L^2(\Omega_1)} ds \nonumber\\
\le&\ C(\delta)\int_0^t \Big[\|(v-\tv)_t\|^2_{L^2(\Omega)} +
\|\nabla_0 (v-\tv) \|^2_{L^2(\Omega)}
+ \|\nabla_0^4 (h-\th)\|^2_{L^4(\Gamma)} \Big]ds \nonumber\\
& + (C(\delta)t^2 + \delta)\int_0^t \|v-\tv\|^2_{H^3(\Omega)} ds +
\delta \int_0^t \|q-\tq\|^2_{H^2(\Omega)} ds \label{uniqueL2H3}
\end{align}
and
\begin{align}
&\ \Big[\|(v-\tv)_t(t)\|^2_{L^2(\Omega)} + \|\nabla_0^2
(h-\th)_t\|^2_{L^2(\Gamma)}\Big]
+ \int_0^t \|\nabla (v-\tv)_t\|^2_{L^2(\Omega)} ds \nonumber\\
\le&\ C(\delta) \int_0^t \Big[\|v-\tv\|^2_{H^1(\Omega)} +
\|\nabla_0^4 (h-\th)\|^2_{L^2(\Gamma)}
+ (1+\|\th_{tt}\|^2_{H^{4.5}(\Gamma)}) \nonumber\\
&\qquad\qquad \times\|\nabla_0^2 (h-\th)_t\|^2_{L^2(\Gamma)}\Big]ds \label{uniqueL2H1vt}\\
& + (C(\delta)(t+t^2) + \delta) \int_0^t \|v-\tv\|^2_{H^3(\Omega)} ds + \delta \|q-\tq\|^2_{L^2(\Omega)} \nonumber\\
& + \delta \int_0^t \Big[\|(v-\tv)_t\|^2_{H^1(\Omega)} +
\|q-\tq\|^2_{H^2(\Omega)} \Big] ds. \nonumber
\end{align}

\noindent
Summing (\ref{uniqueL2L2vt}), (\ref{uniqueL2H3}) and (\ref{uniqueL2H1vt}), we find that
\begin{align}
Y(t) + \int_0^t Z(s) ds \le C(\delta)\int_0^t k(s) Y(s)ds + (C(\delta)(t^2+t) + \delta)\int_0^t Z(s) ds
\end{align}
where
\begin{align*}
k(t)=&\ 1+\|\th_{tt}(t)\|^2_{H^{3.5}(\Gamma)} \\
Y(t)=&\ \Big[\|v-\tv(t)\|^2_{H^1(\Omega)} +
\|\nabla_0^2(v-\tv)(t)\|^2_{L^2(\Omega_1)}
+ \|(v-\tv)_t(t)\|^2_{L^2(\Omega)} \\
& + \|h-\th\|^2_{H^4(\Gamma)} + \|(h-\th)_t\|^2_{H^2(\Gamma)}\Big], \\
Z(t)=&\ \|(v-\tv)_t(t)\|^2_{H^1(\Omega)} + \|\nabla\nabla_0^2
(v-\tv)(t)\|^2_{L^2(\Omega_1)}.
\end{align*}
By letting $\delta=1/4$ and choosing $T_u\le T$ so that $C(\delta)(T_u^2 + T_u)\le 1/4$,
\begin{align}
Y(t) + \int_0^t Z(s) ds \le C \int_0^t k(s) Y(s) ds \label{uniquegronwall}
\end{align}
for all $0<t\le T_u$. Since $Y(0)=0$, the uniqueness of the solution follows from that $Y(t)=0$ for all $0<t\le T_u$.

\appendix
\section{Elliptic regularity}\label{bcelliptic}
We establish a $\kappa$-independent elliptic estimate
for solutions of
\begin{align}
\frac{\ttTheta}{\sqrt{\det(g_0)}}\Big[\Big(\sqrt{\det(g_0)} \ttA^{\alpha\beta\gamma\delta}
\nh_{,\alpha\beta}\Big)_{,\gamma\delta}(-\nabla_0\tth,1)\Big]\circ\tteta^\tau + \kappa \Delta_0^2 \nv = f \label{bc}
\end{align}
where $\nh$ and $\nv$ satisfy (\ref{hevol}) with $\nh\in
H^4(\Gamma)$ $\nv\in H^4(\Gamma)$, and $f\in H^{1.5}(\Gamma)$.
Letting $w = \nv\circ\tteta^{-\tau}$, (\ref{bc}) is equivalent to
\begin{equation}\label{bc2}
\begin{array}{rl}
& \displaystyle{\frac{\ttTheta}{\sqrt{\det(g_0)}}\Big[\sqrt{\det(g_0)} \ttA^{\alpha\beta\gamma\delta}
\nh_{,\alpha\beta}\Big]_{,\gamma\delta} (-\nabla_0\tth,1) + \kappa \Delta_0^2 w} = f\circ\tteta^\tau
\end{array}
\end{equation}
which implies
\begin{equation}\label{bc1}
\begin{array}{rl}
& \displaystyle{\frac{\ttTheta}{\sqrt{\det(g_0)}}\Big[\sqrt{\det(g_0)} \ttA^{\alpha\beta\gamma\delta}
\nh_{,\alpha\beta}\Big]_{,\gamma\delta} + \kappa J_\tth^{-2} \Delta_0^2 w\cdot(-\nabla_0\tth,1)} \vspace{0.2cm}\\
=& J_\tth^{-2} f\circ\tteta^\tau\cdot (-\nabla_0\tth,1) \,.
\end{array}
\end{equation}
Recall that $w\cdot (-\nabla_0\tth,1) = \nh_t$.

Let $D_h$ denote the difference quotients (w.r.t. the surface
coordinate system). Taking the inner-product of (\ref{bc1}) with
$D_{-h}D_h \nabla_0^4 \nh$, by Corollary \ref{ellipticconstant} we
find that
\begin{align*}
\nu_1 \int_0^t \|D_h \nabla_0^4 \nh\|^2_{L^2(\Gamma)} ds \le
C(\epsilon)\int_0^t \Big[\|\nh\|^2_{H^2(\Gamma)} +
\|f\|^2_{H^1(\Gamma)} + \kappa \|w\|^2_{H^4(\Gamma)} \Big] ds.
\end{align*}
Since the right-hand side is independent of difference parameter
$h$, it follows that $\nh\in H^5(\Gamma)$ (as it is already a
$H^4$-function) with the estimate
\begin{align}
\int_0^t \|\nabla_0^5 \nh\|^2_{L^2(\Gamma)} ds \le
C(\epsilon)\int_0^t \Big[\|\nh\|^2_{H^2(\Gamma)} +
\|f\|^2_{H^1(\Gamma)} + \kappa \|w\|^2_{H^4(\Gamma)} \Big] ds.
\label{H5kappa}
\end{align}

Next, we obtain a $\kappa$-independent estimate of $\kappa
\|w\|^2_{H^4(\Gamma)}$. By taking the inner-product of (\ref{bc2})
with $\nabla_0^2 w$ and $\nabla_0^4 w$, we find that
\begin{align}
& \|\nabla_0^3 \nh(t)\|^2_{L^2(\Gamma)} + \kappa \int_0^t \|w\|^2_{H^3(\Gamma)} ds \nonumber\\
\le&\ C(\epsilon)\int_0^t \Big[\|\nabla_0^3 \nh\|^2_{L^2(\Gamma)} +
\|f\|^2_{L^2(\Gamma)} + \|w\|^2_{H^{2.5}(\Omega)}\Big] ds.
\label{H3kappa}
\end{align}
and
\begin{align}
& \|\nabla_0^4 \nh(t)\|^2_{L^2(\Gamma)} + \kappa \int_0^t \|w\|^2_{H^4(\Gamma)} ds \label{H4kappa}\\
\le&\ C(\epsilon,\delta_1)\int_0^t \Big[\|\nabla_0^4
\nh\|^2_{L^2(\Gamma)} + \|f\|^2_{H^{1.5}(\Gamma)} +
\|w\|^2_{H^3(\Omega)}\Big] ds + \delta_1 \int_0^t \|\nabla_0^5
\nh\|^2_{L^2(\Gamma)} dS \nonumber
\end{align}
where we use (\ref{H3kappa}) to estimate $\displaystyle{\kappa
\int_0^t \|w\|_{H^3(\Gamma)} ds}$. (\ref{H4kappa}) provides a
$\kappa$-independent estimate for $\kappa \|w\|^2_{H^4(\Gamma)}$;
hence by choosing $\delta_1>0$ small enough, (\ref{H5kappa}) implies
that for all $t\in [0,T]$,
\begin{align}
\int_0^t \|\nabla_0^2 \nh\|^2_{H^3(\Gamma)} ds \le C' \int_0^t
\Big[\|\nabla_0^4 \nh\|^2_{L^2(\Gamma)} + \|f\|^2_{H^{1.5}(\Gamma)}
+ \|w\|^2_{H^3(\Omega)}\Big] ds \label{H5kappaindep}
\end{align}
for some constant $C'$ depending on $\epsilon$.

\section{Inequalities in the estimates for $\nabla_0^2 v$ near the boundary}\label{L2H3inequality}
\subsection{$\kappa$-independent estimates}
\noindent Since $\zeta_1 \equiv 1$ on $\Gamma$ and
\begin{align*}
(-\nabla_0 \tth\circ\tteta^\tau,1)\cdot \nabla_0^4 \nv =&\ \nabla_0^4 ((-\nabla_0 \tth\circ\tteta^\tau,1)\cdot
\nv) - \nabla_0^4 (-\nabla_0 \tth\circ\tteta^\tau,1)\cdot \nv \\
& - 4\nabla_0^3 (-\nabla_0 \tth\circ\tteta^\tau,1) \cdot \nabla_0 \nv
- 6 \nabla_0^2 (-\nabla_0 \tth\circ\tteta^\tau,1)\cdot\nabla_0^2 \nv \\
& - 4\nabla_0(-\nabla_0 \tth\circ\tteta^\tau,1)\cdot\nabla_0^3 \nv\,,
\end{align*}
we find that
\begin{align*}
& \int_{\Gamma} \ttTheta \Big[ L_\tth(\nh)\circ\tteta^\tau\Big]
((-\nabla_0 \tth\circ\tteta^\tau,1)\cdot \nabla_0^2(\zeta_1^2 \nabla_0^2 \nv)) dS \nonumber\\
=& - \int_{\Gamma} \ttTheta \Big[ L_\tth(\nh)\circ\tteta^\tau\Big]
\Big[\nabla_0^4(-\nabla_0 \tth\circ\tteta^\tau,1)\cdot \nv
+ 4\nabla_0^3(-\nabla_0 \tth\circ\tteta^\tau,1)\cdot\nabla_0 \nv \\
&\qquad\qquad\qquad\qquad\qquad + 6\nabla_0^2(-\nabla_0 \tth\circ\tteta^\tau,1)\cdot \nabla_0^2 \nv \Big] dS \qquad(\equiv I_1)\\
& - 4\int_{\Gamma} \ttTheta
\Big[L_\tth(\nh)\circ\tteta^\tau\Big](\nabla_0 (-\nabla_0
\tth\circ\tteta^\tau,1)\cdot
\nabla_0^3 \nv) dS \qquad(\equiv I_2)\\
& + \int_{\Gamma} \frac{\ttTheta}{\sqrt{\det(g_0)}}
\nabla_0^2\Big[\sqrt{\det(g_0)}\Big(L_1^{\alpha\beta\gamma}\th_{,\alpha\beta\gamma}
+ L_2\Big)\circ\tteta^\tau\Big] \nabla_0^2 (\nh_t \circ\tteta^\tau) dS \qquad(\equiv I_3)\\
& + \int_{\Gamma} \frac{2\nabla_0\ttTheta}{\sqrt{\det(g_0)}}
\nabla_0\Big[\sqrt{\det(g_0)}\Big(L_1^{\alpha\beta\gamma}\th_{,\alpha\beta\gamma}
+ L_2\Big)\circ\tteta^\tau\Big] \nabla_0^2 (\nh_t \circ\tteta^\tau) dS \qquad(\equiv I_4)\\
& + \int_{\Gamma} (\nabla_0^2
\ttTheta)\Big[\Big(L_1^{\alpha\beta\gamma}\th_{,\alpha\beta\gamma} +
L_2\Big)\circ\tteta^\tau\Big]
\nabla_0^2 (\nh_t \circ\tteta^\tau) dS \qquad(\equiv I_5)\\
& + \int_{\Gamma} \frac{\ttTheta}{\sqrt{\det(g_0)}}
\Big[(\sqrt{\det(g_0)}\ttA^{\alpha\beta\gamma\delta}
\nh_{,\alpha\beta})_{,\gamma\delta}\circ\tteta^\tau\Big] \nabla_0^4
(\nh_t\circ\tteta^\tau) dS.
\end{align*}
The last term of the identity above, by a change of coordinates, can be written as
\begin{align*}
& \int_{\Gamma}
\frac{\ttTheta}{\sqrt{\det(g_0)}}\Big[(\sqrt{\det(g_0)}\ttA^{\alpha\beta\gamma\delta}
\nh_{,\alpha\beta})_{,\gamma\delta}\circ\tteta^\tau\Big] \nabla_0^4 (\nh_t\circ\tteta^\tau) dS \\
=& \int_{\Gamma} \frac{B}{\sqrt{\det(g_0)}} \nabla_0^2
(\sqrt{\det(g_0)}\ttA^{\alpha\beta\gamma\delta}
\nh_{,\alpha\beta})_{,\gamma\delta} \nabla_0^2 \nh_t dS + R_1 \\
& + 2\int_{\Gamma} \frac{\nabla_0 \ttTheta}{\sqrt{\det(g_0)}}
\nabla_0 \Big[(\sqrt{\det(g_0)}\ttA^{\alpha\beta\gamma\delta}
\nh_{,\alpha\beta})_{,\gamma\delta}\circ\tteta^\tau\Big] \nabla_0^2 (\nh_t\circ\tteta^\tau) dS \qquad (\equiv J_1) \\
& + \int_{\Gamma} \frac{\nabla_0^2 \ttTheta}{\sqrt{\det(g_0)}}
\Big[(\sqrt{\det(g_0)}\ttA^{\alpha\beta\gamma\delta}
\nh_{,\alpha\beta})_{,\gamma\delta}\circ\tteta^\tau\Big] \nabla_0^2 (\nh_t\circ\tteta^\tau) dS \qquad(\equiv J_2)\\
=&\ \frac{1}{2}\frac{d}{dt}\int_{\Gamma} B
\ttA^{\alpha\beta\gamma\delta} \nabla_0^2 \nh_{,\alpha\beta}
\nabla_0^2 \nh_{,\gamma\delta} dS + R_1'
\end{align*}
where
$B=b^t\otimes b^t\otimes b^t\otimes b^t$ with $b=\nabla_0\tteta^\tau$, and
\begin{align*}
R_1(t) =& \int_{\Gamma} b^t \otimes b^t \otimes(\nabla_0 b^t)\otimes
(\nabla_0 b^t) \nabla_0
(\sqrt{\det(g_0)}\ttA^{\alpha\beta\gamma\delta}
\nh_{,\alpha\beta})_{,\gamma\delta} \nabla_0 \nh_t dS \quad(\equiv J_3) \\
& + \int_{\Gamma} b^t\otimes b^t\otimes b^t\otimes (\nabla_0 b^t)
\nabla_0 (\sqrt{\det(g_0)}\ttA^{\alpha\beta\gamma\delta}
\nh_{,\alpha\beta})_{,\gamma\delta} \nabla_0^2 \nh_t dS \quad(\equiv J_4) \\
& + \int_{\Gamma} b^t\otimes b^t\otimes b^t\otimes (\nabla_0 b^t)
\nabla_0^2 (\sqrt{\det(g_0)}\ttA^{\alpha\beta\gamma\delta}
\nh_{,\alpha\beta})_{,\gamma\delta} \nabla_0 \nh_t dS \quad(\equiv
J_5)
\end{align*}
and
\begin{align*}
R_1'(t) =&\ R_1(t) + J_1(t) + J_2(t) - \frac{1}{2}\int_{\Gamma} (B
\ttA^{\alpha\beta\gamma\delta})_t \nabla_0^2 \nh_{,\alpha\beta}
\nabla_0^2 \nh_{,\gamma\delta} dS \quad(\equiv J_6)\\
& + 2\int_{\Gamma} \frac{B}{\sqrt{\det(g_0)}} \nabla_0
(\sqrt{\det(g_0)}\ttA^{\alpha\beta\gamma\delta}) \nabla_0
\nh_{,\alpha\beta}
\nabla_0^2 \nh_{t,\gamma\delta} dS \quad(\equiv J_7)\\
& + \int_{\Gamma} \frac{B}{\sqrt{\det(g_0)}} \nabla_0^2
(\sqrt{\det(g_0)}\ttA^{\alpha\beta\gamma\delta}) \nh_{,\alpha\beta}
\nabla_0^2 \nh_{t,\gamma\delta} dS \quad(\equiv J_8)\\
& + 2 \int_{\Gamma} \frac{B_{,\gamma}}{\sqrt{\det(g_0)}}
\nabla_0^2(\sqrt{\det(g_0)}\ttA^{\alpha\beta\gamma\delta}
\nh_{,\alpha\beta})
\nabla_0^2 \nh_{t,\delta} dS \quad(\equiv J_9)\\
& + \int_{\Gamma} \frac{B_{,\gamma\delta}}{\sqrt{\det(g_0)}}
\nabla_0^2(\sqrt{\det(g_0)}\ttA^{\alpha\beta\gamma\delta}
\nh_{,\alpha\beta}) \nabla_0^2 \nh_t dS. \quad(\equiv J_{10})
\end{align*}
It follows that
\begin{align*}
|I_1|  \le&\ C(\epsilon)(1+\|\nabla_0^4 \nh\|_{L^2(\Gamma)})\|\nabla_0^2 \nv\|_{H^1(\Omega_1')}\ , \\
|I_3| + |I_4| + |I_5| \le&\
C(M)(1+\|\th\|_{H^5(\Gamma)})\|\nabla_0^2 \nv\|_{H^1(\Omega_1)}\ ,
\end{align*}
and hence that
\begin{align*}
|I_1|+|I_3|+|I_4|+|I_5| \le C(\epsilon)\Big[\|\nabla_0^4
\nh\|^2_{L^2(\Gamma)} + \|\th\|^2_{H^5(\Gamma)} + 1\Big] +
\delta\|\nv\|^2_{H^3(\Omega)}.
\end{align*}
It follows that
\begin{align*}
|J_2| + |J_3| + |J_5| + |J_{10}| \le&\ C(\epsilon)\|\nabla_0^4 \nh\|_{L^2(\Gamma)} \|\nabla_0^2 \nh_t\|_{L^2(\Gamma)} \\
|J_6| \le&\ C(M) (\|\tv\|_{H^3(\Omega)} +
\|\th_t\|_{H^{2.5}(\Gamma)})\|\nabla_0^4 \nh\|^2_{L^2(\Gamma)}.
\end{align*}
We need only obtain $\kappa$-independent estimates for the terms
$I_2$, $J_1$, $J_4$, $J_7$, $J_8$ and $J_9$. By the
$H^{-0.5}(\Gamma)$-$H^{0.5}(\Gamma)$ duality pairing,
\begin{align*}
|I_2| \le&\ C(M)\Big[\|\nabla_0^2 \nh\|_{H^{2.5}(\Gamma)} + 1\Big]
\|\nv\|_{H^{2.5}(\Gamma)}.
\end{align*}
Therefore, by interpolation
%(\ref{interpolation4})
and Young's inequality,
\begin{align}
|I_2| \le&\ C\Big[\|\nh\|^2_{H^4(\Gamma)} + 1\Big] + \delta_1
\|\nabla_0^2 \nh\|^2_{H^3(\Gamma)} + \delta \|\nv\|^2_{H^3(\Omega)}
\label{I2estimate}
\end{align}
for some $C$ depending on $M$, $\delta$ and $\delta_1$.

For $J_1$, $J_4$ and $J_9$, we find that
\begin{align*}
&|J_1| + |J_4| + |J_9| \le C(\epsilon) \|\nh\|_{H^{4.5}(\Gamma)} \|\nv\|_{H^{2.5}(\Gamma)} \\
\le&\ C' \Big[ \|\nabla_0^2 \nh\|^2_{H^2(\Gamma)} + 1\Big] +
\delta_1 \|\nabla_0^2 \nh\|^2_{H^3(\Gamma)} + \delta
\|\nv\|^2_{H^3(\Omega)}
\end{align*}
for some constant $C'$ depending on $M$, $\epsilon$, $\delta$ and $\delta_1$.

For $J_7$ and $J_8$, by the $H^{-1.5}(\Gamma)$-$H^{1.5}(\Gamma)$
duality pairing,
\begin{align*}
|J_7| + |J_8| \le&\ C(M)
\|B\|_{H^{1.5}(\Gamma)}\|\tth\|_{H^{3.5}(\Gamma)}\|\nh\|_{H^{4.5}(\Gamma)}\|\nv\|_{H^{2.5}(\Gamma)}.
\end{align*}
Similarly to the estimate in (\ref{I2estimate}), we find that
\begin{align*}
|J_7|+|J_8| \le&\  C(M)\Big[\|\nh\|^2_{H^4(\Gamma)} + 1\Big] +
\delta_1 \|\nabla_0^2 \nh\|^2_{H^3(\Gamma)} + \delta
\|\nv\|^2_{H^3(\Omega)}.
\end{align*}
Summing all the estimates and then integrating in time from $0$ to $t$, by Corollary \ref{ellipticconstant} and the fact that
$B$ is close to 1 in the uniform norm for $T$ small,
\begin{align*}
& \frac{\nu_1}{2} \|\nabla_0^4 \nh(t)\|^2_{L^2(\Gamma)} \le \int_0^t
\int_{\Gamma} \ttTheta
\Big[[L_\tth(\nh)(-\nabla_0\tth,1)]\circ\tteta^\tau\Big]
\cdot \nabla_0^2(\zeta_1^2 \nabla_0^2 \nv) dS ds\\
&\qquad\quad + C'\int_0^t K(s) \|\nabla_0^4 \nh\|^2_{L^2(\Gamma)} ds
+ C'\int_0^t \Big[\|\th\|^2_{H^5(\Gamma)} + 1\Big] ds \\
&\qquad\qquad\qquad + \delta \int_0^t \|\nv\|^2_{H^3(\Omega)} ds +
\delta_1 \int_0^t \|\nabla_0^2 \nh\|^2_{H^3(\Gamma)} ds
\end{align*}
for some constant $C'$ depending on $M$, $\epsilon$, $\delta$ and $\delta_1$, where
\begin{align*}
K(s) := 1 + \|\tv\|^2_{H^3(\Omega)} + \|\th\|^2_{H^5(\Gamma)} +
\|\th_t\|^2_{H^{2.5}(\Gamma)}.
\end{align*}

\subsection{$\epsilon$-independent estimates}
We next obtain $\epsilon$-independent estimates for the first two terms
of $I_1$, as well as those for $I_2$, $J_1$,
$J_2$, $J_3$, $J_4$, $J_5$, $J_9$ and $J_{10}$ with $\nh$ replaced by $\oh$. Let
\begin{align*}
I_1^1 =&\ -\int_{\Gamma} \ttTheta \Big[L_\th(\oh)\circ\tteta^\tau\Big] \Big[\nabla_0^4(-\nabla_0 \tth\circ\tteta^\tau,1)\cdot \ov\Big] dS,\\
I_1^2 =&\ -4 \int_{\Gamma} \ttTheta
\Big[L_\th(\oh)\circ\tteta^\tau\Big] \Big[\nabla_0^3(-\nabla_0
\tth\circ\tteta^\tau,1)\cdot\nabla_0 \ov\Big] dS
\end{align*}
By the $H^{-1.5}(\Gamma)$-$H^{1.5}(\Gamma)$ duality pairing,
\begin{align*}
|I_1^1| + |I_1^2| \le&\ C(M)\|L_\th(\oh)\|_{H^{1.5}(\Gamma)}
\|\ov\|_{H^{2.5}(\Gamma)} \|(\nabla_0
\th)\circ\tteta^\tau\|_{H^{2.5}(\Gamma)}.
\end{align*}
Therefore, by (\ref{hconvergence}) and (\ref{H55estimate}),
\begin{align}
&\ |I_1^1| + |I_1^2| \le C(M) t^{1/4}\Big[\|\oh\|^2_{H^{5.5}(\Gamma)} + 1\Big] \|\ov\|_{H^3(\Omega)} \label{I1estimate}\\
\le&\ C t^{1/2} \Big[ \|\ov_t\|^2_{H^1(\Omega)} + \|\nabla_0^4
\oh\|^2_{L^2(\Gamma)} + \|F\|^2_{H^1(\Omega)} + 1\Big] + (\delta + C
t^{1/2}) \|\ov\|^2_{H^3(\Omega)} \nonumber
\end{align}
for some constant $C$ depending on $M$ and $\delta$.
%For $I_2$, similar to estimate (\ref{I2estimate}) in Appendix \ref{L2H3inequality}.1 except this time we put
%$H^{0.5}$-norm on $\ttTheta$ and $H^{1.5}$-norm on the others, we obtain
%\begin{align*}
%|I_2|
%\le&\ C t^{1/2} \Big[ \|\nv_t\|^2_{H^1(\Omega)} + \|\nabla_0^4
%\nh\|^2_{L^2(\Gamma)} + \|F\|^2_{H^1(\Omega)} + 1\Big] + (\delta + C t^{1/2}) \|\nv\|^2_{H^3(\Omega)}
%\end{align*}
%for some constant $C$ depending on $M$ and $\delta$.

%\begin{remark}\label{inertiaissue1}
%Taking the inertia term into account, we only have $\eta^\tau\in H^3(\Gamma)$ and $h\in H^{4.5}(\Gamma)$.
%With this regularity alone, we are not able to estimate $I_1$ and this is one of the difficulties of the general problem.
%\end{remark}
%\vspace{0.1in}

For $J_1$, we use  an $L^4$-$L^4$-$L^2$ type of H$\ddot{\rm o}$lder's
inequality and conclude that
\begin{align*}
|J_1| \le C(M)t^{1/2} \|\oh\|_{H^{5.5}(\Gamma)}
\|\ov\|_{H^{2.5}(\Gamma)}
\end{align*}
while for the other $J$ terms, we use the
$H^{0.5}(\Gamma)$-$H^{-0.5}(\Gamma)$ duality pairing to obtain
\begin{align*}
|J_2|+|J_3|+|J_4|+|J_5|+|J_9|+|J_{10}| \le C(M)t^{1/2}
\|\oh\|_{H^{5.5}(\Gamma)} \|\ov\|_{H^{2.5}(\Gamma)} \,,
\end{align*}
and hence all the $J$ terms are bounded by the same right-hand side of the inequality in (\ref{I1estimate}).
Therefore,
\begin{align*}
& \frac{\nu_1}{2} \|\nabla_0^4 \oh(t)\|^2_{L^2(\Gamma)} \le \int_0^t
\int_{\Gamma} \ttTheta
\Big[[L_\tth(\oh)(-\nabla_0\tth,1)]\circ\tteta^\tau\Big]
\cdot \nabla_0^2(\zeta_1^2 \nabla_0^2 \ov) dS ds\\
&\qquad + C N_2(u_0,F) + C\int_0^t K(s) \|\nabla_0^4
\oh\|^2_{L^2(\Gamma)} ds
+ (\delta + C t^{1/2}) \int_0^t \|\ov\|^2_{H^3(\Omega)} ds \\
&\qquad\qquad + (\delta_1 + C t^{1/2}) \int_0^t
\|\ov_t\|^2_{H^1(\Omega)} ds
\end{align*}
for some constant $C$ depending on $M$, $\delta$ and $\delta_1$.

%\begin{remark}\label{inertiaissue2}
%With inertia term, the regularity stated in Remark \ref{inertiaissue1} is not enough to estimate $J_4$, $J_5$ and $J_9$.
%However, further computations can be done and we are able to estimate $J_4$ and $J_9$ when $n=2$. Therefore, the main problem
%of estimating these $J$'s is to estimate $J_5$, which requires $h\in H^5(\Gamma)$.
%\end{remark}

\section{$L^2_tH^1_x$ estimates for $v_t$}\label{L2H1vtapp}
%\subsection{Estimates for the integral over $\Gamma$}
By the chain rule and integrating by parts,
\begin{align*}
& \int_{\Gamma}
\Big[\ttTheta[L_\tth(\nh)(-\nabla_0\tth,1)]\circ\tteta^\tau\Big]_t\cdot
\nv_t dS =
\int_{\Gamma} \ttTheta_t \Big[L_\tth(\nh)\Big]\circ\tteta^\tau (-\nabla_0\tth\circ\tteta^\tau,1)\cdot \nv_t dS \\
&\qquad\qquad + \int_{\Gamma} \ttTheta
\tteta^\tau_t\cdot\Big[\nabla_0[L_\tth(\nh)] (-\nabla_0\tth,1)
\Big]\circ\tteta^\tau
\cdot \nv_t dS \quad (\equiv K_1) \\
&\qquad\qquad + \int_{\Gamma} \ttTheta \Big[[L_\tth(\nh)]
(\nabla_0\tth,-1)]\Big]_t \circ\tteta^\tau \cdot \nv_t dS.\quad
(\equiv K_2)
\end{align*}
The first term is bounded by
\begin{align*}
C(M)\|\ttv\|_{H^3(\Omega)} \Big[\|\nabla_0^4 \nh\|_{L^2(\Gamma)} +
1\Big]\|\nv_t\|_{L^2(\Gamma)} \,.
\end{align*}
After integrating by parts, the most difficult term  to estimate in $K_1$ consists of the integral
\begin{align*}
\int_{\Gamma} \frac{\ttv}{\sqrt{\det(g_0)}}
\Big[[\sqrt{\det(g_0)}\ttA^{\alpha\beta\gamma\delta}
\nh_{,\alpha\beta}]_{,\gamma\delta}(\nabla_0
\tth,-1)\Big]\circ\tteta^\tau \nabla_0 \nv_t dS.
\end{align*}
Integrating from $0$ to $t$ and integrating by parts in time, we find that
\begin{align*}
&\ \int_0^t \int_{\Gamma} \frac{\ttv}{\sqrt{\det(g_0)}}
\Big[[\sqrt{\det(g_0)}\ttA^{\alpha\beta\gamma\delta}
\nh_{,\alpha\beta}]_{,\gamma\delta}(\nabla_0
\tth,-1)\Big]\circ\tteta^\tau \nabla_0 \nv_t dS ds \\
=%&\ \Big(\int_{\Gamma} \frac{\ttv}{\sqrt{\det(g_0)}} \Big[[(\sqrt{\det(g_0)}\ttA^{\alpha\beta\gamma\delta} \nh_{,\alpha\beta}]_{,\gamma\delta}
%(\nabla_0\tth,-1) \Big]\circ\tteta^\tau \nabla_0 \nv dS\Big)(t) \quad(\equiv K_1^1)\\
%& - \int_0^t \int_{\Gamma} \frac{\ttv_t}{\sqrt{\det(g_0)}} \Big[[\sqrt{\det(g_0)}\ttA^{\alpha\beta\gamma\delta}
%\nh_{,\alpha\beta}]_{,\gamma\delta} (\nabla_0 \tth,-1)\Big]\circ\tteta^\tau \nabla_0 \nv dS ds \quad(\equiv K_1^2)\\
%& - \int_0^t \int_{\Gamma} \frac{\ttv\otimes\ttv}{\sqrt{\det(g_0)}} \nabla_0
%\Big\{\Big[[\sqrt{\det(g_0)}\ttA^{\alpha\beta\gamma\delta} \nh_{,\alpha\beta}]_{,\gamma\delta} (\nabla_0\tth,-1)\Big]\circ\tteta^\tau\Big\}
%\nabla_0 \nv dS ds \quad(\equiv K_1^3)\\
%& - \int_0^t \int_{\Gamma} \frac{\ttv}{\sqrt{\det(g_0)}} \Big[[\sqrt{\det(g_0)}\ttA^{\alpha\beta\gamma\delta} \nh_{,\alpha\beta}]_{,\gamma\delta}
%(\nabla_0 \tth_t,0) \Big]\circ\tteta^\tau\nabla_0 \nv dS ds. \quad(\equiv K_1^4) \\
& - \int_0^t \int_{\Gamma} \frac{\ttv}{\sqrt{\det(g_0)}}
\Big[[\sqrt{\det(g_0)}\ttA^{\alpha\beta\gamma\delta}
\nh_{,\alpha\beta}]_{t,\gamma\delta}
(\nabla_0 \tth,-1) \Big]\circ\tteta^\tau\nabla_0 \nv dS ds %. \quad(\equiv K_1^5)
+ R_3
\end{align*}
where $R_3$ is bounded by
\begin{align*}
C \int_0^t \Big[1+\|\tv_t\|^2_{H^1(\Omega)}\Big] \|\nabla_0^4
\nh\|^2_{L^2(\Gamma)} ds
+ \delta_2 \|\nabla_0^4 \nh\|^2_{L^2(\Gamma)} \\
+ \delta\int_0^t \|\nv\|^2_{H^3(\Omega)} ds + (\delta + C
t^{1/2})\int_0^t \|\nv_t\|^2_{H^1(\Omega)} ds
\end{align*}
for some constant $C$ depending on $M$, $\delta$ and $\delta_2$.
%With Young's inequality,
%\begin{align*}
%|K_1^1|
%\le&\ C(M)C(\delta,\delta_2) t \int_0^t \|\nv_t\|^2_{H^1(\Omega)} ds + \delta \int_0^t \Big[\|\nv_t\|^2_{H^1(\Omega)} +
%\|\nv\|^2_{H^3(\Omega)}\Big]ds \\
%& + \delta_2 \|\nabla_0^4 \nh\|^2_{L^2(\Gamma)}.
%\end{align*}
%and
%\begin{align*}
%|K_1^2|
%\le&\ C(M)C(\delta) \int_0^t \|\tv_t\|^2_{H^1(\Omega)} \|\nabla_0^4 \nh\|^2_{L^2(\Gamma)} ds + \delta \int_0^t
%\|\nv\|^2_{H^3(\Omega)} ds
%\end{align*}
%and
%\begin{align*}
%|K_1^4|
%\le&\ C(M)C(\delta) \int_0^t \|\nabla_0^4 \nh\|^2_{L^2(\Gamma)} ds + \delta \int_0^t \|\nv\|^2_{H^3(\Omega)} ds.
%\end{align*}
%Also,
%\begin{align*}
%|K_1^3| \le&\ C(M) \int_0^t\|(\tv\otimes \tv)\nabla_0 \nv\|_{H^1(\Gamma)} \|\nabla_0^4 \nh\|_{L^2(\Gamma)} ds \\
%\le&\ C(M)C(\delta) \int_0^t \|\nabla_0^4 \nh\|^2_{L^2(\Gamma)} ds + \delta \int_0^t \|\nv\|^2_{H^3(\Omega)} ds.
%\end{align*}
Next, using that
\begin{align*}
[(-\nabla_0 \tth,1)\circ\tteta^\tau] \cdot \nabla_0 \nv
=&\ b^t (\nabla_0 \nh_t)\circ\tteta^\tau + b^t (\nabla_0^2 \tth \circ\tteta^\tau,0)\cdot \nv,
\end{align*}
and integrating by parts, we find that the integral on the right-hand side is identical to
\begin{align*}
&\ \frac{1}{2} \int_0^t \int_{\Gamma} \frac{1}{\sqrt{\det(g_0)}}
\nabla_0 \Big[\sqrt{\det(g_0)}\ttTheta \ttv b^t
\ttA^{\alpha\beta\gamma\delta}\Big] \nh_{t,\alpha\beta}
\nh_{t,\gamma\delta} dS ds + R_4
\end{align*}
where
\begin{align*}
|R_4| \le&\ C(M)C(\delta) \int_0^t \|\nabla_0^4
\nh\|^2_{L^2(\Gamma)} ds + \delta \int_0^t \|\nv\|^2_{H^3(\Omega)}
ds.
\end{align*}
By interpolation, the integral part is bounded by
\begin{align*}
C\Big[N(u_0,F) + \int_0^t \|\nabla_0^4 \nh\|^2_{L^2(\Gamma)} ds\Big]
+ \delta \int_0^t \|\nv\|^2_{H^3(\Omega)} ds + C t \int_0^t
\|\nv_t\|^2_{H^1(\Omega)} ds
\end{align*}
for some constant $C$ depending on $M$ and $\delta$.
Therefore, $K_1$ satisfies
\begin{align}
& \Big|\int_0^t K_1 ds\Big| \le C \int_0^t
\Big[K(s)\Big(\|\nabla_0^4 \nh\|^2_{L^2(\Gamma)} + \|\nabla_0^2
\nh_t\|^2_{L^2(\Gamma)}\Big)+1\Big] ds
+ \delta_2 \|\nabla_0^4 \nh\|^2_{L^2(\Gamma)} \nonumber\\
&\qquad\quad + (\delta + C t^{1/2}) \int_0^t \|\nv\|^2_{H^3(\Omega)}
ds + (\delta + C t^{1/2})\int_0^t \|\nv_t\|^2_{H^1(\Omega)} ds
\label{K1estimate}
\end{align}
for some constant $C$ depending on $M$, $\delta$ and $\delta_2$.

For $K_2$, by time differentiating the evolution equation, we find that
\begin{align*}
(-\nabla_0 \tth\circ\tteta^\tau, 1) \nv_t =&\ \nh_{tt}\circ\tteta^\tau
+ \ttv^\tau\cdot(\nabla_0 \nh_t)\circ\tteta^\tau -
\ttv^\tau\cdot(\nabla_0^2\tth\circ\tteta^\tau,0) \cdot \nv \\
& - (\nabla_0 \tth_t\circ\tteta^\tau,0)\cdot \nv
\end{align*}
and hence (after a change of coordinates)
\begin{align*}
K_2 =&\
\int_{\Gamma} [L_\tth(\nh)]_t \nh_{tt} dS + \int_{\Gamma} [L_\tth(\nh)]_t [(\ttv^\tau\circ\tteta^{-\tau})\cdot(\nabla_0 \nh_t)] dS \quad(\equiv K_3) \\
& - \int_{\Gamma} [L_\tth(\nh)]_t [(\nabla_0\tth_t,0)\cdot (\nv\circ\tteta^{-\tau})] dS \quad(\equiv K_4)\\
& -\int_{\Gamma} [L_\tth(\nh)]_t
[(\ttv^\tau\circ\tteta^{-\tau})\cdot (\nabla_0^2 \tth,0)
(\nv\circ\tteta^{-\tau})] dS. \quad(\equiv K_5) \\
& +\int_{\Gamma} [L_\tth(\nh)] [(\nabla_0\tth_t,0) \cdot
(\nv_t\circ\tteta^{-\tau})] dS \quad(\equiv K_6).
\end{align*}
For the first term, we have
\begin{align}
&\ \int_{\Gamma} [L_\tth(\nh)]_t \nh_{tt} dS =
\frac{1}{2}\frac{d}{dt} \int_{\Gamma} \ttA^{\alpha\beta\gamma\delta}
\nh_{t,\alpha\beta}
\nh_{t,\gamma\delta} dS \nonumber \\
%- \frac{1}{2} \int_{\Gamma} (\ttA^{\alpha\beta\gamma\delta})_t \nh_{t,\alpha\beta} \nh_{t,\gamma\delta} dS \nonumber\\
& %+ \int_{\Gamma} \frac{1}{\sqrt{\det(g_0)}}\Big[\sqrt{\det(g_0)}(\ttA^{\alpha\beta\gamma\delta})_t
%\nh_{,\alpha\beta}\Big]_{,\gamma\delta}\nh_{tt} dS
\qquad + \int_{\Gamma}
\frac{1}{\sqrt{\det(g_0)}}\Big[\sqrt{\det(g_0)}(\ttA^{\alpha\beta\gamma\delta})_t\Big]_{,\gamma\delta}
\nh_{,\alpha\beta} \nh_{tt} dS \quad(\equiv K_7) + R_5
\label{L2H1vtid}
%& + \int_{\Gamma} \Big[L_1^{\alpha\beta\gamma} \tth_{,\alpha\beta\gamma}\Big]_t \nh_{tt} dS
%+ \int_{\Gamma} (L_2)_t \nh_{tt} dS. \nonumber
\end{align}
%The second term on the right-hand side verifies
%\begin{align}
%\Big|\int_{\Gamma} (\ttA^{\alpha\beta\gamma\delta})_t \nh_{t,\alpha\beta}\nh_{t,\gamma\delta} dS\Big|
%\le C(M) \|\th_t\|_{H^{2.5}(\Gamma)} \|\nabla_0^2 \nh_t\|^2_{L^2(\Gamma)}, \label{vtestimate1}
%\end{align}
%and the third term satisfies the equality
%\begin{align*}
%& \int_{\Gamma} \frac{1}{\sqrt{\det(g_0)}}\Big[\sqrt{\det(g_0)}(\ttA^{\alpha\beta\gamma\delta})_t
%\nh_{,\alpha\beta}\Big]_{,\gamma\delta} \nh_{tt} dS \\
%=& \int_{\Gamma} \frac{1}{\sqrt{\det(g_0)}}\Big[\sqrt{\det(g_0)}(\ttA^{\alpha\beta\gamma\delta})_t\Big]_{,\gamma\delta} \nh_{,\alpha\beta}
%\nh_{tt} dS \quad(\equiv K_7)\\
%&+ 2 \int_{\Gamma} \frac{1}{\sqrt{\det(g_0)}}\Big[\sqrt{\det(g_0)}(\ttA^{\alpha\beta\gamma\delta})_t\Big]_{,\delta}
%\nh_{,\alpha\beta\gamma} \nh_{tt} dS \quad(\equiv K_8)\\
%& + \int_{\Gamma} \ttA^{\alpha\beta\gamma\delta} \nh_{,\alpha\beta\gamma\delta} \nh_{tt} dS \quad(\equiv
%K_9).
%\end{align*}
%Also,
%$
%\|\nh_{tt}\|_{L^4(\Gamma)} \le C(M)\Big[\|\nv\|_{H^2(\Omega)} + \|\nv_t\|_{H^1(\Omega)}\Big]
%$
%and hence
%\begin{align*}
%|K_8| + |K_9|
%\le&\ C(M)C(\delta,\delta_1) \|\nabla_0^4\nh\|^2_{L^2(\Gamma)} + \delta \|\nv\|^2_{H^2(\Omega)} + \delta_1
%\|\nv_t\|^2_{H^1(\Omega)},
%\end{align*}
%and
where $R_5$ is bounded by
\begin{align*}
& C \Big[1+ \|\th_t\|^2_{H^{2.5}(\Gamma)}\Big] \Big[1 + \|\nabla_0^2
\nh_t\|^2_{L^2(\Gamma)}\Big]
+ \delta \Big[\|\nv\|^2_{H^2(\Omega)} + \|\nabla_0^2 \nv\|^2_{H^1(\Omega_1')}\Big] \\
& + \delta_1 \|\nv_t\|^2_{H^1(\Omega)}
\end{align*}
for some constant $C$ depending on $M$, $\delta$ and $\delta_1$.
Also, by the inequality $\|\nh_{tt}\|_{L^4(\Gamma)} \le
C(M)\Big[\|\nv\|_{H^2(\Omega)} + \|\nv_t\|_{H^1(\Omega)}\Big]$,
\begin{align*}
|K_7| \le&\
C\|[\sqrt{\det(g_0)}(\ttA^{\alpha\beta\gamma\delta})_t]_{,\gamma\delta}\|_{H^{-0.5}(\Gamma)}
\Big\|\frac{1}{\sqrt{\det(g_0)}} \nh_{,\alpha\beta} \nh_{tt}\Big\|_{H^{0.5}(\Gamma)} \\
\le&\
C(M)C(\delta,\delta_1)\|\th_t\|^2_{H^{2.5}(\Gamma)}\|\nabla_0^4
\nh\|^2_{L^2(\Gamma)} + \delta\|\nv\|^2_{H^2(\Omega)} +
\delta_1\|\nv_t\|^2_{H^1(\Omega)}.
\end{align*}
%Therefore,
%\begin{align}
%& \Big|\int_{\Gamma} \frac{1}{\sqrt{\det(g_0)}}\Big[\sqrt{\det(g_0)}(\ttA^{\alpha\beta\gamma\delta})_t
%\nh_{,\alpha\beta}\Big]_{,\gamma\delta} \nh_{tt} dS\Big| \label{vtestimate2}\\
%\le&\ C(M)C(\delta,\delta_1)(1+\|\th_t\|^2_{H^{2.5}(\Gamma)})\|\nabla_0^4 \nh\|^2_{L^2(\Gamma)}
%+ \delta\|\nv\|^2_{H^2(\Omega)} + \delta_1 \|\nv_t\|^2_{H^1(\Omega)}. \nonumber
%\end{align}

\begin{remark}\label{pqr}
The bound for $K_7$ can be refined even further as
$$|K_7|\le C(M)C(\delta)\|\th_t\|^2_{H^{1.5}(\Gamma)} \|\nabla_0^2 \nh\|^2_{H^{1.5}(\Gamma)}
+ \delta \|\nv\|^2_{H^3(\Omega)} + \delta
\|\nv_t\|^2_{H^1(\Omega)};$$ it is this inequality that will be used
in the proof of the fixed-point argument.
\end{remark}

%As for the remaining terms in (\ref{L2H1vtid}), for the fourth term we have
%\begin{align}
%&\ \int_{\Gamma} \Big[L_1^{\alpha\beta\gamma} \tth_{,\alpha\beta\gamma}\Big]_t \nh_{tt} dS \nonumber\\
%=&\ \int_{\Gamma} (L_1^{\alpha\beta\gamma})_t \tth_{,\alpha\beta\gamma} h_{tt} dS
%- \int_{\Gamma} \frac{1}{\sqrt{\det(g_0)}} \tth_{t,\alpha\beta} \Big[\sqrt{\det(g_0)}
%\nh_{tt} L_1^{\alpha\beta\gamma}\Big]_{,\gamma} dS \nonumber\\
%\le &\ C(M)C(\delta,\delta_1)\Big[1 + \|\th_t\|^2_{H^{2.5}(\Gamma)}\Big]
%+ \delta \|\nabla_0^2 \nv\|^2_{H^1(\Omega_1')} + \delta_1\|\nv_t\|^2_{H^1(\Omega)}\,, \label{vtestimate5}
%\end{align}
%while for the last term on the right-hand side in (\ref{L2H1vtid}),
%\begin{align}
%\int_{\Gamma} (L_2)_t \nh_{tt}dS \le C(M)C(\delta,\delta_1)
%+ \delta \|\nabla_0^2 \nv\|^2_{H^1(\Omega_1')} + \delta_1\|\nv_t\|^2_{H^1(\Omega)}.
%\label{vtestimate6}
%\end{align}
It remains to estimate $K_3$ to $K_6$. By proper use of H$\ddot{\text{o}}$lder's inequality,
%Similar to (\ref{vtestimate5}) and (\ref{vtestimate6}), we have that the lower order terms in $K_3$ to $K_5$, i.e,
%terms containing $L_1$ and $L_2$, can be bounded by
%\begin{align}
%C(M)\|\th_t\|_{H^{2.5}(\Gamma)}\|\nv\|_{H^{1.5}(\Gamma)}. \label{vtestimate8}
%\end{align}
%For the highest order term, we note that by (\ref{hconvergence}),
%\begin{align*}
%&\ \|(\sqrt{\det(g_0)}\ttA^{\alpha\beta\gamma\delta}\nh_{,\alpha\beta})_{t,\gamma\delta}\|_{H^{-1.5}(\Gamma)} \le
%\|\sqrt{\det(g_0)}(\ttA^{\alpha\beta\gamma\delta}\nh_{,\alpha\beta})_t\|_{H^{0.5}(\Gamma)} \\
%\le&\ C(M)\Big[t^{1/4} \|\nh_t\|_{H^{2.5}(\Gamma)} + \|\nabla_0^4 \nh\|_{L^2(\Gamma)}\Big].
%\end{align*}
%Therefore, combining with an upper bound (\ref{vtestimate8}) for the lower order terms, we find that
\begin{align*}
|K_3| + |K_5| + |K_6| %\le&\ C(M) \Big[t^{1/4} \|\nh_t\|_{H^{2.5}(\Gamma)} + \|\nabla_0^4 \nh\|_{L^2(\Gamma)}\Big]
%\Big[\|\nh_t\|_{H^{2.5}(\Gamma)} + \|\nv\|_{H^{1.5}(\Gamma)} \Big]\\
%& + C(M)\|\th_t\|_{H^{2.5}(\Gamma)} \|\nv\|_{H^{1.5}(\Gamma)} \\
\le&\ C\Big[1+\|\th_t\|^2_{H^{2.5}(\Gamma)}\Big]\Big[1 + \|\nabla_0^4 \nh\|^2_{L^2(\Gamma)} \Big] \\
& + (\delta + C t^{1/2})\|\nv\|^2_{H^3(\Omega)} + \delta
\|\nv_t\|^2_{H^1(\Omega)}
\end{align*}
for some constant $C$ depending on $M$ and $\delta$. %Also,
%\begin{align}
%|K_6|
%\le&\ C(M)C(\delta)\Big[\|\nabla_0^4 \nh\|^2_{L^2(\Gamma)} + 1\Big] + \delta \|\nv_t\|^2_{H^1(\Omega)}.
%\label{vtestimate7}
%\end{align}
For $K_4$, most of the terms can be estimated in the same fashion except the term
\begin{align*}
\int_{\Gamma}
\frac{1}{\sqrt{\det(g_0)}}\Big[\sqrt{\det(g_0)}\ttA^{\alpha\beta\gamma\delta}\nh_{t,\alpha\beta}\Big]
[(\nabla_0\tth_{t,\gamma\delta},0)\cdot(\nv\circ\tteta^{-\tau})] dS
\end{align*}
which is identical to
\begin{align*}
&\int_{\Gamma} \Big\{\frac{1}{\sqrt{\det(g_0)}}\Big[\sqrt{\det(g_0)}
\ttA^{\alpha\beta\gamma\delta}\nh_{t,\alpha\beta}\Big]
[(\nabla_0\tth_{,\gamma\delta},0)\cdot (\nv\circ\tteta^{-\tau})]
\Big\}_t dS \ (\equiv K_{8}) + R_6
%-& \int_{\Gamma} \frac{1}{\sqrt{\det(g_0)}}\Big[\sqrt{\det(g_0)}\ttA^{\alpha\beta\gamma\delta}\nh_{t,\alpha\beta}\Big]_t
%[(\nabla_0\tth_{,\gamma\delta},0)\cdot (\nv\circ\tteta^{-\tau})] dS \quad(\equiv K_{9})\\
%-& \int_{\Gamma} \frac{1}{\sqrt{\det(g_0)}}\Big[\sqrt{\det(g_0)}\ttA^{\alpha\beta\gamma\delta}\nh_{t,\alpha\beta}\Big]
%[(\nabla_0\tth_{,\gamma\delta},0) \cdot (\nv\circ\tteta^{-\tau})_t] dS .\quad(\equiv K_{10})
\end{align*}
where
\begin{align*}
|R_6| \le C\|\th\|^2_{H^{5.5}(\Gamma)} \Big[\|\nv\|^2_{L^2(\Omega)}
+ \|\nabla_0^2 \nh_t\|^2_{L^2(\Gamma)} \Big] + \delta
\|\nv\|^2_{H^3(\Omega)} + \delta_1 \|\nv_t\|^2_{H^1(\Omega)}
\end{align*}
for some constant $C$ depending on $M$, $\delta$ and $\delta_1$.
Time integrating $K_{8}$ and use the interpolation inequality %(\ref{interpolation1a}) (or (\ref{interpolation1b}) if $n=2$)
together with Young's inequality, we find that
\begin{align}
&\ \Big|\int_0^t K_{8}(s) ds\Big| \le
C(M)\Big[\|u_0\|^2_{H^{2.5}(\Omega)} + \|\nabla_0^2
\nh_t\|_{L^2(\Omega)}
\|\nv\|_{L^4(\Omega)}\Big] \nonumber\\
\le&\ C(M)C(\delta_1,\delta_2)N_3(u_0,F) + \delta_2
\|\nabla_0^2\nh_t\|^2_{L^2(\Gamma)} + \delta_1 \int_0^t
\|\nv_t\|^2_{H^1(\Omega)} ds \label{K10}
\end{align}
where
\begin{align*}
N_3(u_0,F) :=&\ \|u_0\|^2_{H^{2.5}(\Omega)} + \|u_0\|^2_{H^{4.5}(\Gamma)} + \|F\|^2_{L^2(0,T;H^1(\Omega))} \\
& + \|F_t\|^2_{L^2(0,T;H^1(\Omega)')} + \|F(0)\|^2_{H^1(\Omega)} + 1
\end{align*}
and we use $\displaystyle{\|\nv\|^2_{H^1(\Omega)} \le C\Big[\int_0^t
\|\nv_t\|^2_{H^1(\Omega)} ds +
\|u_0\|^2_{H^1(\Omega)}\Big]}$ to obtain (\ref{K10}). %The worst term in $K_{9}$ is
%\begin{align*}
%\int_{\Gamma} \sqrt{\det(g_0)}\ttA^{\alpha\beta\gamma\delta}\nh_{tt,\alpha\beta}
%[(\nabla_0\tth_{,\gamma\delta},0)\cdot (\nv\circ\tteta^{-\tau})] dS
%\end{align*}
%which, by $H^{-1.5}$-$H^{1.5}$ pairing, can be bounded by
%\begin{align*}
%&\ C(M) \|\nh_{tt}\|_{H^{0.5}(\Gamma)} \|\th\|_{H^{4.5}(\Gamma)} \|\nv\|_{H^{1.5}(\Gamma)} \\
%\le&\ C(M)C(\delta,\delta_1) \|\th\|^2_{H^{5.5}(\Gamma)} \|\nv\|^2_{L^2(\Omega)} + \delta
%\|\nv\|^2_{H^3(\Omega)} + \delta_1 \|\nv_t\|^2_{H^1(\Omega)}.
%\end{align*}
%Therefore,
%\begin{align*}
%|K_{9}| \le&\ C \Big[ \|\th\|^2_{H^{5.5}(\Gamma)}\|\nv\|^2_{L^2(\Omega)} + \|\nabla_0^2
%\nh_t\|^2_{L^2(\Gamma)}\Big] + \delta \|\nv\|^2_{H^3(\Omega)} + \delta_1 \|\nv_t\|^2_{H^1(\Omega)}
%\end{align*}
%for some constant $C$ depending on $M$, $\delta$ and $\delta_1$. Also,
%\begin{align*}
%|K_{12}| \le&\ C(M)C(\delta)\|\th\|^2_{H^{4.5}(\Gamma)}\|\nabla_0^2 \nh_t\|^2_{L^2(\Gamma)} + \delta
%\|\nv_t\|^2_{H^1(\Omega)}
%\end{align*}
and hence
\begin{align}
\sum_{i=3}^6 |K_i| \le&\ C \Big[1+ \|\th\|^2_{H^{5.5}(\Gamma)} +
\|\th_t\|^2_{H^{2.5}(\Gamma)}\Big]
\Big[1+\|\nv\|^2_{L^2(\Omega)} + \|\nabla_0^4 \nh\|^2_{L^2(\Gamma)}\Big] \nonumber\\
& + (\delta + C t^{1/2})\|\nv\|^2_{H^3(\Omega)} + \delta_1
\|\nv_t\|^2_{H^1(\Omega)} + K_{8} \label{vtestimate9}
\end{align}
with $K_{8}$ satisfying inequality (\ref{K10}).
Finally, combining all the estimates,
\begin{align}
& \int_0^t \|\nabla_0^2 \nh_t\|^2_{L^2(\Gamma)} ds \le \int_0^t
\int_{\Gamma} \Big[[L_\tth(\nh)
(\nabla_0\tth,-1)]\circ\tteta^\tau\Big]_t \cdot \nv_t dS
+ C N_3(u_0,F) \nonumber\\
& + C \int_0^t K(s) \Big[\|\nv\|^2_{L^2(\Omega)} + \|\nabla_0^4
\nh\|^2_{L^2(\Gamma)}
+ \|\nabla_0^2 \nh_t\|^2_{L^2(\Gamma)} \Big]ds \label{mainvtestimate2}\\
& + (\delta + C t^{1/2})\int_0^t \|\nv\|^2_{H^3(\Omega)} ds +
(\delta_1 + C t^{1/2}) \int_0^t \|\nv_t\|^2_{H^1(\Omega)} ds +
\delta_2 \|\nabla_0^4 \nh\|^2_{L^2(\Gamma)} \nonumber
\end{align}
for some constant $C$ depending on $M$, $\delta$, $\delta_1$ and $\delta_2$.

\section*{Acknowledgments}
AC, DC, and SS were supported by the National Science Foundation
under grant NSF ITR-0313370.


\begin{thebibliography}{50}

\bibitem{AuBe} F.~Auricchio, L. Beir${\tilde{\rm{a}}}$o da Veiga
and C. Lovadina, {\scshape Remarks on the asymptotic behaviour of Koiter
shells,} Computers and Structures, {\bf 80}
(2002), 735-745

\bibitem{DaVe} H. Beir\~{a}o. da Veiga, {\scshape On the existence of strong
solutions to a coupled fluid-structure evolution problem,}
J. Math. Fluid Mech.,  {\bf 6}  2004,  21--52.

\bibitem{Bo2005} M.~Boulakia, {\scshape Existence of weak solutions for an interaction
problem between an elastic structure and a compressible viscous
fluid,} J. Math. Pures Appl. (9) {\bf 84} (2005),  no. 11,
1515--1554.

\bibitem{ChDeEsGr} A.~Chambolle, B.~Desjardins, M.J.~Esteban, C.~Grandmont,
 {\scshape  Existence of weak solutions for an unsteady fluid-plate
interaction problem,} Preprint.

\bibitem{RSCh2002} R.S.~Chadwick, {\scshape Axisymmetric indentation of a
thin incompressible
elastic layer,} SIAM J. Appl. Math. {\bf 62} (2002) 1520--1530.

\bibitem{ChCoSh2006} C.H.~Cheng, D.~Coutand and S. Shkoller,
{\scshape Navier-Stokes equations interacting with a nonlinear
elastic solid shell,} preprint

\bibitem{PGC1} P.G.~Ciarlet, {\scshape Introduction to linear shell theory,}
Series in Applied
Mathematics (Paris), vol. 1, Gauthier-Villars, Editions
Scientifiques et M edicales Elsevier, Paris, 1998.

\bibitem{PCG2} P.G.~Ciarlet, {\scshape Mathematical elasticity, vol. III,}
Studies in
Mathematics and its Applications, vol. {\bf 29}, North-Holland,
Amsterdam, 2000, Theory of shells.

\bibitem{CoSh2002} D.~Coutand and S. Shkoller, {\scshape Unique solvability
of the free-boundary
Navier-Stokes equations with surface tension,}

\bibitem{CoSh2005} D.~Coutand and S. Shkoller, {\scshape On the motion of an
elastic solid
inside of an incompressible viscous fluid,} to appear in Arch. Rational
Mech. Anal.

\bibitem{CoSh2006} D.~Coutand and S. Shkoller, {\scshape On the interaction
between quasilinear
elastodynamics and the Navier-Stokes equations,}

\bibitem{DeEs} B.~Desjardins, M.J.~Esteban,
 {\scshape  Existence of weak solutions for the motion of rigid bodies in a
viscous fluid,} Arch. Rational Mech. Anal.,
{\bf 146} (1999), 59--71.

\bibitem{DeEsGrTa} B.~Desjardins, M.J.~Esteban, C.~Grandmont, P.~Le Tallec,
 {\scshape Weak solutions for a fluid-structure interaction problem,} Rev.
Mat. Complut., {\bf 14} (2001), 523--538.

\bibitem{LCE1} L.C.~Evans, {\scshape Partial Differential Equations},
Graduate Studies in Mathematics, {\bf 19} American Mathematical
Society, Providence, RI, 1998.

\bibitem{YCFu1981} Y.C.~Fung, {\scshape Biomechanics: Mechanical Properties
of Living
Tissues,} Springer, New York, 1981.

\bibitem{GPG1} G.P.~Galdi, {\scshape An Introduction to the
Mathematical Theory of the Navier-Stokes Equations Volume I,}
Springer Tracts in Natural Philosophy, Vol {\bf 38}.

\bibitem{GeKrMa1996} Z.~Ge, H.P. Kruse and J.E. Marsden, {\scshape The
limits of Hamiltonian
structures in three-dimensional elasticity, shells, and rods,}
J. Nonlinear Sci. Vol. {\bf 6} (1996), 19-57.

\bibitem{WMGe1995} W.M.~Gelbart, {\scshape Micelles, Membranes, Microemulsions, and
Monolayers,} Springer-Verlag New York, 1995.

\bibitem{EG} E.~Givelberg, {\scshape Modeling Elastic Shells Immersed in
Fluid,} Comm. Pure Appl.
Math., {\bf 57} (2004), no. 3, 283-309.

\bibitem{GrMa} C.~Grandmont, Y.~Maday,
 {\scshape  Existence for unsteady fluid-structure interaction problem,}
Math. Model. Numer. Anal., {\bf 34} (2000), 609--636.

\bibitem{LePeLa} R.J.~Leveque, C.S. Peskin and P.D. Lax, {\scshape Solution
of
a two-dimensional cochlea model with fluid viscosity,} SIAM J. Appl.
Math., {\bf 45} (1985), no. 3, 450-464.

\bibitem{LiWa} C. Liu, N.J. Walkington, {\scshape An Eulerian description of
fluids containing visco-elastic particles}, Arch. Rational Mech. Anal., {\bf
159} (2001), 229-252.

\bibitem{MiSeWoDo1994} L.~Miao, U.~Seifert, M.~Wortis, H.~Dobereiner,
{\scshape Budding
transitions of fluid-bilayer vesicle: the effect of area-difference
elasticity,} Phys. Rev. E {\bf 49} (1994) 5389--5407.

\bibitem{OuHe1989} Z.~Ou-Yabg, W.~Helfrich, {\scshape Bending energy of
vesicle membranes:
general expressions for the first, second and third variation of the
shape energy and applications to spheres and cylinders,} Phys. Rev.
E {\bf 39} (1989) 5280--5288.

\bibitem{Serre} D.~Serre,
{\scshape  Chute libre d'un solide dans un fluide visqueux
incompressible: Existence,} Japan J. Appl. Math., {\bf 4} (1987),
33--73.

\bibitem{USe1991} U.~Seifert, {\scshape Adhesion of vesicles in two
dimensions,} Phys. Rev. A
{\bf 43} (1991) 6803--6814.

\bibitem{USe1993} U.~Seifert, {\scshape Curvature-induced lateral phase
segregation in
two-component vesicles,} Phys. Rev. Lett. {\bf 70} (1993) 1335-1338.

\bibitem{Wein1972} H.F.~Weinberger, {\scshape Variational properties of
steady fall in Stokes
flow,} J. Fluid Mech., {\bf 52} (1972), 321--344.

\end{thebibliography}
\end{document}